\documentclass[12pt]{article}

\usepackage{epsf,epsfig,amsmath,amsfonts}
\usepackage{graphicx}

\newtheorem{thm}{Theorem}
\newtheorem{exa}{Example}
\newtheorem{cor}[thm]{Corollary}

\newtheorem{lem}[thm]{Lemma}
\newtheorem{defn}{Definition}
\newtheorem{rem}{Remark}
\newtheorem{assm}{Assumption}

\newcommand{\be}{\begin{equation*}}
\newcommand{\ben}{\begin{equation}}
\newcommand{\ee}{\end{equation*}}
\newcommand{\een}{\end{equation}}

\newcommand{\A}{{\mathcal A}}
\newcommand{\F}{{\mathcal F}}
\newcommand{\G}{{\mathcal G}}
\newcommand{\T}{{\mathcal T}}
\newcommand{\Lop}{{\mathcal L}}
\newcommand{\Hc}{{\mathcal H}}
\newcommand{\Om}{\Omega}
\newcommand{\om}{\omega}
\newcommand{\R}{{\mathbb R}}

\newcommand{\intI}{\operatorname{int} {\mathcal I}}

\newcommand{\half}{{\frac{1}{2}}}
\newcommand{\CI}{{\mathcal I}}
\newcommand{\ac}{\mathrm{ac}}
\newcommand{\s}{\mathrm{s}}
\newcommand{\e}{{\mathbb{E}}}
\newcommand{\p}{{\mathbb{P}}}
\newcommand{\str}{{\mathbb{S}}}

\begin{document}

\title{\sc On the Optimal Stopping of a One-dimensional
Diffusion\footnote{Research supported
by EPSRC grant no.\,GR/S22998/01 and the Isaac Newton Institute,
Cambridge}}

\author{
{\sc Damien Lamberton\footnote{Laboratoire d'Analyse
et de Math\'{e}matiques Appliqu\'{e}es, Universit\'{e}
Paris-Est, 5 Boulevard Descartes, 77454
Marne-la-Vall\'{e}e Cedex 2, France,
\texttt{damien.lamberton@univ-mlv.fr}}
\ and Mihail Zervos\footnote{Department of Mathematics,
London School of Economics, Houghton Street, London
WC2A 2AE, UK, \texttt{m.zervos@lse.ac.uk}}}}

\maketitle

\begin{abstract}
We consider the one-dimensional diffusion $X$ that satisfies
the stochastic differential equation
\ben
dX_t = b(X_t) \, dt + \sigma (X_t) \, dW_t \label{SDE}
\een
in the interior $\intI = \mbox{} ]\alpha, \beta[$ of a given interval
$\CI \subseteq [-\infty, \infty]$, where $b, \sigma: \intI \rightarrow
\R$ are Borel-measurable functions and $W$ is a standard
one-dimensional Brownian motion.
We allow for the endpoints $\alpha$ and $\beta$ to be inaccessible
or absorbing.
Given a Borel-measurable function $r: \CI \rightarrow \R_+$
that is uniformly bounded away from 0, we establish a new analytic
representation of the $r(\cdot)$-potential of a continuous
additive functional of $X$.
Furthermore, we derive a complete characterisation of differences
of two convex functions in terms of appropriate $r(\cdot)$-potentials,
and we show that a function $F: \CI \rightarrow \R_+$
is $r(\cdot)$-excessive if and only if it is the difference of two
convex functions and $- \bigl(\half \sigma ^2 F'' + bF' - rF \bigr)$
is a positive measure.
We use these results to study the optimal stopping problem that
aims at maximising the performance index
\ben
\e_x \left[ \exp \left( - \int _0^\tau r(X_t) \, dt \right) f(X_\tau)
{\bf 1} _{\{ \tau < \infty \}} \right] \label{PI}
\een
over all stopping times $\tau$, where $f: \CI \rightarrow \R_+$
is a Borel-measurable function that may be unbounded.
We derive a simple necessary and sufficient condition for the value
function $v$ of this problem to be real-valued.
In the presence of this condition, we show that $v$ is the
difference of two convex functions, and we prove that it
satisfies the variational inequality
\ben
\max \left\{ \half \sigma ^2 v'' + bv' - rv , \ \overline{f} - v \right\} = 0
\label{QVIa}
\een
in the sense of distributions, where $\overline{f}$ identifies with
the upper semicontinuous envelope of $f$ in the interior $\intI$
of $\CI$.
Conversely, we derive a simple necessary and sufficient
condition for a solution to (\ref{QVIa}) to identify with the value
function $v$.
Furthermore, we establish several other characterisations of
the solution to the optimal stopping problem, including a
generalisation of the so-called ``principle of smooth fit''.
In our analysis, we also make a construction that is concerned
with pasting weak solutions to (\ref{SDE}) at appropriate hitting
times, which is an issue of fundamental importance to dynamic
programming.
\end{abstract}

\section{Introduction}

We consider the one-dimensional diffusion $X$ that
satisfies the SDE (\ref{SDE}) in the interior $\intI = \mbox{}
]\alpha, \beta[$ of a given interval $\CI \subseteq [-\infty,
\infty]$.
We assume that $b, \sigma : \intI \rightarrow \R$ are
Borel-measurable functions satisfying appropriate local
integrability and non-degeneracy conditions ensuring that
(\ref{SDE}) has a weak solution that is unique in the sense
of probability law up to a possible explosion time at which
$X$ hits the boundary $\{ \alpha, \beta \}$ of $\CI$
(see Assumption~\ref{A1} in Section~\ref{sec:pr-form}).
If the boundary point $\alpha$ (resp., $\beta$) is inaccessible,
then the interval $\CI$ is open from the left (resp., open from
the right), while, if $\alpha$ (resp., $\beta$) is not inaccessible,
then it is absorbing and the interval $\CI$ is closed from the
left (resp., closed from the right).

In the presence of Assumption~\ref{A1}, a weak solution to
(\ref{SDE}) can be obtained by first time-changing a standard
one-dimensional Brownian motion and then making an
appropriate state space transformation.
This construction can be used to prove all of the results
that we obtain by first establishing them assuming that the
diffusion $X$ identifies with a standard one-dimensional
Brownian motion.
However, such an approach would hardly simplify the formalism
because the data $b$ (resp., $\sigma$) appear in all of the analysis
exclusively (resp., mostly) though the operators $\Lop$,
$\Lop_\ac$ defined by (\ref{lop})--(\ref{lopac}) below.
Furthermore, deriving the general results, which are important
because many applications assume specific functional forms
for the data $b$ and $\sigma$, by means of this approach
would require several time changes and state space transformations,
which would lengthen the paper significantly.

Given a point $z \in \intI$, we denote by $L^z$ the right-sided
local time process of $X$ at level $z$ (see Revuz and
Yor~\cite[Section~VI.1]{RY} for the precise definition of $L^z$
and its properties).
Also, we denote by ${\mathcal B} ({\mathcal J})$ the Borel
$\sigma$-algebra on any given interval ${\mathcal J}
\subseteq [-\infty, \infty]$.
With each signed Radon measure $\mu$ on $\bigl( \intI ,
{\mathcal B} (\intI) \bigr)$ such that $\sigma^{-2}$ is locally
integrable with respect to $|\mu|$, we associate the continuous
additive functional
\ben
\A_t^\mu = \int _\alpha^\beta \frac{L_t^z}{\sigma ^2 (z)} \,
\mu (dz) , \quad t \in [0, T_\alpha \wedge T_\beta[ ,
\label{A}
\een
where $T_\alpha$ (resp., $T_\beta$) is the first hitting time of
$\alpha$ (resp., $\beta$).
It is worth noting that (\ref{A}) provides a one-to-one
correspondence between the continuous additive functionals
of the Markov process $X$ and the signed Radon measures
on $\bigl( \intI , {\mathcal B} (\intI) \bigr)$ (see Theorem~X.2.9,
Corollary~X.2.10 and the comments on Section~2 at the end
of Chapter~X in Revuz and Yor~\cite[Section~X.2]{RY}).
We also consider a discounting rate function $r : \CI \rightarrow
\R_+$, we assume that this is a Borel-measurable
function that is uniformly bounded away from 0 and satisfies
a suitable local integrability condition (see Assumption~\ref{A2}
in Section~\ref{sec:pr-form}), and we define
\ben
\Lambda _t \equiv \Lambda _t (X) = \int _0^t r(X_s) \, ds .
\label{Lambda}
\een

Given a signed Radon measure $\mu$ on $\bigl( \intI ,
{\mathcal B} (\intI) \bigr)$, we consider the $r(\cdot)$-potential
of the continuous additive functional $A^\mu$, which is
defined by
\ben
R_\mu (x) = \e_x \left[ \int _0^{T_\alpha \wedge T_\beta}
e^{-\Lambda_t} \, dA_t^{\mu} \right] . \label{R-prob-intro}
\een
We recall that a function $F: \intI \rightarrow \R$ is the
difference of two convex functions if and only if its
left-hand side derivative $F_-'$ exists and its second
distributional derivative is a measure, and we define
the measure $\Lop F$ by
\be
\Lop F (dx) = \half \sigma^2 (x) F'' (dx) + b(x) F_-' (x) \,
dx - r(x) F(x) \, dx .
\ee
In the presence of a general integrability condition ensuring
that the potential $R_\mu$ is well-defined, we show that
it is the difference of two convex functions, the measures
$\Lop R_\mu$ and $-\mu$ are equal, and
\begin{align}
R_\mu (x) & = \frac{2}{C} \varphi (x) \int _{]\alpha, x[}
\frac{\psi (s)} {\sigma^2 (s) p' (s)} \, \mu (ds) + \frac{2}{C}
\psi (x) \int _{[x, \beta[} \frac{\varphi (s)} {\sigma^2 (s) p' (s)}
\, \mu (ds) \nonumber \\
& = \int _{]\alpha, \beta[} \frac{2\varphi (x) \psi (x)}{C\sigma^2
(s) p'(s)} \min \left\{ \frac{\psi (s)}{\psi (x)} , \frac{\varphi (s)}
{\varphi (x)} \right\} \mu (ds) , \label{R-ana-intro}
\end{align}
where $C>0$ is an appropriate constant, $p: \intI \rightarrow
\R$ is the scale function of $X$, and $\varphi , \psi : \intI
\rightarrow \mbox{} ]0,\infty[$ are $C^1$ functions with absolutely
continuous with respect to the Lebesgue measure derivatives
spanning the solution space of the ODE
\be
\half \sigma^2 (x) g'' (x) + b(x) g' (x) - r(x) g(x) = 0 ,
\ee
and such that $\varphi$ (resp., $\psi$) is decreasing
(resp., increasing) (see Theorem~\ref{prop:ODE-special}).
If the signed measure $\mu^h$ is absolutely continuous with
respect to the Lebesgue measure with Radon-Nikodym
derivative given by a function $h$, then the potential
$R_{\mu^h}$ admits the expressions
\begin{align}
R_{\mu^h} (x) & = \e_x \left[ \int _0^{T_\alpha \wedge
T_\beta} e^{-\Lambda _t} h(X_t) \, dt \right] \nonumber \\
& = \frac{2}{C} \varphi (x) \int _\alpha^x \frac{\psi (s)}
{\sigma^2 (s) p' (s)} h(s) \, ds + \frac{2}{C} \psi (x) \int
_x^\beta \frac{\varphi (s)} {\sigma^2 (s) p' (s)} h(s) \, ds
\label{RAC-inftro}
\end{align}
(see Corollary~\ref{cor:R_mu_h} for this and other
related results).
Conversely, we show that, under a general growth condition,
a difference of two convex functions $F: \intI \rightarrow
\R$ is such that (a) both limits $\lim _{y \downarrow \alpha}
F(y) / \varphi (y)$ and $\lim _{y \uparrow \beta} F(y)
/ \psi (y)$ exist, (b) $F$ admits the characterisation
\ben
F(x) = \lim _{y \downarrow \alpha} \frac{F(y)}{\varphi (y)}
\varphi (x) + R_{-\Lop F} (x) + \lim _{y \uparrow \beta}
\frac{F(y)}{\psi (y)} \psi (x) , \label{F-repr-intro}
\een
and (c) an appropriate form of Dynkin's formula holds
true (see Theorem~\ref{cor:F-repres}).
With a view to optimal stopping, we use these results to show
that a function $F: \CI \rightarrow \R_+$ is $r(\cdot)$-excessive
if and only if it is the difference of two convex functions and
$- \Lop F$ is a positive measure
(see Theorem~\ref{prop:excessive} for the precise result).

If $r$ is constant, then general theory of Markov processes
implies the existence of a transition kernel $u^r$ such that
$R_\mu (x) = \int _{]\alpha, \beta[} u^r (x,s) \, \mu (ds)$
(see Meyer~\cite{M} and Revuz~\cite{R}).
If $X$ is a standard Brownian motion, then
\be
u^r (x,s) = \frac{1}{\sqrt{2r}} e^{-\sqrt{2r} |x-s|}
\ee
(see Revuz and Yor~\cite[Theorem~X.2.8]{RY}).
The general expression for this kernel provided by
(\ref{R-ana-intro}) is one of the contributions of this paper.
On the other hand, the identity in (\ref{RAC-inftro}) is well-known
and can be found in several references (e.g., see
Borodin and Salminen~\cite[II.4.24]{BS}).
Also, Johnson and Zervos~\cite{JZ} prove that the potential
given by (\ref{R-prob-intro}) admits the analytic expression
(\ref{R-ana-intro}) and show that the measures $\Lop R_\mu$
and $-\mu$ are equal when both of the endpoints $\alpha$
and $\beta$ are assumed to be inaccessible.

The representation of differences of two convex functions
given by (\ref{F-repr-intro}) is also new.
Such a result is important for the solution to one-dimensional
infinite time horizon stochastic control as well as optimal
stopping problems using dynamic programming.
Indeed, the analysis of several explicitly solvable problems
involve such a representation among their assumptions.
For constant $r$, Salminen~\cite{Sa} considered more general
one-dimensional linear diffusions than the one given by
(\ref{SDE}) and used Martin boundary theory to show that
every $r$-excessive function admits a representation that is
similar to but much less straightforward than the one in
(\ref{F-repr-intro}).
Since a function on an open interval is the difference of
two convex functions if and only if it is the difference of two
excessive functions (see \c{C}inlar, Jacod, Protter and
Sharpe~\cite{CJPS}), the representation derived by
Salminen~\cite{Sa} can be extended to differences of two
convex functions.
However, it is not straightforward to derive such an extension
of the representation in Salminen~\cite{Sa} from (\ref{F-repr-intro})
or vice-versa when the underlying diffusion satisfies (\ref{SDE})
and $r$ is constant.

The result that a function $F$ is $r(\cdot)$-excessive
if and only if it is the difference of two convex functions and
$-\Lop F$ is a positive measure is perhaps the simplest
possible characterisation of excessive functions because
it involves only derivative operators.
In fact, we show that this result is equivalent to the
characterisations of excessive functions derived by
Dynkin~\cite{Dy65} and Dayanik~\cite{Da}
(see Corollary~\ref{cor:excessive}).

We use the results that we have discussed above to
analyse the optimal stopping problem that aims at
maximising the performance criterion given by (\ref{PI})
over all stopping times $\tau$, assuming that the reward
function $f$ is a positive Borel-measurable function
that may be unbounded (see Assumption~\ref{A2} in
Section~\ref{sec:pr-form}).
We first prove that the value function $v$ is the difference
of two convex functions and satisfies the variational inequality
(\ref{QVIa}) in the sense of distributions, where $\overline{f}$
is defined by
\ben
\overline{f} (x) = \begin{cases} \limsup _{y \rightarrow x}
f(y) , & \text{if } x \in \intI , \\ f(\alpha) , & \text{if } \alpha
\text{ is absorbing and } x=\alpha , \\ f(\beta) , & \text{if }
\beta \text{ is absorbing and } x=\beta \end{cases}
\label{f-over}
\een
(see Definition~\ref{sense} and Theorem~\ref{thm:main1}.(I)--(II)
in Section~\ref{sec:solution}).
This result provides simple criteria for deciding which parts
of the interval $\CI$ must be subsets of the so-called
waiting region.
Indeed, the derived regularity of $v$ implies that all points at
which the reward function $f$ is discontinuous as well as
all ``minimal'' intervals in which $f$ cannot be expressed as
the difference of two convex functions (e.g., intervals
throughout which $f$ has the regularity of a Brownian
sample path) should be parts of the closure of the waiting
region.
Similarly, the support of the measure $(\Lop f)^+$ in all
intervals in which $\Lop f$ is well-defined should also
be a subset of the closure of the waiting region.

We then establish a verification theorem that is the
strongest one possible because it involves only the
optimal stopping problem's data.
In particular, we derive a simple necessary and sufficient
condition for a solution $w$ to (\ref{QVIa}) in the sense
of distributions to identify with the problem's value function
(see Theorem~\ref{thm:main2}.(I)--(II)).

These results establish a complete characterisation of the
value function $v$ in terms of the variational inequality
(\ref{QVIa}).
Indeed, they imply that the restriction of the optimal
stopping problem's value function $v$ in $\intI$
identifies with the unique solution to the variational inequality
(\ref{QVIa}) in the sense of Definition~\ref{sense} that
satisfies the boundary conditions
\be
\lim _{y \in \intI , \, y \downarrow \alpha} \frac{v(y)}{\varphi (y)}
= \limsup _{y \downarrow \alpha} \frac{f(y)}{\varphi (y)} 
\quad \text{and} \quad
\lim _{y \in \intI , \, y \uparrow \beta} \frac{v(y)}{\psi (y)}
= \limsup _{y \uparrow \beta} \frac{f(y)}{\psi (y)} .
\ee
It is worth noting that, if $\alpha$ (resp., $\beta$) is absorbing,
then the corresponding boundary condition is equivalent to
\be
\lim _{y \in \intI , \, y \downarrow \alpha} v(y) = \limsup
_{y \downarrow \alpha} f(y) \qquad \left( \text{resp., }
\lim _{y \in \intI , \, y \uparrow \beta} \psi (y) = \limsup
_{y \uparrow \beta} f(y) \right)
\ee
(see (\ref{phi-psi-aabs})--(\ref{phi-psi-babs})).
Also, it is worth stressing the precise nature of these boundary
conditions.
The limits on the left-hand sides are taken from inside
the interior $\intI$ of $\CI$ and they indeed exist.
On the other hand, the limsups on the right-hand sides
are taken from inside $\CI$ itself. 
Therefore, if, e.g., $\alpha$ is absorbing, then we are faced
either with
\be
v(\alpha) = f(\alpha) = \lim _{y \in \intI , \, y \downarrow \alpha}
v(y) = \limsup _{y \downarrow \alpha}  f(y) , \quad \text{if }
f(\alpha) = \limsup _{y \downarrow \alpha}  f(y) \geq \limsup
_{y \in \intI , \, y \downarrow \alpha} f(y) ,
\ee
or with
\be
v(\alpha) = f(\alpha) < \lim _{y \in \intI , \, y \downarrow \alpha}
v(y) = \limsup _{y \downarrow \alpha} f(y) , \quad \text{if }
f(\alpha) < \limsup _{y \downarrow \alpha}  f(y) = \limsup
_{y \in \intI , \, y \downarrow \alpha} f(y) .
\ee

Furthermore, we prove that
\ben
v(x) = \inf \bigl\{ A \varphi (x) + B \psi (x) \mid \ A, B \geq 0
\text{ and } A \varphi + B \psi \geq f \bigr\} \label{vLPintro}
\een
for all $x \in \intI$ (see Theorem~\ref{thm:main2}.(III)).
In fact, this characterisation can be used as a
verification theorem as well (see also the discussion further
below).

In the generality that we consider, an optimal stopping time
might not exist (see Examples~\ref{ex1}--\ref{ex3} in
Section~\ref{sec:ex}).
Moreover, the hitting time of the so-called ``stopping region'',
which is given by
\ben
\tau^\star = \inf \bigl\{ t \geq 0 \mid \ v(X_t) = f(X_t)
\bigr\} , \label{tau*intro}
\een
may not be optimal (see Examples~\ref{ex2}
and~\ref{ex3}).
In particular, Example~\ref{ex2} shows that $\tau^*$ may not
be optimal and that an optimal stopping time may not exist
at all unless $f$ satisfies appropriate boundary / growth
conditions.
Also, Example~\ref{ex3} reveals that $\tau^\star$ is not in
general optimal if $f \neq \overline{f}$.
In Theorem~\ref{thm:main1}.(III), we obtain a simple
sequence of $\varepsilon$-optimal stopping times if
$f$ is assumed to be upper semicontinuous, and we show
that $\tau^\star$ is an optimal stopping time if $f$ satisfies
an appropriate growth condition.

Building on the general theory, we also consider a number
of related results and characterisations.
In particular, we obtain a generalisation of the so-called
``principle of smooth fit'' (see part~(III) of Corollaries
\ref{call}, \ref{put} and \ref{straddle} in Section~\ref{sec:sf}).

In view of the version of Dynkin's formula (\ref{Rmuh-Dynkin})
in Corollary~\ref{cor:R_mu_h}, we can see that, if $h$ is
any function such that $R_{\mu^h}$ given by
(\ref{RAC-inftro}) is well-defined, then
\begin{align}
\sup _\tau \e_x & \left[ \int _0^{\tau \wedge T_\alpha \wedge
T_\beta} e^{-\Lambda_t} h(X_t) \, dt + e^{-\Lambda_{\tau
\wedge T_\alpha \wedge T_\beta}} f(X_\tau) {\bf 1} _{\{ \tau
< \infty \}} \right] \nonumber \\
& \mbox{} = R_{\mu^h} (x) + \sup _\tau \e_x \left[
e^{-\Lambda _{\tau \wedge T_\alpha \wedge T_\beta}}
\left( f - R_{\mu^h} \right) (X_{\tau \wedge T_\alpha \wedge
T_\beta}) {\bf 1} _{\{ \tau < \infty \}} \right] \nonumber \\
& \mbox{} = R_{\mu^h} (x) + \sup _\tau \e_x \left[
e^{-\Lambda _{\tau \wedge T_\alpha \wedge T_\beta}}
\left( f - R_{\mu^h} \right) ^+ (X_{\tau \wedge T_\alpha \wedge
T_\beta}) {\bf 1} _{\{ \tau < \infty \}} \right] . \label{OS-gen}
\end{align}
Therefore, all of the results on the optimal stopping problem
that we consider generalise most trivially to account for
the apparently more general optimal stopping problem
associated with (\ref{OS-gen}).

The various aspects of the optimal stopping theory have been
developed in several monographs, including
Shiryayev~\cite{Sh},
Friedman~\cite[Chapter~16]{F},
Krylov~\cite{K},
Bensoussan and Lions~\cite{BL},
El~Karoui~\cite{ElK},
{\O}ksendal~\cite[Chapter 10]{O}
and Peskir and Shiryaev~\cite{PS}.
In particular, the solution of optimal stopping problems using
classical solutions to variational inequalities has been extensively
studied (e.g., see Friedman~\cite[Chapter~16]{F},
Krylov~\cite{K} and Bensoussan and Lions~\cite{BL}).
Results in this direction typically make strong regularity
assumptions on the problem data (e.g., the diffusion
coefficients are assumed to be Lipschitz continuous).
To relax such assumptions, {\O}ksendal and Reikvam~\cite{OR}
and Bassan and Ceci~\cite{BC} have considered
viscosity solutions to the variational inequalities associated
with the optimal stopping problems that they study.
Closer to the spirit of this paper, Lamberton~\cite{L}
proved that the value function of the finite version of the
problem we consider here satisfies its associated
variational inequality in the sense of distributions.

Relative to the optimal stopping problem that we
consider here when $r$ is constant, Dynkin~\cite{Dy63} and
Shiryaev~\cite[Theorem~3.3.1]{Sh} prove that the
value function $v$ identifies with the smallest $r$-excessive
function that majorises the reward function $f$ if $f$ is
assumed to be lower semicontinuous.
Also, Shiryaev~\cite[Theorem~3.3.3]{Sh} proves that
the stopping time $\tau^\star$ defined by (\ref{tau*intro})
is optimal if $f$ is assumed to be continuous and
bounded, while Salminen~\cite{Sa} establishes the optimality
of $\tau^*$ assuming that the smallest $r$-excessive
majorant of $f$ exists and $f$ is upper semicontinuous.
Later, Dayanik and Karatzas~\cite{DK} and Dayanik~\cite{Da},
who also considers random discounting instead of
discounting at a constant rate $r$, addressed the solution
of the optimal stopping problem by means of a certain
concave characterisation of excessive functions.
In particular, they established a generalisation of
the so-called ``principle of smooth fit'' that is similar to,
though not the same as, the one we derive here.

There are numerous special cases of the general
optimal stopping problem we consider that have been
explicitly solved in the literature.
Such special cases have been motivated by applications
or have been developed as illustrations of various general
techniques.
In all cases, their analysis relies on some sort of a
verification theorem.
Existing verification theorems for solutions using
dynamic programming and variational inequalities
typically make strong assumptions that are either
tailor-made or difficult to verify in practice.
For instance, Theorem~10.4.1 in {\O}ksendal~\cite{O}
involves Lipschitz as well as uniform integrability assumptions,
while, Theorem~I.2.4 in Peskir and Shiryaev~\cite{PS}
assumes the existence of an optimal stopping time,
for which, a sufficient condition is provided by Theorem~I.2.7.
Alternatively, they assume that the so-called stopping region
is a set of a simple specific form (e.g., see R\"{u}schendorf
and Urusov~\cite{RU} or Gapeev and Lerche~\cite{GL}).

Using martingale and change of measure techniques,
Beibel and Lerche~\cite{BL1, BL2}, Lerche and Urusov~\cite{LU}
and Christensen and Irle~\cite{CI} developed an approach
to determining an optimal stopping strategy at any given point
in the interval $\CI$.
Similar techniques have also been extensively used by
Alvarez~\cite{A1, A2, A3}, Lempa~\cite{Lempa} and references
therein.
To fix ideas, we consider the following representative cases
that can be associated with any given initial condition
$x \in \CI$.
If there exists a point $d_1 > x$ such that
\ben
C_1 := \sup _{x \in \CI} \frac{f(x)}{\psi (x)} = \frac{f(d_1)}
{\psi (d_1)} , \label{BL1intro}
\een
then $v(x) = C_1 \psi (x)$ and the first hitting time of
$\{ d_1 \}$ is optimal.
Alternatively, if there exist points $\kappa \in \mbox{} ]0,1[$
and $c_2 < x < d_2$ such that
\ben
C_2 := \sup _{x \in \CI} \frac{f(x)}{\kappa \psi (x) + (1-\kappa)
\varphi (x)} = \frac{f(c_2)}{\kappa \psi (c_2) + (1-\kappa)
\varphi (c_2)} = \frac{f(d_2)}{\kappa \psi (d_2) + (1-\kappa)
\varphi (d_2)} , \label{BL2intro}
\een
then $v(x) = \kappa C_2 \psi (x) + (1-\kappa) C_2 \varphi (x)$
and the first hitting time of $\{ c_2, d_2 \}$ is optimal.
On the other hand, if $x$ is a global maximiser of the function
$f / (A \psi + B \varphi)$, for some $A,B \geq 0$, then
$x$ is in the stopping region and $v(x) = f(x)$.
It is straightforward to see that the conclusions associated
with each of these cases follow immediately from the
representation (\ref{vLPintro}) of the value function $v$
(see also Corollary~\ref{no-stopping} and part~(II)
of Corollaries~\ref{call}, \ref{put} and~\ref{straddle}).
Effectively, this approach, which is summarised by
(\ref{vLPintro}), is a verification theorem of a local character.
Indeed, its application invariably involves ``guessing''
the structure of the waiting and the stopping regions.
Also, e.g., (\ref{BL1intro}) on its own does not allow
for any conclusions for initial conditions $x > d_1$
(see Example~\ref{ex5}).
It is also worth noting that, if $f$ is $C^1$, then this
approach is effectively the same as application of the
so-called ``principle of smooth fit'': first order conditions
at $d_1$ (resp., $c_2$, $d_2$) and (\ref{BL1intro}) (resp.,
(\ref{BL2intro})) yield the same equations for
$d_1$, $C_1$ (resp.\, $c_2$, $d_2$, $\kappa$, $C_2$)
as the one that the ``principle of smooth fit'' yields
(see also the generalisations in part~(III) of
Corollaries~\ref{call}, \ref{put} and~\ref{straddle}).

In stochastic analysis, a filtration can be viewed as a model for
an information flow.
Such an interpretation gives rise to the following modelling
issue.
Consider an observer whose information flow identifies with
a filtration $({\mathcal H}_t)$.
At an $({\mathcal H}_t)$-stopping time $\tau$, the observer
gets access to an additional information flow, modelled by a
filtration $({\mathcal G}_t)$, that ``switches on'' at time $\tau$.
In this context, we construct a filtration that aggregates the
two information sources available to such an observer
(see Theorem~(\ref{prop:augm-filtr})).
Building on this construction, we address the issue of pasting
weak solutions to (\ref{SDE}), or, more, generally, the issue of
pasting stopping strategies for the optimal stopping problem that
we consider, at an appropriate stopping time (see
Theorem~(\ref{prop:pasting}) and Corollary~\ref{cor:paste}).
Such a rather intuitive result is fundamental to dynamic
programming and has been assumed by several authors
in the literature (e.g., see the proof of Proposition~3.2 in
Dayanik and Karatzas~\cite{DK}).

The paper is organised as follows.
In Section~\ref{sec:pr-form}, we develop the context within
which the optimal stopping problem that we study is defined
and we list all of the assumptions we make.
Section~\ref{sec:v-finite} is concerned with a number
of preliminary results that are mostly of a technical nature.
In Section~\ref{sec:F-representations}, we derive the
representation (\ref{R-ana-intro}) for $r(\cdot)$-potentials and
the characterisation (\ref{F-repr-intro}) of differences of two
convex functions as well as a number of related results.
In Section~\ref{sec:excessive}, we consider analytic
characterisations of $r(\cdot)$-excessive functions,
while, in Section~\ref{sec:solution}, we establish our
main results on the optimal stopping problem that
we consider.
In Section~\ref{sec:sf}, we present several ramifications
of our general results on optimal stopping, including a
generalisation of the ``principle of smooth fit''.
In Section~\ref{sec:ex}, we consider a number
of illustrating examples.
Finally, we develop the theory concerned with pasting weak
solutions to (\ref{SDE}) in the Appendix.


\section{The underlying diffusion and the optimal
stopping problem}
\label{sec:pr-form}

We consider a one-dimensional diffusion with state
space an interval of the form
\ben
\CI = \mbox{} ]\alpha, \beta[ \quad \text{or} \quad
\CI = [\alpha, \beta[ \quad \text{or} \quad
\CI = \mbox{} ]\alpha, \beta] \quad \text{or} \quad
\CI = [\alpha, \beta] , \label{I}
\een
for some endpoints $-\infty \leq \alpha < \beta \leq \infty$.
Following Definition~5.20 in Karatzas and
Shreve~\cite[Chapter~5]{KS}, a weak solution to the SDE
(\ref{SDE})
in the interval $\CI$ is a collection $\str_x = (\Om, \F, \F_t,
\p_x, W, X)$ such that $(\Om, \F, \F_t, \p_x)$ is a filtered probability
space satisfying the usual conditions and supporting a standard
one-dimensional $(\F_t)$-Brownian motion $W$ and a
continuous $(\F_t)$-adapted $\CI$-valued process $X$.
The process $X$ satisfies
\begin{gather}
\int _0^{t \wedge T_{\bar{\alpha}} \wedge T_{\bar{\beta}}}
\left[ |b(X_u)| + \sigma^2 (X_u) \right] du < \infty
\label{X-integr} \\
\intertext{and}
X_{t \wedge T_{\bar{\alpha}} \wedge T_{\bar{\beta}}} = x +
\int _0^{t \wedge T_{\bar{\alpha}} \wedge T_{\bar{\beta}}}
b(X_u) \, du + \int _0^{t \wedge T_{\bar{\alpha}} \wedge
T_{\bar{\beta}}} \sigma (X_u) \, dW_u
\end{gather}
for all $t \geq 0$ and $\alpha < \bar{\alpha} < x < \bar{\beta}
< \beta$, $\p_x$-a.s..
Here, as well as throughout the paper, we denote by $T_y$
the first hitting time of the set $\{y\}$, which is defined by
\be
T_y = \inf \left\{ t \geq 0 \mid \ X_t = y \right\} ,
\quad \text{for } y \in [\alpha, \beta] ,
\ee
with the usual convention that $\inf \emptyset = \infty$.
The actual choice of the interval $\CI$ from among the four
possibilities in (\ref{I}) depends on the choice of the data
$b$ and $\sigma$ through the resulting properties of the explosion
time $T_\alpha \wedge T_\beta$ at which the process $X$ hits
the boundary $\{ \alpha, \beta \}$ of the interval $\CI$.
If the boundary point $\alpha$ (resp., $\beta$) is inaccessible,
i.e., if
\be
\p_x \bigl( T_\alpha < \infty \bigr) = 0 \quad \left( \text{resp., }
\p_x \bigl( T_\beta < \infty \bigr) = 0 \right) ,
\ee
then the interval $\CI$ is open from the left (resp., open from
the right).
If $\alpha$ (resp., $\beta$) is not inaccessible, then it is absorbing
and the interval $\CI$ is closed from the left (resp., closed from the
right).
In particular,
\ben
X_t = \begin{cases} \alpha , & \text{if } \lim _{u \rightarrow
T_\alpha \wedge T_\beta} X_u = \alpha , \\ \beta , & \text{if }
\lim _{u \rightarrow T_\alpha \wedge T_\beta} X_u = \beta ,
\end{cases} \quad \text{for all } t \geq T_\alpha \wedge T_\beta
. \label{Xafter}
\een

The following assumption ensures that the SDE (\ref{SDE})
has a weak solution in $\CI$, as described above, which is
unique in the sense of probability law (see Theorem~5.15 in
Karatzas and Shreve~\cite[Chapter~5]{KS}).

\begin{assm} \rm{
The functions $b, \sigma: \intI \rightarrow \R$ are Borel-measurable,
\begin{gather}
\sigma ^2(x) > 0 \quad \text{for all } x \in \intI \equiv \mbox{}
]\alpha , \beta [ , \label{ND} \\
\intertext{and}
\int _{\bar{\alpha}}^{\bar{\beta}} \frac{1+|b(s)|} {\sigma ^2 (s)} \,
ds < \infty \quad \text{for all } \alpha < \bar{\alpha} < \bar{\beta}
< \beta . \label{LI}
\end{gather}
} \mbox{}\hfill$\Box$ \label{A1} \end{assm}
This assumption also implies that, given $c \in \intI$ fixed,
the scale function $p$, given by
\ben
p(x) = \int _c^x \exp \left( -2 \int _c^s \frac{b(u)}{\sigma^2(u)}
\, du \right) ds, \quad \text{for } x \in \intI , \label{scale}
\een
is well-defined, and the speed measure $m$ on $\bigl(
\intI, {\mathcal B} (\CI) \bigr)$, given by
\ben
m(dx) = \frac{2}{\sigma ^2 (x) p' (x)} \, dx , \label{speed}
\een
is a Radon measure.
At this point, it is worth noting that Feller's test for explosions
provides necessary and sufficient conditions that determine
whether the solution of (\ref{SDE}) hits one or the other
or both of the boundary points $\alpha$, $\beta$ in finite time
with positive probability (see Theorem~5.29 in Karatzas and
Shreve~\cite[Chapter~5]{KS}).

We consider the optimal stopping problem, the value
function of which is defined by
\ben
v(x) = \sup _{(\str_x, \tau) \in \T_x} \e_x \left[ e^{-\Lambda
_\tau} f(X_\tau) {\bf 1} _{\{ \tau < \infty \}} \right] = \sup
_{(\str_x, \tau) \in \T_x} J(\str_x, \tau) , \quad \text{for }
x \in \CI , \label{v}
\een
where
\be
J(\str_x, \tau) = \e_x \left[ e^{-\Lambda _{\tau \wedge T_\alpha
\wedge T_\beta}} f(X_{\tau \wedge T_\alpha \wedge T_\beta})
{\bf 1} _{\{ \tau < \infty \}} \right] ,
\ee
the discounting factor $\Lambda$ is defined by (\ref{Lambda})
in the introduction,  and the set of all stopping strategies $\T_x$
is the collection
of all pairs $(\str_x, \tau)$ such that $\str_x$ is a weak solution to
(\ref{SDE}), as described above, and $\tau$ is an associated
$(\F_t)$-stopping time.

We make the following assumption, which also implies the identity
in (\ref{v}).

\begin{assm} {\rm
The reward function $f: \CI \rightarrow \R_+$ is Borel-measurable.
The discounting rate function $r : \CI \rightarrow \R_+$ is
Borel-measurable and uniformly bounded away from 0, i.e., $r(x)
\geq r_0$ for all $x \in \CI$, for some $r_0 > 0$.
Also,
\ben
\int _{\bar{\alpha}}^{\bar{\beta}} \frac{r(s)} {\sigma ^2 (s)} \,
ds < \infty \quad \text{for all } \alpha < \bar{\alpha} < \bar{\beta}
< \beta . \label{rLI}
\een
} \mbox{}\hfill$\Box$ \label{A2} \end{assm}

In the presence of Assumptions~\ref{A1} and~\ref{A2}, there
exists a pair of $C^1$ with absolutely continuous first derivatives
functions  $\varphi, \psi: \CI \rightarrow \R_+$ such
that $\varphi$ (resp., $\psi$) is strictly decreasing (resp., increasing),
and
\begin{align}
\varphi (x) & = \varphi (y) \e_x \left[ e^{-\Lambda _{T_y}}
\right] \equiv \varphi (y) \e_x \left[ e^{-\Lambda _{T_y}}
{\bf 1} _{\{ T_y < T_\beta \}} \right] \quad \text{for all }
y<x , \label{phi} \\
\psi (x) & = \psi (y) \e_x \left[ e^{-\Lambda _{T_y}} \right]
\equiv \psi (y) \e_x \left[ e^{-\Lambda _{T_y}} {\bf 1}
_{\{ T_y < T_\alpha \}} \right] \quad \text{for all } x<y
, \label{psi}
\end{align}
for every solution $\str_x$ to (\ref{SDE}).
Also,
\begin{gather}
\text{if } \alpha \text{ is absorbing, then } \varphi (\alpha) :=
\lim _{x \downarrow \alpha} \varphi (x) < \infty \text{ and }  \psi
(\alpha) := \lim _{x \downarrow \alpha} \psi (x) = 0 ,
\label{phi-psi-aabs} \\
\text{if } \beta \text{ is absorbing, then } \varphi (\beta) :=
\lim _{x \uparrow \beta} \varphi (x) = 0 \text{ and } \psi (\beta)
:= \lim _{x \uparrow \beta} \psi (x) < \infty ,
\label{phi-psi-babs} \\
\text{and, if } \alpha \text{ (resp., } \beta \text{) is inaccessible,
then } \lim _{x \downarrow \alpha} \varphi (x) = \infty \text{ (resp., }
\lim _{x   \uparrow \beta} \psi (x) = \infty \text{)}
. \label{phi-psi-nat}
\end{gather}
An inspection of these facts reveals that, in all cases,
\ben
\lim _{y \downarrow \alpha} \frac{\psi (y)}{\varphi (y)} =
\lim _{y \uparrow \beta} \frac{\varphi (y)}{\psi (y)} = 0 .
\label{phi-psi->0}
\een
The functions $\varphi$ and $\psi$ are classical solutions to
the homogeneous ODE
\ben
\half \sigma ^2 (x) g'' (x) + b (x) g' (x) - r(x) g(x) = 0 ,
\label{ODEhom}
\een
and satisfy
\ben
\varphi (x) \psi' (x) - \varphi' (x) \psi (x) = C p' (x) \quad
\text{for all } x \in \CI , \label{C}
\een
where $C = \varphi (c) \psi' (c) - \varphi' (c) \psi (c)$ and $p$
is the scale function defined by (\ref{scale}).
Furthermore, given any solution $\str_x$ to (\ref{SDE}),
\ben
\text{the processes } \left( e^{-\Lambda_t} \varphi (X_t) \right)
\text{ and } \left( e^{-\Lambda_t} \psi (X_t) \right) \text{ are
local martingales} . \label{phi-psi-lms}
\een
The existence of these functions and their properties that we have
listed can be found in several references, including Borodin and
Salminen~\cite[Section~II.1]{BS}, Breiman~\cite[Chapter~16]{B},
and It\^{o} and McKean~\cite[Chapter~4]{IM}.

\section{Preliminary considerations}
\label{sec:v-finite}

Throughout this section, we assume that a weak solution
$\str_x$ to (\ref{SDE}) has been associated with each initial
condition $x \in \intI$.
We first need to introduce some notation.
To this end, we recall that, if $g: \intI \rightarrow \R$ is a
function that is the difference of two convex functions, then
its left-hand side first derivative $g_-'$ exists and is a function
of finite variation, and its second distributional derivative $g''$
is a measure.
We denote by
\ben
g''(dx) = g_\ac '' (x) \, dx + g_\s (dx) \label{Leb-dec}
\een
the Lebesgue decomposition of the second distributional
derivative $g''(dx)$ into the measure $g_\ac '' (x) \, dx$
that is absolutely continuous with respect to the Lebesgue
measure and the measure $g_\s '' (dx)$ that is mutually
singular with the Lebesgue measure.
Note that the function $g_\ac ''$ identifies with the
``classical'' sense second derivative of $g$, which exists
Lebesgue-a.e..
In view of these observations and notation, we define the
measure $\Lop g$ on $\bigl( \intI , {\mathcal B} (\intI)
\bigr)$ and the function $\Lop _\ac g : \intI \rightarrow \R$
by
\begin{gather}
\Lop g (dx) = \half \sigma ^2 (x) g'' (dx) + b (x) g_-' (x)
\, dx - r(x) g(x) \, dx \label{lop} \\
\intertext{and}
\Lop_\ac g (x) = \half \sigma ^2 (x) g_\ac'' (x) + b (x) g_-' (x)
- r(x) g(x) . \label{lopac}
\end{gather}

Given a Radon measure $\mu$ on $\bigl( \intI , {\mathcal B}
(\intI) \bigr)$ such that $\sigma^{-2}$ is locally integrable with
respect to $|\mu|$, we consider the continuous additive functional
$A^\mu$ defined by (\ref{A}) in the introduction.
Given any $t < T_\alpha \wedge T_\beta$, $A_t^\mu$ is
well-defined and real-valued because $\alpha < \inf _{s \leq t}
X_s < \sup _{s \leq t} X_s < \beta$ and the process $L^z$
increases on the set $\{ X_s = z \}$.
Also, since $L^z$ is an increasing process, $A^\mu$
(resp., $-A^\mu$) is an increasing process if $\mu$
(resp., $-\mu$) is a positive measure.
The following result is concerned with various properties of the
process $A^\mu$ that we will need.

\begin{lem} \label{lem:Amu}
Let $\mu$ be a Radon measure on $\bigl( \intI, {\mathcal B}
(\intI) \bigr)$ such that $\sigma^{-2}$ is locally integrable with
respect to $|\mu|$, consider any increasing sequence of real-valued
Borel-measurable functions $(\zeta_n)$ on $\CI$ such that
\ben
0 \leq \zeta_n (z) \leq 1 \quad \text{and} \quad \lim _{n \rightarrow
\infty} \zeta_n (z) = 1 , \quad \mu \text{-a.e.} , \label{zeta}
\een
and denote by $\mu_n$ the measure defined by
\ben
\mu_n (\Gamma) = \int _\Gamma \zeta_n (z) \, \mu (dz) , \quad
\text{for } \Gamma \in {\mathcal B} (\intI) . \label{mun-defn}
\een
$A^{|\mu|}$ is a continuous increasing process,
\ben
A^{\mu} = -A^{-\mu} = A^{\mu^+} - A^{\mu^-} , \qquad
A^{|\mu|} = A^{\mu^+} + A^{\mu^-} , \label{A-lin}
\een
and
\ben
\lim _{n \rightarrow \infty} \e_x \left[ \int _0^{T_\alpha \wedge
T_\beta} e^{-\Lambda_t} \, dA_t^{|\mu_n|} \right] = \e_x \left[
\int _0^{T_\alpha \wedge T_\beta} e^{-\Lambda_t} \, dA_t^{|\mu|}
\right] \quad \text{for all } x \in \intI . \label{A-mun}
\een
\end{lem}
{\bf Proof.}
The process $A^{|\mu|}$ is continuous and increasing
because this is true for the local time process $L^z$ for all
$z \in \CI$.
Also, (\ref{A-lin}) can be seen by a simple inspection of the
definition (\ref{A}) of $A^\mu$.
To prove (\ref{A-mun}), we have to show that, given any
$x \in \intI$,
\ben
\lim _{n \rightarrow \infty} \e_x \bigl[ I_{T_\alpha \wedge
T_\beta}^{(n)} \bigr] = \e_x \bigl[ I_{T_\alpha \wedge T_\beta}
\bigr] , \label{In-conv}
\een
where
\be
I_t^{(n)} = \int _0^t e^{-\Lambda _u} \, dA_u^{|\mu_n|}
\quad \text{and} \quad I_t = \int _0^t e^{-\Lambda_u}
\, dA_u^{|\mu|} , \quad \text{for } t \in [0, T_\alpha \wedge
T_\beta] .
\ee
To this end, we note that (\ref{zeta}) and the monotone convergence
theorem imply that the sequence $( A_t^{|\mu_n|} )$ increases to
$A_t^{|\mu|}$ for all $t < T_\alpha \wedge T_\beta$ as
$n \rightarrow \infty$,
because
\be
A_t^{|\mu_n|} = \int _\alpha^\beta \frac{L_t^z} {\sigma ^2 (z)}
\, |\mu_n| (dz) = \int _\alpha^\beta \frac{L_t^z} {\sigma ^2 (z)}
\zeta_n (z) \, |\mu| (dz) , \quad \text{for } t \in [0, T_\alpha
\wedge T_\beta[ .
\ee
Also, we use the integration by parts formula to calculate
\ben
\int _0^t e^{-\Lambda_u} \, dA_u^{|\mu_n|} = e^{-\Lambda _t}
A_t^{|\mu_n|} + \int _0^t e^{-\Lambda_u} r(X_u) A_u^{|\mu_n|}
\, du , \quad \text{for } t \in [0, T_\alpha \wedge T_\beta[ .
\label{A-IBP}
\een
In view of these observations and the monotone convergence
theorem, we can see that
\begin{gather}
0 \leq I_t^{(n)} \leq I_t^{(n+1)} \quad \text{for all } t \in
[0, T_\alpha \wedge T_\beta] \text{ and } n \geq 1 ,
\label{Iprop1} \\
\intertext{and}
\lim _{n \rightarrow \infty} I_t^{(n)} = I_t \quad \text{for all }
t \in [0, T_\alpha \wedge T_\beta[ , \label{Iprop2}
\end{gather}
Combining these results with the fact that the positive
processes $I^{(n)}$ are increasing, we can see that
\begin{gather}
I_{T_\alpha \wedge T_\beta} = \lim _{t \rightarrow T_\alpha
\wedge T_\beta} I_t \geq \lim _{t \rightarrow T_\alpha \wedge
T_\beta} I_t^{(n)} = I_{T_\alpha \wedge T_\beta}^{(n)} \quad
\text{for all } n \geq 1 \nonumber \\
\intertext{and}
I_{T_\alpha \wedge T_\beta} =  \lim _{t \rightarrow T_\alpha
\wedge T_\beta} I_t = \lim _{t \rightarrow T_\alpha
\wedge T_\beta} \lim _{n \rightarrow \infty} I_t^{(n)} \leq
\lim _{n \rightarrow \infty} I_{T_\alpha \wedge T_\beta}^{(n)}
. \nonumber
\end{gather}
It follows that
$\lim _{n \rightarrow \infty} I_{T_\alpha \wedge T_\beta}^{(n)}
= I_{T_\alpha \wedge T_\beta}$, which, combined with
monotone convergence theorem, implies (\ref{In-conv})
and the proof is complete.
\mbox{}\hfill$\Box$
\vspace{5mm}

We will need the results derived in the following lemma, the
proof of which is based on the It\^{o}-Tanaka-Meyer formula.

\begin{lem} \label{lem:Ito}
If $F: \intI \rightarrow \R$ is a function that is the difference
of two convex functions, then the following statements are
true:

\noindent
{\rm (I)}
The increasing process $A^{|\Lop F|}$ is real-valued, and
\ben
e^{-\Lambda_t} F(X_t) = F(x) + \int _0^t e^{-\Lambda _u}
\, d A_u^{\Lop F} + \int _0^t e^{-\Lambda _u} \sigma (X_u) F_-'
(X_u) \, dW_u , \quad \text{for } t \in [0, T_\alpha \wedge
T_\beta] . \label{ITO!}
\een

\noindent
{\rm (II)}
If $F$ is $C^1$ with absolutely continuous with respect to the
Lebesgue measure first derivative, i.e., if $\Lop F (dx) =
\Lop_\ac F(x) \, dx$ in the notation of (\ref{lop})--(\ref{lopac}),
then
\ben
\int _0^t e^{-\Lambda_u} \, dA_u^{\Lop F} = \int _0^t
e^{-\Lambda_u} \Lop_\ac F (X_u) \, du , \quad \text{for } t
\in [0, T_\alpha \wedge T_\beta] . \label{OTF}
\een
\end{lem}
{\bf Proof.}
In view of the Lebesgue decomposition of the second distributional
derivative $F''(dx)$ of $F$ as in (\ref{Leb-dec}) and the occupation
times formula
\be
\int _\alpha^\beta L_t^z F_\ac'' (z) \, dz = \int _0^t \sigma
^2 (X_u) F_\ac'' (X_u) \, du ,
\ee
we can see that the It\^o-Tanaka-Meyer formula
\ben
F (X_t) = F(x) + \int _0^t b(X_u) F_-' (X_u) \, du + \half
\int _\alpha^\beta L_t^z \, F''(dz ) + \int _0^t \sigma (X_u)
F_-' (X_u) \, dW_u \nonumber 
\een
implies that
\begin{align}
F(X_t) = \mbox{} & F(x) + \int _0^t \left[ \half \sigma ^2
(X_u) F_\ac'' (X_u) + b(X_u) F_-' (X_u) \right] du + \half \int
_\alpha^\beta L_t^z \, F_\s'' (dz) \nonumber \\
& + \int _0^t \sigma (X_u) F_-' (X_u) \, dW_u . \label{ITO2}
\end{align}
Combining this expression with the definition (\ref{lopac}) of
$\Lop_\ac$, we can see that
\begin{align}
F(X_t) = \mbox{} & F(x) + \int _0^t r(X_u) F(X_u) \, du + \int
_0^t \Lop _\ac F(X_u) \, du + \half \int _\alpha^\beta L_t^z \,
F_\s'' (dz) \nonumber \\
& + \int _0^t \sigma (X_u) F_-' (X_u) \, dW_u . \label{ITO3}
\end{align}
Using the occupation times formula once again and the
definitions (\ref{lop}), (\ref{lopac}) of $\Lop$, $\Lop_\ac$, we
can see that
\begin{align}
\int _0^t \Lop _\ac F(X_u) \, du + \half \int _\alpha^\beta
L_t^z \, F_\s'' (dz) & = \int _\alpha^\beta \frac{L_t^z}{\sigma^2
(z)} \Lop _\ac F(z) \, dz + \int _\alpha^\beta \frac{L_t^z}
{\sigma^2 (z)} \half \sigma^2 (z) \, F_\s'' (dz) \nonumber \\
& = \int _\alpha^\beta \frac{L_t^z}{\sigma^2 (z)} \, \Lop F
(dz) \nonumber \\
& = A_t^{\Lop F} . \label{OT}
\end{align}
The validity of It\^{o}-Tanaka-Meyer's and the occupation times
formulae and (\ref{ITO3})--(\ref{OT}) imply that the process
$A^{\Lop F}$ is well-defined and real-valued.
Also, (\ref{ITO!}) follows from the definition (\ref{Lambda}) of
the process $\Lambda$, (\ref{ITO3})--(\ref{OT}) and an application
of the integration by parts formula.

If $\Lop F (dx) = \Lop_\ac F(x) \, dx$, the definition of
$A^{\Lop F}$ and the occupation times formula imply that
\be
A_t^{\Lop F} = \int_0^t \Lop_\ac F (X_u) \, du ,
\ee
and (\ref{OTF}) follows.
\mbox{}\hfill$\Box$
\vspace{5mm}

The next result is concerned with a form of Dynkin's formula
that the functions $\varphi$, $\psi$ satisfy as well as with a
pair of expressions that become useful when explicit solutions
to special cases of the general optimal stopping problem
are explored (see Section~\ref{sec:sf}).

\begin{lem} \label{lem:phi-psi-expect}
The functions $\varphi$, $\psi$ introduced by (\ref{phi}),
(\ref{psi}) satisfy
\ben
\varphi (x) = \e_x \left[ e^{-\Lambda _{\tau \wedge T_{\bar{\alpha}}
\wedge T_{\bar{\beta}}}} \varphi (X_{\tau \wedge T_{\bar{\alpha}}
\wedge T_{\bar{\beta}}}) \right] \quad \text{and} \quad
\psi (x) = \e_x \left[ e^{-\Lambda _{\tau \wedge T_{\bar{\alpha}}
\wedge T_{\bar{\beta}}}} \psi (X_{\tau \wedge T_{\bar{\alpha}}
\wedge T_{\bar{\beta}}}) \right] \label{phi-psi-expect}
\een
for all stopping times $\tau$ and all points $\bar{\alpha} < x <
\bar{\beta}$ in $\CI$.
Furthermore,
\begin{gather}
\e_x \left[ e^{-\Lambda _{T_{\bar{\alpha}}}} {\bf 1} _{\{
T_{\bar{\alpha}} < T_{\bar{\beta}} \}} \right] = \frac{\varphi
(\bar{\beta}) \psi (x) - \varphi (x) \psi (\bar{\beta})} {\varphi
(\bar{\beta}) \psi (\bar{\alpha}) - \varphi (\bar{\alpha}) \psi
(\bar{\beta})} \label{ELam_alpha} \\
\intertext{and}
\e_x \left[ e^{-\Lambda _{T_{\bar{\beta}}}} {\bf 1} _{\{
T_{\bar{\beta}} < T_{\bar{\alpha}} \}} \right] = \frac{\varphi
(x) \psi (\bar{\alpha}) - \varphi (\bar{\alpha}) \psi (x)} {\varphi
(\bar{\beta}) \psi (\bar{\alpha}) - \varphi (\bar{\alpha}) \psi
(\bar{\beta})} . \label{ELam_beta}
\end{gather}
\end{lem}
{\bf Proof.}
Combining (\ref{ITO!}) with the fact that $\Lop \varphi = 0$, we
can see that
\ben
e^{-\Lambda _{\tau \wedge T_{\bar{\alpha}} \wedge T_{\bar{\beta}}}}
\varphi (X_{\tau \wedge T_{\bar{\alpha}} \wedge T_{\bar{\beta}}})
= \varphi (x) + M_{\tau \wedge T_{\bar{\alpha}} \wedge
T_{\bar{\beta}}} , \label{phi-Ito}
\een
where
\be
M_t = \int _0^t e^{-\Lambda _u} \sigma (X_u) \varphi ' (X_u)
\, dW_u .
\ee
In view of (\ref{phi-psi-aabs}) and the fact that the positive
function $\varphi$ is decreasing, we can see that $\sup _{y \in
[\bar{\alpha}, \bar{\beta}]} \varphi (y) < \infty$.
Therefore, $M^{T_{\bar{\alpha}} \wedge T_{\bar{\beta}}}$ is a
uniformly integrable martingale because it is a uniformly bounded
local martingale.
It follows that $\e_x \bigl[ M_{\tau \wedge T_{\bar{\alpha}} \wedge
T_{\bar{\beta}}} \bigr] = 0$ and (\ref{phi-Ito}) implies the first
identity in (\ref{phi-psi-expect}).
The second identity in (\ref{phi-psi-expect}) can be established
using similar arguments.

Finally, (\ref{ELam_alpha}) and (\ref{ELam_beta}) follow
immediately once we observe that they are equivalent to
the system of equations
\begin{gather}
\varphi (x) = \varphi (\bar{\alpha}) \e_x \left[ e^{-\Lambda
_{T_{\bar{\alpha}}}} {\bf 1} _{\{ T_{\bar{\alpha}} < T_{\bar{\beta}}
\}} \right] + \varphi (\bar{\beta}) \e_x \left[ e^{-\Lambda
_{T_{\bar{\beta}}}} {\bf 1} _{\{ T_{\bar{\beta}} < T_{\bar{\alpha}}
\}} \right] \nonumber \\
\intertext{and}
\psi (x) = \psi (\bar{\alpha}) \e_x \left[ e^{-\Lambda
_{T_{\bar{\alpha}}}} {\bf 1} _{\{ T_{\bar{\alpha}} < T_{\bar{\beta}}
\}} \right] + \psi (\bar{\beta}) \e_x \left[ e^{-\Lambda
_{T_{\bar{\beta}}}} {\bf 1} _{\{ T_{\bar{\beta}} < T_{\bar{\alpha}}
\}} \right] , \nonumber
\end{gather}
which holds true thanks to (\ref{phi-psi-expect}) for $\tau \equiv
\infty$.
\mbox{}\hfill$\Box$
\vspace{5mm}

We conclude this section with a necessary and sufficient
condition for the value function of our optimal stopping problem
to be finite.

\begin{lem} \label{lem:v<oo}
Consider the optimal stopping problem formulated in
Section~\ref{sec:pr-form}, and let $\overline{f}$ be defined
by (\ref{f-over}) in the introduction.
If
\ben
\overline{f} : \CI \rightarrow \R_+ \text{ is real-valued} , \quad
\limsup _{y \downarrow \alpha} \frac{f(y)}{\varphi (y)} < \infty
\quad \text{and} \quad
\limsup _{y \uparrow \beta} \frac{f(y)}{\psi (y)} < \infty ,
\label{v<oo-cond}
\een
then $v(x) < \infty$ for all $x \in \CI$,
\ben
\limsup _{y \downarrow \alpha} \frac{v(y)}{\varphi (y)} =
\limsup _{y \downarrow \alpha} \frac{f(y)}{\varphi (y)}
\quad \text{and} \quad
\limsup _{y \uparrow \beta} \frac{v(y)}{\psi (y)} =
\limsup _{y \uparrow \beta} \frac{f(y)}{\psi (y)} .
\label{v-phi-psi-lims}
\een
If any of the conditions in (\ref{v<oo-cond}) is not true, then
$v(x) = \infty$ for all $x \in \intI$.
\end{lem}
{\bf Proof.}
If (\ref{v<oo-cond}) is true, then we can see that
\be
\sup _{u \leq y} \frac{f(u)}{\varphi (u)} < \infty
\quad \text{and} \quad
\sup _{u \geq y} \frac{f(u)}{\psi (u)} < \infty
\quad \text{for all } y \in \CI .
\ee
Also,
\be
f(x) \leq \sup _{u \leq y} \frac{f(u)}{\varphi (u)} \varphi
(x) + \sup _{u \geq y} \frac{f(u)}{\psi (u)} \psi (x) \quad
\text{for all } x,y \in \CI .
\ee
In view of (\ref{phi-psi-lms}), the processes $\left( e^{-\Lambda_t}
\varphi (X_t) \right)$ and $\left( e^{-\Lambda_t} \psi (X_t) \right)$
are positive supermartingales.
It follows that, given any stopping strategy $(\str_x, \tau)
\in \T_x$,
\begin{align}
J(\str_x, \tau) \leq \mbox{} & \sup _{u \leq y} \frac{f(u)}{\varphi (u)}
\e_x \left[ e^{-\Lambda_{\tau \wedge T_\alpha \wedge T_\beta}}
\varphi (X_{\tau \wedge T_\alpha \wedge T_\beta}) {\bf 1}
_{\{ \tau < \infty \}} \right] \nonumber \\
& + \sup _{u \geq y} \frac{f(u)}{\psi (u)} \e_x \left[ e^{-\Lambda
_{\tau \wedge T_\alpha \wedge T_\beta}} \psi (X_{\tau \wedge
T_\alpha \wedge T_\beta}) {\bf 1} _{\{ \tau < \infty \}} \right]
\nonumber \\
\leq \mbox{} & \sup _{u \leq y} \frac{f(u)}{\varphi (u)} \varphi
(x) + \sup _{u \geq y} \frac{f(u)} {\psi (u)} \psi (x) , \label{v<oo1}
\end{align}
which implies that $v(x) < \infty$.

To show the first identity in (\ref{v-phi-psi-lims}), we
note that (\ref{v<oo1}) implies that
\be
\frac{v(x)}{\varphi (x)} \leq \sup _{u \leq y} \frac{f(u)}
{\varphi (u)} + \sup _{u \geq y} \frac{f(u)}{\psi (u)}
\frac{\psi (x)}{\varphi (x)} .
\ee
Combining this calculation with (\ref{phi-psi->0}), we obtain
\be
\limsup _{x \downarrow \alpha} \frac{v(x)}{\varphi (x)} \leq
\sup _{u \leq y} \frac{f(u)}{\varphi (u)} ,
\ee
which implies that $\limsup _{y \downarrow \alpha} v(y) / 
\varphi (y) \leq \limsup _{y \downarrow \alpha} f(y) / \varphi
(y)$.
The reverse inequality follows immediately from the fact
that $v \geq f$.
The second identity in (\ref{v-phi-psi-lims}) can be established
using similar arguments.

If the problem data is such that the first limit in
(\ref{v<oo-cond}) is infinite, then we consider any initial
condition $x \in \intI$ and any sequence $(y_n)$ in $\CI$
such that $y_n < x$ for all $n \geq 1$ and
$\lim _{n \rightarrow \infty} f(y_n) / \varphi (y_n) = \infty$.
We can then see that
\be
v(x) \geq \lim _{n \rightarrow \infty} J(\str_x, T_{y_n}) \geq
\lim _{n \rightarrow \infty} f(y_n) \e_x \left[ e^{-\Lambda
_{T_{y_n}}} \right] \stackrel{(\ref{phi})}{=} \lim _{n \rightarrow
\infty} \frac{f(y_n) \varphi (x)} {\varphi (y_n)} = \infty ,
\ee
where $\str_x$ is any solution to (\ref{SDE}).
Similarly, we can see that $v(x) = \infty$ for all $x \in \intI$
if the second limit in (\ref{v<oo-cond}) is infinite or if there
exists a point $y \in \intI$ such that $\overline{f} (y) = \infty$.
\mbox{}\hfill$\Box$

\section{$\pmb{r(\cdot)}$-potentials and differences of two
convex functions}
\label{sec:F-representations}

Throughout this section, we assume that a weak solution
$\str_x$ to (\ref{SDE}) has been associated with each initial
condition $x \in \intI$.
Accordingly, whenever we consider a stopping time $\tau$,
we refer to a stopping time of the filtration in the solution
$\str_x$.

We first characterise the limiting behaviour at the boundary
of $\CI$ of a difference of two convex functions on $\intI$,
and we show that such a function satisfies Dynkin's formula
under appropriate assumptions.

\begin{lem} \label{prop:Dynkin}
Consider any function $F: \intI \rightarrow \R$ that is a difference
of two convex functions and is such that
\ben
\limsup _{y \downarrow \alpha} \frac{|F(y)|}{\varphi (y)} <
\infty \quad \text{ and } \quad \limsup _{y \uparrow \beta}
\frac{|F(y)|} {\psi (y)} < \infty . \label{F(a,b)}
\een

\noindent
{\rm (I)}
If $-\Lop F$ is a positive measure, then
\ben
\e_x \left[ \int _0^{T_\alpha \wedge T_\beta} e^{-\Lambda_t}
\, dA_t^{|\Lop F|} \right] < \infty \quad \text{for all } x
\in \intI . \label{F-LF-IC}
\een

\noindent
{\rm (II)}
If $F$ satisfies
\ben
\e_x \left[ \int _0^{T_\alpha \wedge T_\beta} e^{-\Lambda_t}
\, dA_t^{|\Lop F|} \right] < \infty , \quad \text{for some } x
\in \intI , \label{F-LF-IC-x}
\een
then both of the limits
$\lim _{y \downarrow \alpha} F(y) / \varphi (y)$ and \/
$\lim _{y \uparrow \beta} F(y) / \psi (y)$ exist.

\noindent
{\rm (III)}
Suppose that $F$ satisfies (\ref{F-LF-IC-x}),
\ben
\lim _{y \downarrow \alpha} \frac{F(y)}{\varphi (y)} = 0
\quad \text{and } \quad \lim _{y \uparrow \beta} \frac{F(y)}
{\psi (y)} = 0  . \label{natMthm}
\een
If $x \in \intI$ is an initial condition such that (\ref{F-LF-IC-x})
is true, then
\begin{align}
\e_x \left[ e^{-\Lambda _\tau} F(X_\tau) {\bf 1} _{\{ \tau <
T_\alpha \wedge T_\beta \}} \right] & = F(x) + \e_x \left[ \int
_0^{\tau  \wedge T_\alpha \wedge T_\beta} e^{-\Lambda _t} \,
dA_t^{\Lop F} \right] \nonumber \\
& = \e_x \left[ e^{-\Lambda _{\tau \wedge T_\alpha \wedge T_\beta}}
F(X_{\tau \wedge T_\alpha \wedge T_\beta}) {\bf 1}
_{\{ \tau \wedge T_\alpha \wedge T_\beta < \infty \}} \right]
\label{Dynkin}
\end{align}
for every stopping time $\tau$; in the last identity here, we
assume that 
\be
F(\alpha) = \lim _{y \downarrow \alpha} F(y) = 0 \quad
\left( \text{resp., } F(\beta) = \lim _{y \uparrow \beta} F(y) = 0
\right)
\ee
if $\alpha$ (resp., $\beta$) is absorbing, namely, if $\p_x (T_\alpha
< \infty) > 0$ (resp., $\p_x (T_\beta < \infty) > 0$), consistently
with (\ref{natMthm}).
\end{lem}
{\bf Proof.}
Throughout the proof, $\tau$ denotes any stopping time.
Recalling (\ref{ITO!}) in Lemma~\ref{lem:Ito}, we write
\ben
e^{-\Lambda_t} F(X_t) = F(x) + \int _0^t e^{-\Lambda _u}
\, d A_u^{\Lop F} + M_t , \label{ITO-CFR}
\een
where $M$ is the stochastic integral defined by
\be
M_t = \int _0^t e^{-\Lambda _u} \sigma (X_u) F_-' (X_u) \,
dW_u .
\ee
We consider any decreasing sequence $(\alpha_n)$ and any
increasing sequence $(\beta_n)$ such that
\ben
\alpha < \alpha_n < x < \beta_n < \beta \, \text{ for all } n \geq 1 ,
\quad \lim _{n \rightarrow \infty} \alpha _n = \alpha \quad \text{and}
\quad \lim _{n \rightarrow \infty} \beta_n = \beta . \label{ab-seq}
\een
Also, we define
\ben
\tau_\ell (\bar{\alpha}, \bar{\beta}) = \inf \left\{ t \geq 0 \,
\bigg| \ \int _0 ^{t \wedge T_{\bar{\alpha}} \wedge T_{\bar{\beta}}}
\sigma^2 (X_u) \, du \geq \ell \right\} \wedge T_{\bar{\alpha}}
\wedge T_{\bar{\beta}} , \label{taul}
\een
where we adopt the usual convention that $\inf \emptyset =
\infty$, and we note that the definition and the construction of a
weak solution to (\ref{SDE}) (see Definition~5.5.20 in
Karatzas and Shreve~\cite{KS}) imply that these stopping times
satisfy
\ben
\tau_\ell (\bar{\alpha}, \bar{\beta}) > 0 \text{ for all } \ell \geq 1
\quad \text{and} \quad \lim _{\ell \rightarrow \infty} \tau
_\ell (\bar{\alpha}, \bar{\beta}) = T_{\bar{\alpha}} \wedge
T_{\bar{\beta}} . \label{taulim}
\een
The function $F_-'$ is locally bounded because it is of finite
variation.
Therefore, we can use It\^{o}'s isometry to calculate
\begin{align}
\e_x \left[ M_{\tau \wedge \tau_\ell (\alpha_m, \beta_n)}^2 \right]
& = \e_x \left[ \int _0^{\tau \wedge \tau_\ell (\alpha_m, \beta_n)}
\left[ e^{-\Lambda _u} \sigma (X_u) F_-' (X_u) \right]^2 du \right]
\nonumber \\
& \leq \sup _{y \in [\alpha_n, \beta_n]} \left[ F_-' (y) \right] ^2
\e_x \left[ \int _0^{\tau_\ell (\alpha_m, \beta_n)}  \sigma ^2 (X_u)
\, du \right] \nonumber \\
& \leq \ell \sup _{y \in [\alpha_n, \beta_n]} \left[ F_-' (y) \right] ^2
\nonumber \\
& < \infty ,
\end{align}
which implies that the stopped process $M^{\tau \wedge \tau_\ell
(\alpha_m, \beta_n)}$ is a uniformly integrable martingale.
Combining this observation with (\ref{ITO-CFR}), we can see
that
\be
\e_x \left[ e^{-\Lambda _{\tau \wedge \tau_\ell (\alpha_m,
\beta_n)}} F (X_{\tau \wedge \tau_\ell (\alpha_m, \beta_n)}) \right]
= F(x) + \e_x \left[ \int _0^{\tau \wedge \tau_\ell (\alpha_m, \beta_n)}
e^{-\Lambda_u} \, dA _u^{\Lop F} \right] .
\ee
In view of (\ref{taulim}) and the local boundedness of $F$,
we can pass to the limit using the dominated convergence theorem
to obtain
\begin{align}
F(x) + \lim _{\ell \rightarrow \infty} \e_x & \left[ \int _0^{\tau \wedge
\tau_\ell (\alpha_m, \beta_n)} e^{-\Lambda_u} \, dA _u^{\Lop F}
\right] \nonumber \\
= \mbox{} & \e_x \left[ e^{-\Lambda _{\tau \wedge T_{\alpha_m}
\wedge T_{\beta_n}}} F(X_{\tau \wedge T_{\alpha_m} \wedge
T_{\beta_n}}) \right] \nonumber \\
= \mbox{} & \e_x \left[ e^{-\Lambda _\tau} F(X_\tau) {\bf 1}
_{\{ \tau \leq T_{\alpha_m} \wedge T_{\beta_n} \}} \right]
+ F(\alpha_m) \e_x \bigl[ e^{-\Lambda _{T_{\alpha_m}}}
{\bf 1} _{\{ T_{\alpha_m} < \tau \wedge T_{\beta_n} \}} \bigr]
 \nonumber \\
& + F(\beta_n) \e_x \bigl[ e^{-\Lambda _{T_{\beta_n}}} {\bf 1}
_{\{ T_{\beta_n} < \tau \wedge T_{\alpha_m} \}} \bigr]
\nonumber \\
= \mbox{} & \e_x \left[ e^{-\Lambda _\tau} F(X_\tau) {\bf 1}
_{\{ \tau \leq T_{\alpha_m} \wedge T_{\beta_n} \}} \right] + \varphi
(x) \frac{F(\alpha_m)} {\varphi (\alpha_m)} \frac{\e_x \bigl[
e^{-\Lambda _{T_{\alpha_m}}} {\bf 1} _{\{ T_{\alpha_m} <
\tau \wedge T_{\beta_n} \}} \bigr]} {\e_x \left[ e^{-\Lambda
_{T_{\alpha_m}}} \right]} \nonumber \\
& + \psi (x) \frac{F(\beta_n)} {\psi (\beta_n)} \frac{\e_x \bigl[
e^{-\Lambda _{T_{\beta_n}}} {\bf 1} _{\{ T_{\beta_n} <
\tau \wedge T_{\alpha_m} \}} \bigr]} {\e_x \bigl[ e^{-\Lambda
_{T_{\beta_n}}} \bigr]} , \label{Dynkin-CFR}
\end{align}
the last identity following thanks to (\ref{phi})--(\ref{psi}).

{\em Proof of} (I).
If $-\Lop F$ is a positive measure, then $-A^{\Lop F} = A^{-\Lop F}
= A^{|\Lop F|}$ is an increasing process.
Therefore, we can use (\ref{taulim}), (\ref{ab-seq}) and the
monotone convergence theorem to calculate
\begin{align}
\lim _{m,n \rightarrow \infty} \lim _{\ell \rightarrow \infty}
\e_x \left[ \int _0^{\tau_\ell (\alpha_m, \beta_n)} e^{-\Lambda_u}
\, dA _u^{\Lop F} \right] & = \lim _{m,n \rightarrow \infty} \e_x
\left[ \int _0^{T_{\alpha_m} \wedge T_{\beta_n}} e^{-\Lambda_u}
\, dA _u^{\Lop F} \right] \nonumber \\
& = \e_x \left[ \int _0^{T_\alpha \wedge T_\beta} e^{-\Lambda_u}
\, dA _u^{\Lop F} \right] . \nonumber
\end{align}
Combining this with assumption (\ref{F(a,b)}), the inequalities
\ben
0 < \frac{\e_x \bigl[ e^{-\Lambda _{T_{\alpha_m}}} {\bf 1}
_{\{ T_{\alpha_m} < T_{\beta_n} \}} \bigr]} {\e_x \bigl[
e^{-\Lambda _{T_{\alpha_m}}} \bigr]} \leq 1
\quad \text{and} \quad
0 < \frac{\e_x \bigl[ e^{-\Lambda _{T_{\beta_n}}} {\bf 1}
_{\{ T_{\beta_n} < T_{\alpha_m} \}} \bigr]} {\e_x \bigl[
e^{-\Lambda _{T_{\beta_n}}} \bigr]} \leq 1 , \label{exp-lims}
\een
and  (\ref{Dynkin-CFR}) for $\tau = \infty$, we can see that
\begin{align}
0 \leq \mbox{} & \e_x \left[ \int _0^{T_\alpha \wedge T_\beta}
e^{-\Lambda_t} \, dA_t^{|\Lop F|} \right] \nonumber \\
= \mbox{} & - \e_x \left[ \int _0^{T_\alpha \wedge T_\beta}
e^{-\Lambda_u} \, dA _u^{\Lop F} \right] \nonumber \\
= \mbox{} & \lim _{m,n \rightarrow \infty} \biggl( F(x) -
\varphi (x) \frac{F(\alpha_m)} {\varphi (\alpha_m)}
\frac{\e_x \bigl[ e^{-\Lambda _{T_{\alpha_m}}} {\bf 1}
_{\{ T_{\alpha_m} < T_{\beta_n} \}} \bigr]}
{\e_x \bigl[ e^{-\Lambda _{T_{\alpha_m}}} \bigr]} - \psi (x)
\frac{F(\beta_n)} {\psi (\beta_n)} \frac{\e_x \bigl[
e^{-\Lambda _{T_{\beta_n}}} {\bf 1} _{\{ T_{\beta_n} <
T_{\alpha_m} \}} \bigr]} {\e_x \bigl[ e^{-\Lambda _{T_{\beta_n}}}
\bigr]} \biggr) \nonumber \\
\leq \mbox{} & |F(x)| + \varphi (x) \limsup _{m \rightarrow
\infty} \frac{|F(\alpha_m)|} {\varphi (\alpha_m)} + \psi (x)
\limsup _{n \rightarrow \infty} \frac{|F(\beta_n)|} {\psi
(\beta_n)} \nonumber \\
< \mbox{} & \infty .
\end{align}

{\em Proof of} (II).
We now fix any initial condition $x \in \intI$ such that
(\ref{F-LF-IC-x}) is true and we assume that the sequence
$(\alpha_m)$ has been chosen so that
\ben
\lim _{m \rightarrow \infty} \frac{F(\alpha_m)}{\varphi (\alpha_m)}
\ \text{ exists} . \label{Fab-lim}
\een

In light of (\ref{A-lin}) in Lemma~\ref{lem:Amu} and (\ref{taulim}),
we can see that the dominated convergence theorem implies
that
\ben
\lim _{m,n \rightarrow \infty} \lim _{\ell \rightarrow \infty} \e_x
\left[ \int _0^{\tau \wedge \tau_\ell (\alpha_m, \beta_n)}
e^{-\Lambda _u} \, dA _u^{\Lop F} \right] = \e_x \left[ \int
_0^{\tau \wedge T_\alpha \wedge T_\beta} e^{-\Lambda _u}
\, dA _u^{\Lop F} \right] . \label{lim-A}
\een
The continuity of $F$ and (\ref{F(a,b)}) imply that there
exists a constant $C_1 > 0$ such that
\be
|F(y)| \leq C_1 \left[ \varphi (y) + \psi (y) \right] .
\ee
Also, (\ref{phi-psi-lms}) implies that the processes $\left(
e^{-\Lambda_t} \varphi (X_t) \right)$ and $\left( e^{-\Lambda_t}
\psi (X_t) \right)$ are positive supermartingales, therefore,
\be
\e_x \left[ e^{-\Lambda_\tau} \left[ \varphi (X_\tau) + \psi (X_\tau)
\right] {\bf 1} _{\{ \tau < \infty \}} \right] \leq C_1 \left[
\varphi (x) + \psi (x) \right]< \infty .
\ee
Since
\be
e^{-\Lambda _\tau} \left| F(X_\tau) \right| {\bf 1} _{\{ \tau \leq
T_{\alpha_m} \wedge T_{\beta_n} \}}
\leq C_1 e^{-\Lambda_\tau} \left[ \varphi (X_\tau) + \psi (X_\tau)
\right] {\bf 1} _{\{ \tau < \infty \}} \quad \text{for all }
m,n \geq 1 ,
\ee
we can see that the dominated convergence theorem implies that
\begin{gather}
\lim _{m \rightarrow \infty} \e_x \bigl[ e^{-\Lambda _\tau}
F(X_\tau) {\bf 1} _{\{ \tau \leq T_{\alpha_m} \wedge T_{\beta_n}
\}} \bigr] = \e_x \bigl[ e^{-\Lambda _\tau} F(X_\tau) {\bf 1}
_{\{ \tau < T_\alpha \} \cap \{ \tau \leq T_{\beta_n} \}} \bigr]
\nonumber \\
\intertext{and}
\lim _{m,n \rightarrow \infty} \e_x \bigl[ e^{-\Lambda _\tau}
F(X_\tau) {\bf 1} _{\{ \tau \leq T_{\alpha_m} \wedge T_{\beta_n}
\}} \bigr] = \e_x \bigl[ e^{-\Lambda _\tau} F(X_\tau) {\bf 1}
_{\{ \tau < T_\alpha \wedge T_\beta \}} \bigr] . \label{S3-lim3}
\end{gather}
In view of these results, we can pass to the limit $m \rightarrow
\infty$ in (\ref{Dynkin-CFR}) to obtain
\begin{align}
F(x) & + \e_x \left[ \int _0^{\tau \wedge T_\alpha \wedge
T_{\beta_n}} e^{-\Lambda_u} \, dA _u^{\Lop F} \right]
\nonumber \\
= \mbox{} & \e_x \left[ e^{-\Lambda _\tau} F(X_\tau) {\bf 1}
_{\{ \tau < T_\alpha \} \cap \{ \tau \leq T_{\beta_n} \}} \right]
+ \lim _{m \rightarrow \infty} \varphi (x) \frac{F(\alpha_m)}
{\varphi (\alpha_m)} \frac{\e_x \bigl[ e^{-\Lambda _{T_{\alpha_m}}}
{\bf 1} _{\{ T_{\alpha_m} < \tau \wedge T_{\beta_n} \}} \bigr]}
{\e_x \bigl[ e^{-\Lambda _{T_{\alpha_m}}} \bigr]} \nonumber \\
& + \lim _{m \rightarrow \infty} \psi (x) \frac{F(\beta_n)}
{\psi (\beta_n)} \frac{\e_x \bigl[ e^{-\Lambda _{T_{\beta_n}}}
{\bf 1} _{\{ T_{\beta_n} < \tau \wedge T_{\alpha_m} \}} \bigr]}
{\e_x \bigl[ e^{-\Lambda _{T_{\beta_n}}} \bigr]} \nonumber \\
= \mbox{} & \e_x \left[ e^{-\Lambda _\tau} F(X_\tau) {\bf 1}
_{\{ \tau < T_\alpha \} \cap \{ \tau \leq T_{\beta_n} \}} \right]
+ \varphi (x) \frac{\e_x \bigl[ e^{-\Lambda _{T_\alpha}} {\bf 1}
_{\{ T_\alpha \leq \tau \wedge T_{\beta_n} \}} \bigr]}
{\e_x \bigl[ e^{-\Lambda _{T_\alpha}} \bigr]} \lim _{m \rightarrow
\infty} \frac{F(\alpha_m)} {\varphi (\alpha_m)}\nonumber \\
& + \psi (x) \frac{F(\beta_n)} {\psi (\beta_n)} \frac{\e_x \bigl[
e^{-\Lambda _{T_{\beta_n}}} {\bf 1} _{\{ T_{\beta_n} < \tau
\wedge T_\alpha \}} \bigr]} {\e_x \bigl[ e^{-\Lambda
_{T_{\beta_n}}} \bigr]} , \label{Dynkin-pr-n}
\end{align}
the second equality following by an application of the
dominated convergence theorem.
These identities prove  that the limit $\lim _{y \downarrow \alpha}
F(y) / \varphi (y)$ exists because $(\alpha_m)$ has been an
arbitrary sequence satisfying (\ref{Fab-lim}) and the function
$F/\varphi$ is continuous.

Proving that the limit  $\lim _{y \uparrow \beta} F(y) / \psi (y)$
exists follows similar symmetric arguments.

{\em Proof of} (III).
The event $\{ T_\alpha < \infty \}$ has strictly positive probability
if and only if $\alpha$ is an absorbing boundary point, in which
case, (\ref{phi-psi-aabs}) and (\ref{natMthm}) imply that
$\lim _{y \downarrow \alpha} F(y) = 0$.
In view of this observation and a similar one concerning the
boundary point $\beta$, we can see that the first identity
in (\ref{Dynkin}) holds true.
Finally, the second identity in (\ref{Dynkin}) follows immediately
once we combine (\ref{natMthm}) with (\ref{exp-lims}) and
(\ref{lim-A})--(\ref{Dynkin-pr-n}).
\mbox{}\hfill$\Box$
\vspace{5mm}

The assumptions of the previous lemma involve the measure
$\Lop F$ that we can associate with a function on $\intI$ that
is the difference of two convex functions.
We now address the following inverse problem: given a
signed measure $\mu$ on $\bigl( \intI, {\mathcal B} (\intI) \bigr)$,
determine a function $F$ on $\intI$ such that $F$ is the
difference of two convex functions and $\Lop F = - \mu$.
Plainly, the solution to this problem is not unique because
$\Lop \varphi = \Lop \psi = 0$.
In view of this observation, the solution $R_\mu$ that we
now derive and identifies with the $r(\cdot)$-potential of the
continuous additive functional $A^\mu$ is ``minimal'' in the
sense that it has the limiting behaviour captured by
(\ref{Rab}).

\begin{thm} \label{prop:ODE-special}
A signed Radon measure $\mu$ on $\bigl( \intI, {\mathcal B}
(\intI) \bigr)$ satisfies
\ben
\int _{]\alpha, x[} \frac{\psi (s)} {\sigma^2 (s) p' (s)} \, |\mu| (ds)
+ \int _{[x, \beta[} \frac{\varphi (s)} {\sigma^2 (s) p' (s)} \, |\mu|
(ds) < \infty \label{mu-intA}
\een
for all $x \in \CI$, if and only if
\ben
\int _{\bar{\alpha}}^{\bar{\beta}} \frac{1}{\sigma^2 (s)} \, |\mu|
(ds) < \infty \quad \text{and} \quad
\e_x \left[ \int _0^{T_\alpha \wedge T_\beta} e^{-\Lambda_t} \,
dA_t^{|\mu|} \right] < \infty \label{mu-intP}
\een
for all $\alpha < \bar{\alpha} < \bar{\beta} < \beta$ and all
$x \in \CI$.
In the presence of these integrability conditions, the function
$R_\mu : \intI \rightarrow \R$ defined by
\ben
R_\mu (x) = \frac{2}{C} \varphi (x) \int _{]\alpha, x[} \frac{\psi (s)}
{\sigma^2 (s) p' (s)} \, \mu (ds) + \frac{2}{C} \psi (x) \int
_{[x, \beta[} \frac{\varphi (s)} {\sigma^2 (s) p' (s)} \, \mu (ds)
, \label{Rmu-an}
\een
where $C>0$ is the constant appearing in (\ref{C}), identifies
with the $r(\cdot)$-potential of $A^\mu$, namely,
\ben
R_\mu (x) = \e_x \left[ \int _0^{T_\alpha \wedge T_\beta}
e^{-\Lambda_t} \, dA_t^{\mu} \right] , \label{Rmu-pr}
\een
it is the difference of two convex functions, and
\ben
\Lop R_\mu (dx) = - \mu(dx) \quad \text{and} \quad
\Lop R_{|\mu|} (dx) = - |\mu| (dx) . \label{LopRmu}
\een
Furthermore,
\begin{align}
\lim _{y \downarrow \alpha} \frac{|R_\mu (y)|} {\varphi (y)} =
\lim _{y \uparrow \beta} \frac{|R_\mu (y)|} {\psi (y)} =
\lim _{y \downarrow \alpha} \frac{R_{|\mu|} (y)} {\varphi (y)} = 
\lim _{y \uparrow \beta} \frac{R_{|\mu|} (y)} {\psi (y)} = 0 .
\label{Rab}
\end{align}
\end{thm}
{\bf Proof.}
First, we note that, if the integrability condition (\ref{mu-intA})
is true for some $x \in \CI$, then it is true for all $x \in \CI$.
If $\mu$ is a measure on $\bigl( \intI, {\mathcal B} (\intI) \bigr)$
satisfying (\ref{mu-intA}), then the function $R_\mu$ given by
(\ref{Rmu-an}) is well-defined, it is the difference of two convex
functions, and it satisfies the corresponding identity in
(\ref{LopRmu}).
To see these claims, we consider the left-continuous function
$H: \intI \rightarrow \R$ given by $H(\gamma) = 0$ and
\be
H(x) = \begin{cases} - \int _{]x,\gamma[} \frac{2}{C \sigma
^2 (s) p'(s)} \, \mu (ds) , & \text{if } x \in \mbox{} ]\alpha ,
\gamma[ , \\ \int _{[\gamma , x[} \frac{2}{C \sigma ^2 (s) p'(s)}
\, \mu (ds) , & \text{if } x \in \mbox{} ]\gamma , \beta[ ,
\end{cases}
\ee
where $\gamma$ is any constant in $\intI$.
Given any points $\bar{\alpha} , \bar{\beta} \in \intI$ such that
$\bar{\alpha} < \gamma < \bar{\beta}$, we can use the integration
by parts formula to see that
\begin{gather}
- H (\bar{\alpha}) \psi (\bar{\alpha}) - \int _{\bar{\alpha}}^x \psi
' (s) H(s) \, ds = - H(x) \psi (x) + \int _{[\bar{\alpha} ,x[} \frac{2
\psi (s)}{C \sigma ^2 (s) p'(s)} \, \mu (ds) , \nonumber \\
H(\bar{\beta}) \varphi (\bar{\beta}) - \int _x^{\bar{\beta}}
\varphi '(s) H(s) \, ds = H(x) \varphi (x) + \int _{[x, \bar{\beta}[}
\frac{2 \varphi (s)} {C \sigma ^2 (s) p'(s)} \, \mu (ds) \nonumber
\end{gather}
for all $x \in [\bar{\alpha}, \bar{\beta}]$.
It follows that the function $R_\mu$ defined by (\ref{Rmu-an})
admits the expression
\begin{align}
R_\mu (x) = \mbox{} & \left[ \frac{2}{C} \int _{]\alpha ,
\bar{\alpha}[} \frac{\psi (s)}{\sigma ^2 (s) p'(s)} \, \mu (ds)
- H(\bar{\alpha}) \psi (\bar{\alpha}) \right] \varphi (x)
\nonumber \\
& + \left[ \frac{2}{C} \int _{[\bar{\beta} , \beta [}
\frac{\varphi (s)}{\sigma ^2 (s) p'(s)} \, \mu (ds) +
H (\bar{\beta}) \varphi (\bar{\beta}) \right] \psi (x) \nonumber \\
& - \varphi (x) \int _{\bar{\alpha}}^x \psi '(s) H(s) \, ds - \psi (x)
\int _x^{\bar{\beta}} \varphi '(s) H(s) \, ds \label{TH}
\end{align}
for all $\alpha < \bar{\alpha} \leq x \leq \bar{\beta} < \beta$. 
This result, the left-continuity of $H$ and (\ref{C})
imply that
\begin{align}
(R_\mu) _-' (x) = \mbox{} & \left[ \frac{2}{C} \int _{]\alpha ,
\bar{\alpha}[} \frac{\psi (s)} {\sigma ^2 (s) p'(s)} \, \mu (ds)
- H(\bar{\alpha}) \psi (\bar{\alpha}) \right] \varphi ' (x)
\nonumber \\
& + \left[ \frac{2}{C} \int _{[\bar{\beta} , \beta [} \frac{\varphi (s)}
{\sigma ^2 (s) p'(s)} \, \mu (ds) + H(\bar{\beta}) \varphi
(\bar{\beta}) \right] \psi '(x) - C p' (x) H(x) \nonumber \\
& - \varphi '(x) \int _{\bar{\alpha}}^x \psi '(s) H(s) \, ds
- \psi ' (x) \int _x^{\bar{\beta}} \varphi '(s) H(s) \, ds
\label{TH'}
\end{align}
for all $\alpha < \bar{\alpha} \leq x \leq \bar{\beta} < \beta$. 
Furthermore, we can see that the restriction of the measure
$(R_\mu)''$ in $\bigl( ]\bar{\alpha}, \bar{\beta}[, {\mathcal B}
(]\bar{\alpha}, \bar{\beta}[) \bigr)$ has Lebesgue decomposition
that is given by
\begin{align}
(R_\mu) _\ac '' (x) = \mbox{} & \left[ \frac{2}{C} \int _{]\alpha,
\bar{\alpha}[} \frac{\psi (s)} {\sigma ^2 (s) p'(s)} \, \mu (ds)
- H(\bar{\alpha}) \psi (\bar{\alpha}) \right] \varphi '' (x)
\nonumber \\
& + \left[ \frac{2}{C} \int _{[\bar{\beta} , \beta [}
\frac{\varphi (s)} {\sigma ^2 (s) p'(s)} \, \mu (ds) +
H(\bar{\beta}) \varphi (\bar{\beta}) \right] \psi ''(x) - C p'' (x)
H(x) - \frac{2 \mu_\ac (x)} {\sigma ^2 (x)} \nonumber \\
& - \varphi ''(x) \int _{\bar{\alpha}}^x \psi '(s) H(s) \, ds - \psi
''(x) \int _x^{\bar{\beta}} \varphi '(s) H(s) \, ds , \nonumber \\
(R_\mu) _\s '' (dx) = \mbox{} & - \frac{2}{\sigma ^2 (x)} \,
\mu _\s (dx) , \nonumber
\end{align}
in the notation of (\ref{Leb-dec}).
Combining these expressions with (\ref{TH})--(\ref{TH'})
and the definition (\ref{scale}) of the scale function $p$,
we can see that the restrictions of the measures $\Lop
R_\mu$ and $-\mu$ in $\bigl( ]\bar{\alpha}, \bar{\beta}[,
{\mathcal B} (]\bar{\alpha}, \bar{\beta}[) \bigr)$ are equal.
It follows that the measures $\Lop R_\mu$ and $-\mu$ on
$\bigl( \intI, {\mathcal B} (\intI) \bigr)$ are equal because
$\bar{\alpha} < \bar{\beta}$ have been arbitrary points
in $\intI$.
Similarly, we can check that the function $R_{|\mu|}$ that is
defined by (\ref{Rmu-an}) with $|\mu|$ in place of $\mu$ is the
difference of two convex functions and satisfies the
corresponding identity in (\ref{LopRmu}).

To proceed further, we consider any Radon measure $\mu$
on $\bigl( \intI, {\mathcal B} (\intI) \bigr)$.
Given monotone sequences $(\alpha_n)$ and $(\beta_n)$
as in (\ref{ab-seq}), we define
\be
\zeta _n (z) = \begin{cases} 0 , & \text{if } z < \alpha_n
\text{ or } z > \beta_n , \\ 1 , & \text{if } \sigma^2 (z) \geq
\frac{1}{n} \text{ and } \alpha_n \leq z \leq \beta_n , \\
\sigma^2 (z) , & \text{if } \sigma^2 (z) < \frac{1}{n} \text{ and }
\alpha_n \leq z \leq \beta_n , \end{cases}
\ee
and we consider the sequence of measures $(\mu_n)$ that
are defined by (\ref{mun-defn}).
The functions $R_{|\mu_n|}$, defined by (\ref{Rmu-an}) with
$|\mu_n|$ in place of $\mu$, are real-valued and satisfy
\be
R_{|\mu_n|} (x) = \begin{cases} \frac{2}{C} \psi (x) \int
_{[\alpha_n, \beta_n]} \frac{\varphi (s)} {\sigma^2 (s) p' (s)}
\zeta_n (s) \, |\mu| (ds) , & \text{if } x < \alpha_n , \\
\frac{2}{C} \varphi (x) \int _{[\alpha_n, \beta_n]} \frac{\psi (s)}
{\sigma^2 (s) p' (s)} \zeta_n (s) \, |\mu| (ds) , & \text{if }
x > \beta_n . \end{cases}
\ee
Combining this calculation with (\ref{phi-psi->0}), we can see
that $R_{|\mu_n|}$ satisfies the corresponding limits in
(\ref{Rab}).
Since $-\Lop R_{|\mu_n|} = |\mu_n| = |\Lop R_{\mu_n}|$ is a
positive measure, part~(I) of Lemma~\ref{prop:Dynkin}
implies that
\be
\e_x \left[ \int _0^{T_\alpha \wedge T_\beta} e^{-\Lambda_t}
\, dA_t^{|\mu_n|} \right] < \infty \quad \text{for all } x \in \CI ,
\ee
while, (\ref{Dynkin}) in Lemma~\ref{prop:Dynkin} with
$\tau = T_\alpha \wedge T_\beta$ implies that
\be
R_{|\mu_n|} (x) = - \e_x \left[ \int _0^{T_\alpha \wedge
T_\beta} e^{-\Lambda _t} \, dA_t^{\Lop R_{|\mu_n|}} \right] .
\ee
This identity, the fact that $\Lop R_{|\mu_n|} = - |\mu_n|$
and (\ref{A-lin}) imply that the function $R_{|\mu_n|}$
that is defined as in (\ref{Rmu-an}) satisfies
\ben
R_{|\mu_n|} (x) = \e_x \left[ \int _0^{T_\alpha \wedge T_\beta}
e^{-\Lambda _t} \, dA_t^{|\mu_n|} \right] . \label{Rmu-CS}
\een

Since the sequence of functions $(\zeta_n)$ is monotonically
increasing to the identity function, the monotone convergence
theorem implies that
\begin{align}
R_{|\mu|} (x) & = \lim _{n \rightarrow \infty} \biggl( \frac{2}{C}
\varphi (x) \int _{]\alpha, x[} \frac{\psi (s)} {\sigma^2 (s)
p' (s)} \zeta_n (s) \, |\mu| (ds) \nonumber \\
& \mbox{} \hspace{17mm} + \frac{2}{C} \psi (x) \int
_{[x, \beta[} \frac{\varphi (s)} {\sigma^2 (s) p' (s)} \zeta_n (s)
\, |\mu| (ds) \biggr) . \label{Rmun-lim/an}
\end{align}
If (\ref{mu-intA}) is satisfied, then $\sigma^{-2}$ is locally
integrable with respect to $|\mu|$, namely, the first condition
in (\ref{mu-intP}) holds true, thanks to the continuity of the
functions $\varphi$, $\psi$ and $p'$.
In this case, (\ref{A-mun}) in Lemma~\ref{lem:Amu} and
(\ref{Rmu-CS}) imply that
\ben
\lim _{n \rightarrow \infty} R_{|\mu_n|} (x)
= \lim _{n \rightarrow \infty} \e_x \left[ \int _0^{T_\alpha \wedge
T_\beta} e^{-\Lambda_t} \, dA_t^{|\mu_n|} \right]
= \e_x \left[ \int _0^{T_\alpha \wedge T_\beta} e^{-\Lambda_t}
\, dA_t^{|\mu|} \right] \label{Rmun-lim/pr}
\een
because $(\zeta_n)$ satisfies (\ref{zeta}).
Combining this result with (\ref{Rmun-lim/an}) and the fact that
(\ref{mu-intA}) implies that $R_{|\mu|} (x) < \infty$, we can
see that
\ben
R_{|\mu|} (x) = \e_x \left[ \int _0^{T_\alpha \wedge T_\beta}
e^{-\Lambda_t} \, dA_t^{|\mu|} \right] < \infty , \label{Rmu-an/pr}
\een
and the second condition in (\ref{mu-intP}) follows.
Thus, we have proved that (\ref{mu-intA}) implies 
(\ref{mu-intP}).
Conversely, if (\ref{mu-intP}) is satisfied, then (\ref{A-mun})
in Lemma~\ref{lem:Amu} and (\ref{Rmu-CS}) imply that
(\ref{Rmun-lim/pr}) is true.
Combining (\ref{mu-intP}) with (\ref{Rmun-lim/an}) and
(\ref{Rmun-lim/pr}), we can see that $R_{|\mu|} (x) < \infty$,
and (\ref{mu-intA}) follows.

If $\mu$ satisfies the integrability conditions
(\ref{mu-intA})--(\ref{mu-intP}), then the function $R_\mu$ given by
(\ref{Rmu-an}) is well-defined and real-valued.
Furthermore, it satisfies (\ref{Rmu-pr}) thanks to (\ref{A-lin}),
(\ref{Rmu-an/pr}) with $\mu^+$ and $\mu^-$ in place of $|\mu|$,
and the linearity of integrals.

To establish (\ref{Rab}), we consider any sequences
$(\alpha_n)$, $(\beta_n)$ as in (\ref{ab-seq}), and we calculate
\begin{align}
0 & \stackrel{(\ref{Rmu-pr})}{=} R_{|\mu|} (x) - \lim
_{m, n \rightarrow \infty} \e_x \left[ \int _0 ^{T_{\alpha_m}
\wedge T_{\beta_n}} e^{-\Lambda_u} \, dA _u^{|\mu|} \right]
\nonumber \\
& \stackrel{(\ref{LopRmu})}{=} R_{|\mu|} (x) + \lim
_{m, n \rightarrow \infty} \e_x \left[ \int _0^{T_{\alpha_m}
\wedge T_{\beta_n}} e^{-\Lambda_u} \, dA _u ^{\Lop R_{|\mu|}}
\right] \nonumber \\
& = \lim _{m, n \rightarrow \infty} \e_x \left[ e^{-\Lambda
_{T_{\alpha_m} \wedge T_{\beta_n}}} R_{|\mu|}
(X_{T_{\alpha_m} \wedge T_{\beta_n}}) \right]
\nonumber \\
& = \lim _{m, n \rightarrow \infty} R_{|\mu|} (\alpha_m)
\e_x \left[ e^{-\Lambda _{T_{\alpha_m}}} {\bf 1}
_{\{ T_{\alpha_m} < T_{\beta_n} \}} \right] + \lim
_{m, n \rightarrow \infty} R_{|\mu|} (\beta_n) \e_x \left[
e^{-\Lambda _{T_{\beta_n}}} {\bf 1} _{\{ T_{\beta_n} <
T_{\alpha_m} \}} \right] , \nonumber
\end{align}
the third identity following from (\ref{Dynkin}) for
$\tau = T_{\alpha_m} \wedge T_{\beta_n}$.
Since $R_{|\mu|}$ is a positive function, each of the two limits on
the right-hand side of this expression is equal to 0.
We can therefore see that the first of these limits
implies that
\begin{align}
0 & = \lim _{m \rightarrow \infty} \lim _{n \rightarrow \infty}
R_{|\mu|} (\alpha_m) \e_x \left[ e^{-\Lambda _{T_{\alpha_m}}}
{\bf 1} _{\{ T_{\alpha_m} < T_{\beta_n} \}} \right] \nonumber \\
& = \lim _{m \rightarrow \infty} R_{|\mu|} (\alpha_m)
\e_x \left[ e^{-\Lambda _{T_{\alpha_m}}} {\bf 1} _{\{
T_{\alpha_m} < T_{\beta} \}} \right] \nonumber \\
& \stackrel{(\ref{phi})}{=} \lim _{m \rightarrow \infty}
\frac{R_{|\mu|} (\alpha_m) \varphi (x)} {\varphi (\alpha_m)} ,
\nonumber
\end{align}
which proves that $\lim _{y \downarrow \alpha} R_{|\mu|}
(y) / \varphi (y) = 0$ because $(\alpha_m)$ has been arbitrary.
We can show that $\lim _{x \uparrow \beta} R_{|\mu|} (x) / \psi
(x) = 0$ using similar arguments.
Finally, the function $|R_\mu|$ satisfies the corresponding limits in
(\ref{Rab}) because $|R_\mu| \leq
R_{|\mu|}$.
\mbox{} \hfill$\Box$
\vspace{5mm}

The result we have just established and
Lemma~\ref{prop:Dynkin} imply the following representation
of differences of two convex functions that involves the
operator $\Lop$ and the functions $\varphi$, $\psi$.

\begin{thm} \label{cor:F-repres}
Consider any function $F: \intI \rightarrow \R$ that is the difference
of two convex functions, and suppose that
\be
\limsup _{y \downarrow \alpha} \frac{|F(y)|}{\varphi (y)}
< \infty \quad \text{and} \quad \limsup _{y \uparrow \beta}
\frac{|F(y)|}{\psi (y)} < \infty ,
\ee
and that the measure $\Lop F$ satisfies the equivalent
integrability conditions (\ref{mu-intA})--(\ref{mu-intP}) (see
also Remark~\ref{ACremRLF} below).
In this context, the limits $\lim _{y \downarrow \alpha}
F(y) / \varphi (y)$ and $\lim _{y \uparrow \beta} F(y) / \psi (y)$
both exist, and the function $F$ admits the representation
\ben
F(x) = \lim _{y \downarrow \alpha} \frac{F(y)}{\varphi (y)}
\varphi (x) + R_{-\Lop F} (x) + \lim _{y \uparrow \beta}
\frac{F(y)}{\psi (y)} \psi (x) , \label{F-repr}
\een
where $R_{-\Lop F}$ is given by (\ref{Rmu-an})--(\ref{Rmu-pr}).
Furthermore, given any points $\bar{\alpha} < x < \bar{\beta}$ in
$\CI$ and any stopping time $\tau$,
\ben
\e_x \left[ e^{-\Lambda _{\tau \wedge T_{\bar{\alpha}} \wedge
T_{\bar{\beta}}}} F(X_{\tau \wedge T_{\bar{\alpha}}
\wedge T_{\bar{\beta}}}) \right] = F(x) + \e_x \left[ \int
_0^{\tau \wedge T_{\bar{\alpha}} \wedge T_{\bar{\beta}}}
e^{-\Lambda_u} \, dA_u^{\Lop F} \right] ,
\label{Dynkin-genF}
\een
in which expression, we denote
\be
F(\alpha) = \lim _{y \downarrow \alpha} F(y) \quad
\left( \text{resp., } F(\beta) = \lim _{y \uparrow \beta} F(y)
\right)
\ee
if $\alpha$ (resp., $\beta$) is absorbing, namely, if $\p_x (T_\alpha
< \infty) > 0$ (resp., $\p_x (T_\beta < \infty) > 0$).
\end{thm}
{\bf Proof.}
In the presence of the assumption that $\Lop F$ satisfies
(\ref{mu-intA})--(\ref{mu-intP}), Lemma~\ref{prop:Dynkin}.(II)
implies that the limits $\lim _{y \downarrow \alpha} F(y) / \varphi (y)$
and $\lim _{y \uparrow \beta} F(y) / \psi (y)$ exist, while
Theorem~\ref{prop:ODE-special} implies that the function
$R_{-\Lop F}$ is well-defined.
In particular, (\ref{LopRmu}) implies that $\Lop \left( F -
R_{-\Lop F} \right) = 0$.
It follows that
\be
F - R_{-\Lop F} = A \varphi + B \psi ,
\ee
for some constants $A, B \in \R$.
Combining (\ref{phi-psi->0}) with (\ref{Rab}), we can see that
the constants $A$ and $B$ are as in (\ref{F-repr}).
Finally, (\ref{Dynkin-genF}) follows from the representation
(\ref{F-repr}) of $F$, (\ref{phi-psi-expect}) in
Lemma~\ref{lem:phi-psi-expect}, (\ref{Dynkin}) in
Lemma~\ref{prop:Dynkin} and (\ref{Rab}) in
Theorem~\ref{prop:ODE-special}.
\mbox{}\hfill$\Box$

\begin{rem} \label{ACremRLF}
{\rm
In view of Lemma~\ref{prop:Dynkin}.(I), the positivity
of the measure $-\Lop F$ is a sufficient condition for
$\Lop F$ to satisfy the integrability conditions
(\ref{mu-intA})--(\ref{mu-intP}).
Also, if $F$ is $C^1$ with first derivative that is absolutely
continuous with respect to the Lebesgue measure, then
$\Lop F (dx) = \Lop_\ac F (x) \, dx$, where $\Lop_\ac$ is
defined by (\ref{lopac}).
This observation and part~(II) of Lemma~\ref{lem:Ito}
imply that, in this case (\ref{mu-intA})--(\ref{mu-intP}) are
equivalent to (\ref{R|h|-ana})--(\ref{R|h|-prob}) below
for $h = \Lop_\ac F$.
Furthermore, $R_{-\Lop F}$ admits the expressions
(\ref{Rh-a})--(\ref{Rh-p}) below for $h = - \Lop_\ac F$.
} \mbox{} \hfill $\Box$
\end{rem}

The measure $\Lop F$ and the potential $R_{-\Lop F}$
have central roles in the characterisation of differences
of two convex functions we have established above.
The following result is concerned with the potential
$R_{-\Lop F}$ when $\Lop F$ is absolutely continuous
with respect to the Lebesgue measure.

\begin{cor} \label{cor:R_mu_h}
Consider any function $h: \CI \rightarrow \R$ that is locally
integrable with respect to the Lebesgue measure, and let
$\mu^h$ be the measure on $\bigl( \intI, {\mathcal B} (\intI)
\bigr)$ defined by
\be
\mu^h (\Gamma) = \int _\Gamma h(s) \, ds , \quad \text{for }
\Gamma \in {\mathcal B} (\intI) .
\ee
If $\mu^h$ satisfies the equivalent integrability conditions
(\ref{mu-intA})--(\ref{mu-intP}), which are equivalent to
\begin{gather}
\int _\alpha^x \frac{\psi (s)} {\sigma^2 (s) p' (s)} |h(s)| \, ds +
\int _x^\beta \frac{\varphi (s)} {\sigma^2 (s) p' (s)} |h(s)| \, ds
< \infty , \label{R|h|-ana} \\
\e_x \left[ \int _0^{T_\alpha \wedge T_\beta} e^{-\Lambda _t}
\bigl| h(X_t) \bigr| \, dt \right] < \infty , \label{R|h|-prob}
\end{gather}
then the function $R_{\mu^h}: \intI \rightarrow \R$ defined by
(\ref{Rmu-an}) or, equivalently, by
\ben
R_{\mu^h} (x) = \frac{2}{C} \varphi (x) \int _\alpha^x
\frac{\psi (s)} {\sigma^2 (s) p' (s)} h(s) \, ds + \frac{2}{C}
\psi (x) \int _x^\beta \frac{\varphi (s)} {\sigma^2 (s) p' (s)}
h(s) \, ds , \label{Rh-a}
\een
admits the probabilistic expression
\ben
R_{\mu^h} (x) = \e_x \left[ \int _0^{T_\alpha \wedge T_\beta}
e^{-\Lambda _t} h(X_t) \, dt \right] . \label{Rh-p}
\een
This function, as well as the function defined by
\ben
\tilde{R}_{\mu^h} (x) = \e_x \left[ \int _0^\infty e^{-\Lambda _t}
h(X_t) \, dt \right] , \quad \text{for } x \in \intI ,
\een
is $C^1$ with absolutely continuous first derivative
and satisfies the ODE
\ben
\Lop_\ac g (x) + h(x) \equiv \half \sigma ^2 (x) g'' (x) +
b(x) g' (x) - r(x) g(x) + h(x) = 0 . \label{ODE}
\een
The functions $R_{\mu^h}$ and $\tilde{R}_{\mu^h}$ satisfy
\begin{gather}
\tilde{R}_{\mu^h} (x) = \frac{h(\alpha)}{r(\alpha)} \frac{\varphi
(x)}{\varphi (\alpha)} {\mathbb I}_\alpha + R_{\mu^h} (x) +
\frac{h(\beta)}{r(\beta)} \frac{\psi (x)}{\psi (\beta)} {\mathbb I}
_\beta , \label{Rmuh-tilde} \\
\lim _{y \downarrow \alpha} \frac{\tilde{R}_{\mu^h} (y)}
{\varphi (y)} = \lim _{y \downarrow \alpha} \frac{R_{\mu^h}
(y)} {\varphi (y)} + \frac{h(\alpha)} {r(\alpha) \varphi (\alpha)}
{\mathbb I}_\alpha = \frac{h(\alpha)} {r(\alpha) \varphi (\alpha)}
{\mathbb I}_\alpha , \label{Rmuh-lima} \\
\lim _{y \uparrow \beta} \frac{\tilde{R}_{\mu^h} (y)} {\psi (y)}
= \lim _{y \uparrow \beta} \frac{R_{\mu^h} (y)} {\psi (y)} +
\frac{h(\beta)}{r(\beta) \psi (\beta)} {\mathbb I}_\beta =
\frac{h(\beta)} {r(\beta) \psi (\beta)} {\mathbb I}_\beta ,
\label{Rmuh-limb} \\
\intertext{where}
{\mathbb I} _\alpha = \begin{cases} 1 , & \text{if } \alpha
\text{ is absorbing} , \\ 0 , & \text{if } \alpha \text{ is inaccessible}
, \end{cases} \quad \text{and} \quad {\mathbb I} _\beta
= \begin{cases} 1 , & \text{if } \beta \text{ is absorbing} , \\
0 , & \text{if } \beta \text{ is inaccessible} . \end{cases}
\nonumber
\end{gather}
Furthermore,
\ben
R_{\mu^h} (x) = \e_x \left[ \int _0^{\tau \wedge T_\alpha \wedge
T_\beta} e^{-\Lambda_t} h(X_t) \, dt + e^{-\Lambda _{\tau
\wedge T_\alpha \wedge T_\beta}} R_{\mu^h} (X_{\tau \wedge
T_\alpha \wedge T_\beta}) {\bf 1} _{\{ \tau \wedge T_\alpha
\wedge T_\beta < \infty \}} \right] \label{Rmuh-Dynkin}
\een
for every stopping time $\tau$ and all initial conditions $x \in
\intI$, in which expression, $R_{\mu^h} (\alpha) = 0$ (resp.,
$R_{\mu^h} (\beta) = 0$) if $\alpha$ (resp., $\beta$) is absorbing,
consistently with (\ref{Rmuh-lima})--(\ref{Rmuh-limb}).
\end{cor}
{\bf Proof.}
It is straightforward to check that the function $R_{\mu^h}$
defined by (\ref{Rh-a}) is $C^1$ with absolutely continuous
first derivative and satisfies the ODE (\ref{ODE}).
This observation and (\ref{Rmuh-tilde}) imply the corresponding
statements for $\tilde{R}_{\mu^h}$.
The equivalence of (\ref{mu-intP}) (resp., (\ref{Rmu-pr}))
with (\ref{R|h|-prob}) (resp., (\ref{Rh-p})) is a consequence
of part~(II) of Lemma~\ref{lem:Ito} and the identities
$\mu^h (dx) = -\Lop R_{\mu^h} (dx) = - \Lop_\ac R_{\mu^h}
(x) \, dx = h(x) \, dx$.
Also, these identities, part~(II) of Lemma~\ref{lem:Ito} and
(\ref{Dynkin}) imply (\ref{Rmuh-Dynkin}), while the limits in
(\ref{Rmuh-lima})--(\ref{Rmuh-limb}) follow from (\ref{Rab})
and (\ref{Rmuh-tilde}).

To prove (\ref{Rmuh-tilde}), we first note that
\be
\tilde{R}_{\mu^h} (x) = \e_x \left[ {\bf 1} _{\{ T_\alpha < T_\beta
\}} \int _{T_\alpha}^\infty e^{-\Lambda _t} \, dt \right] h(\alpha)
+ R_{\mu^h} (x) + \e_x \left[ {\bf 1} _{\{ T_\beta < T_\alpha \}}
\int _{T_\beta}^\infty e^{-\Lambda _t} \, dt \right] h(\beta) .
\ee
In view of the definition (\ref{Lambda}) of $\Lambda$, we can see
that, if $\alpha$ is absorbing, then
\begin{align}
\e_x \left[ {\bf 1} _{\{ T_\alpha < T_\beta \}} \int _{T_\alpha}
^\infty e^{-\Lambda _t} \, dt \right] & = \e_x \left[ {\bf 1}
_{\{ T_\alpha < T_\beta \}} e^{-\Lambda _{T_\alpha}} \int _{T_\alpha}
^\infty e^{- r(\alpha) (t-T_\alpha)} \, dt \right] \nonumber \\
& = \frac{1}{r(\alpha)} \e_x \left[ e^{-\Lambda _{T_\alpha}} {\bf 1}
_{\{ T_\alpha < T_\beta \}} \right] \nonumber \\
& \stackrel{(\ref{phi})}{=} \frac{1}{r(\alpha)} \frac{\varphi (x)}
{\varphi (\alpha)} , \nonumber
\end{align}
otherwise, this expectation is plainly 0.
Similarly, we can see that
\be
\e_x \left[ {\bf 1} _{\{ T_\beta < T_\alpha \}} \int _{T_\beta}^\infty
e^{-\Lambda _t} \, dt \right] = \frac{1}{r(\beta)} \frac{\psi (x)}
{\psi (\beta)} {\mathbb I}_\beta ,
\ee
and (\ref{Rmuh-tilde}) follows.
\mbox{}\hfill$\Box$

\section{Analytic characterisations of $\pmb{r(\cdot)}$-excessive
functions}
\label{sec:excessive}

The following is the main result of this section.

\begin{thm} \label{prop:excessive}
A function $F: \CI \rightarrow \R_+$ is $r(\cdot)$-excessive,
namely, it satisfies
\ben
\e_x \left[ e^{-\Lambda _\tau} F(X_\tau) {\bf 1} _{\{\tau < \infty \}}
\right] \leq F(x) \label{F-exc}
\een
for all stopping times $\tau$ and all initial conditions $x \in \CI$,
if and only if the following statements are both true:

\noindent
{\rm (I)}
the restriction of $F$ in the interior $\intI$ of $\CI$  is the
difference of two convex functions and the associated measure
$-\Lop F$ on $\bigl( \intI, {\mathcal B} (\intI) \bigr)$ is positive;

\noindent
{\rm (II)}
if $\alpha$ (resp., $\beta$) is an absorbing boundary point,
then $F(\alpha) \leq \liminf _{y \in \intI , \, y \downarrow \alpha}
F(y)$ (resp., $F(\beta) \leq \liminf _{y \in \intI , \, y \uparrow \beta}
F(y)$).
\end{thm}
{\bf Proof.}
First, we consider any function $F: \CI \rightarrow \R_+$
with the properties listed in (I)--(II).
The assumption that $-\Lop F$ is a positive measure
implies that $- A^{\Lop F} = A^{- \Lop F}$ is an increasing
process.
Therefore, (\ref{Dynkin-genF}) in Theorem~\ref{cor:F-repres}
implies that, given any points $\bar{\alpha} < x < \bar{\beta}$
in $\CI$  and any stopping time $\tau$ such that $\bar{\alpha}
= \alpha$ and $\tau = \tau \wedge T_\alpha$ (resp., $\bar{\beta}
= \beta$ and $\tau = \tau \wedge T_\beta$) if $\alpha$ (resp.,
$\beta$) is absorbing,
\begin{align}
F(x) \geq \mbox{} & \e_x \left[ e^{-\Lambda _\tau} F(X_\tau)
{\bf 1} _{\{ \tau < T_{\bar{\alpha}} \wedge T_{\bar{\beta}} \}}
\right] + F(\bar{\alpha}) \e_x \left[ e^{-\Lambda
_{T_{\bar{\alpha}}}} {\bf 1} _{\{ T_{\bar{\alpha}} \leq \tau
\wedge T_{\bar{\beta}} \}} \right] (1 - {\mathbb I} _\alpha)
\nonumber \\
& + \lim _{y \downarrow \alpha} F(y) \e_x \left[ e^{-\Lambda
_{T_{\bar{\alpha}}}} {\bf 1} _{\{ T_{\bar{\alpha}} \leq \tau \wedge
T_{\bar{\beta}} \}} \right] {\mathbb I} _\alpha  + F(\bar{\beta})
\e_x \left[ e^{-\Lambda _{T_{\bar{\beta}}}} {\bf 1} _{\{
T_{\bar{\beta}} \leq \tau \wedge T_{\bar{\alpha}} \}} \right]
(1 - {\mathbb I} _\beta) \nonumber \\
& + \lim _{y \uparrow \beta} F(y) \e_x \left[ e^{-\Lambda
_{T_{\bar{\beta}}}} {\bf 1} _{\{ T_{\bar{\beta}} \leq \tau
\wedge T_{\bar{\alpha}} \}} \right] {\mathbb I} _\beta
\nonumber \\
\geq \mbox{} & \e_x \left[ e^{-\Lambda _{\tau \wedge
T_{\bar{\alpha}} \wedge T_{\bar{\beta}}}} F(X_{\tau
\wedge T_{\bar{\alpha}} \wedge T_{\bar{\beta}}}) \right] ,
\end{align}
the second inequality following from the assumption that
$F$ satisfies the inequalities in (II).
If $\alpha$ (resp., $\beta$) is inaccessible, then we can pass
to the limit $\bar{\alpha} \downarrow \alpha$ (resp., $\bar{\beta}
\uparrow \beta$) using Fatou's lemma to obtain (\ref{F-exc})
thanks to the choices of $\bar{\alpha}$ and $\bar{\beta}$
that we have made.
It follows that $F$ is $r(\cdot)$-excessive.

To establish the reverse implication, we first show that an
$r(\cdot)$-excessive function is lower semicontinuous and
its restriction in $\intI$ is continuous.
Given an initial condition $x \in \intI$ and a point $y \in \CI$,
we can use (\ref{F-exc}) to calculate
\be
F(x) \geq \e_x \left[ e^{-\Lambda _{T_y}} \right] F(y)
\stackrel{(\ref{phi})-(\ref{psi})}{=} \min \left\{ \frac{\psi (x)}
{\psi (y)} , \frac{\varphi (x)} {\varphi (y)} \right\} F(y) .
\ee
This calculation and the continuity of the functions $\varphi$,
$\psi$ imply that $F(x) \geq \limsup _{y \rightarrow x} F(y)$,
which proves that $F$ is upper semicontinuous in $\intI$.
The same arguments but with points $x \in \CI$ and
$y \in \intI$ and their roles reversed imply that
\be
F(y) \geq \min \left\{ \frac{\psi (y)}{\psi (x)} , \frac{\varphi
(y)}{\varphi (x)} \right\} F(x) .
\ee
It follows that $F(x) \leq \liminf _{y \in \intI , \, y \rightarrow x}
F(y)$, and the lower semicontinuity of $F$ in $\CI$ has been
established.
In particular, part~(II) of the proposition is true.

To prove that an $r(\cdot)$-excessive function satisfies (I),
we define the function $F_q$ by
\ben
F_q (x) = q \e_x \left[ \int _0^\infty e^{-qt - \Lambda_t} F(X_t)
\, dt \right] , \quad \text{for } x \in \CI , \label{vq1}
\een
where $q>0$ is a constant, and we note that
\ben
0 \leq F_q (x) \stackrel{(\ref{F-exc})}{\leq} q \int _0^\infty
e^{-qt} F(x) \, dt = F(x) \quad \text{for all } x \in \CI . \label{vq2}
\een
If we consider the change of variables $u = qt$, then we can see
that
\be
F_q (x) =  \e_x \left[ \int _0^\infty e^{-u - \Lambda_{u/q}}
F(X_{u/q}) \, du \right] .
\ee
In view of (\ref{vq2}), the continuity properties of the function
$F$ and the continuity of the process $X$, this expression
implies that
\ben
\lim _{q \rightarrow \infty} F_q (x) = F(x) \quad \text{for all }
x \in \CI . \label{vq-lim}
\een

Given its definition in (\ref{vq1}), Corollary~\ref{cor:R_mu_h}
implies that the function $F_q$ is $C^1$ with absolutely
continuous first derivative and that it satisfies the ODE
\be
\half \sigma ^2 (x) F_q'' (x) + b(x) F_q'(x) - \left( q + r(x) \right)
F_q(x) + qF(x) = 0
\ee
in the interior of $\CI$.
In view of (\ref{vq2}), we can see that
\be
\half \sigma ^2 (x) F_q'' (x) + b(x) F_q'(x) - r(x) F_q(x) =
- q\left[ F(x) - F_q (x) \right] \leq 0 .
\ee
This inequality implies that
\ben
\frac{d}{dx} \left( \frac{d}{dx} \left( \frac{F_q (x)}{p' (x)}
\right) - F_q (x)\frac{d}{dx}\frac{1}{p' (x)} \right) -
\frac{2 r(x) F_q (x)}{\sigma ^2 (x) p' (x)} \leq 0 , \label{conv1}
\een
where $p$ is the scale function of the diffusion $X$, which is
defined by (\ref{scale}).

To proceed further, we introduce the antiderivatives $\A^1$ and
$\A^2$ of a function $g$ that is locally integrable in $\CI$, which
are defined by
\be
\A^1 g(x) = \int _c^x g(y) \, dy \quad \text{and} \quad
\A^2 g(x) = \int _c^x \int _c^y g(z) \, dz \, dy ,
\ee
respectively, where $c \in \CI$ is a fixed point that we can take to
be the same as the point appearing in the definition (\ref{scale}) of
the scale function $p$.
Inequality (\ref{conv1}) then implies that the function $F_q/p' -
\A^1 \left( (1/p')' F_q \right) - \A^2 \left( (2rF_q) / (\sigma ^2 p')
\right)$ is concave, which, combined with (\ref{vq-lim}), implies that
the function $G := F/p' - \A^1 \left( (1/p')' F \right) - \A^2 \left(
(2rF) / (\sigma ^2 p') \right)$ is concave.
The concavity of $G$ and the equality
\be
\frac{F}{p'} = G + \A^1 \left( \left( \frac{1}{p'} \right) ' F
\right) + \A^2 \left( \frac{2rF}{\sigma^2 p'} \right)
\ee
imply that $F/p'$ is absolutely continuous and
\be
F_- ' (x) = p' (x) \left( G_-' (x) + \A^1 \left( \frac{2rF}{\sigma^2
p'} \right) (x) \right) .
\ee
This expression shows that $F'$ has finite variation.
Furthermore, taking distributional derivatives, we can see that
\be
\frac{2}{\sigma^2 (x)} \Lop F(dx) \equiv F'' (dx) + \frac{2b (x)}
{\sigma ^2 (x)} F_-' (x) \, dx - \frac{2r (x)}{\sigma ^2 (x)} F(x)
\, dx = p' (x) G'' (dx) ,
\ee
which proves that $F$ has the properties listed in part~(I)
thanks to the concavity of $G$.
\mbox{}\hfill$\Box$
\vspace{5mm}

In the spirit of Dynkin~\cite[Theorems 15.10 and 16.4]{Dy65},
Dayanik~\cite{Da} proves that a function
$F$ is $r(\cdot)$-excessive if and only if the function $F/\varphi$
is $(\psi / \varphi)$-concave (equivalently, the function $F/\psi$
is $(- \varphi / \psi)$-concave), and he shows that such
concavity assumptions imply that the function
$- D_{\psi / \varphi}^- (F/\varphi)$ defined by (\ref{D-(F/phi)})
is increasing (equivalently, the right-continuous modification
$D_{\varphi / \psi}^+ (F/\psi)$ of the function defined by
(\ref{D-(F/psi)}) is increasing)
(see Proposition 3.1 and Remarks~3.1--3.3 of
Dayanik~\cite{Da} for the precise statements).
Such a result, which focuses on the functions
$- D_{\psi / \varphi}^- (F/\varphi)$, $D_{\varphi / \psi}^+ (F/\psi)$,
follows immediately from our analysis above.

\begin{cor} \label{cor:excessive}
A function $F: \CI \rightarrow \R_+$ is $r(\cdot)$-excessive if
and only if the following statements are both true:

\noindent
{\rm (I)}
the function $- D_{\psi / \varphi}^- (F/\varphi) $ given by
\ben
- D_{\psi / \varphi}^- (F/\varphi) (x) = - \lim _{y \uparrow x}
\frac{(F/\varphi) (x) - (F/\varphi) (y)} {(\psi / \varphi) (x) -
(\psi / \varphi) (y)} , \quad \text{for } x \in \intI ,
\label{D-(F/phi)}
\een
is well-defined, real-valued and increasing; equivalently,
the function $D_{\varphi / \psi} ^- (F/\psi)$ given by
\ben
D_{\varphi / \psi}^- (F/\psi) (x) = \lim _{y \uparrow x}
\frac{(F/\psi) (y) - (F/\psi) (x)} {(\varphi / \psi) (y) -
(\varphi / \psi) (x)} , \quad \text{for } x \in \intI ,
\label{D-(F/psi)} 
\een
is well-defined, real-valued and increasing, and

\noindent
{\rm (II)}
if $\alpha$ (resp., $\beta$) is an absorbing boundary point,
then $F(\alpha) \leq \liminf _{y \in \intI , \, y \downarrow \alpha}
F(y)$ (resp., $F(\beta) \leq \liminf _{y \in \intI , \, y \uparrow \beta}
F(y)$).
\end{cor}
{\bf Proof.}
Given a measure $\mu$ on $\bigl( \intI, {\mathcal B}
(\intI) \bigr)$, we mean that $-\mu$ is a positive measure
whenever we write $\mu (dx) \leq 0$ in the proof below.
In view of Theorem~\ref{prop:excessive}, the result will
follow if we show that either of the functions given by
(\ref{D-(F/phi)}), (\ref{D-(F/psi)}) is well-defined, real-valued
and increasing if and only if the restriction of $F$ in
$\intI$ is the difference of two convex functions and
$\Lop F \leq 0$.
To this end, we note that the functions given by
(\ref{D-(F/phi)}), (\ref{D-(F/psi)}) are well-defined and
real-valued if and only if $F_-'$ exists and is real-valued,
in which case,
\begin{align}
- D_{\psi / \varphi}^- (F/\varphi) (x) & = - \frac{\varphi (x)
F_-' (x) - \varphi' (x) F(x)} {\varphi (x) \psi' (x) - \varphi' (x)
\psi (x)} \stackrel{(\ref{C})}{=} - \frac{\varphi (x) F_-' (x) -
\varphi' (x) F(x)} {Cp' (x)} , \nonumber \\
D_{\varphi / \psi}^- (F/\psi) (x) & = \frac{\psi (x) F_-' (x)
- \psi' (x) F(x)} {\varphi' (x) \psi (x) - \varphi (x) \psi' (x)}
\stackrel{(\ref{C})}{=} - \frac{\psi (x) F_-' (x) - \psi' (x) F(x)}
{Cp' (x)} . \nonumber
\end{align}
The function $- D_{\psi / \varphi}^- (F/\varphi)$ is increasing
if and only if its first distributional derivative is a positive
measure, namely, if and only if the second distributional
derivative of $F$ is a measure and
\be
\frac{\varphi (x)}{Cp' (x)} F'' (dx) - \frac{\varphi'' (x)}{Cp' (x)}
F(x) \, dx - \bigl[ \varphi (x) F_-' (x) - \varphi' (x) F(x) \bigr]
\frac{p''(x)}{C \bigl[ p'(x) \bigr]^2} \, dx \leq 0 .
\ee
In view of the definition (\ref{scale}) of the scale function
$p$ and the fact that $p$ and $C$ are both strictly positive,
we can see that this is true if and only if
\be
\varphi (x) \half \sigma^2 (x) F'' (dx) + \varphi (x) b(x) F_-' (x)
\, dx - \left[ \half \sigma^2 (x) \varphi'' (x) + b(x) \varphi' (x)
\right] F(x) \, dx \leq 0 ,
\ee
which is true if and only if $-\Lop F \geq 0$, thanks to the fact
that $\varphi > 0$ satisfies the ODE (\ref{ODEhom}).
Similarly, we can see that the function $D_{\varphi / \psi}^-
(F/\psi)$ is increasing if and only if $-\Lop F \geq 0$.
\mbox{}\hfill$\Box$

\section{The solution of the optimal stopping problem}
\label{sec:solution}

Before addressing the main results of the section, we prove
that the value function $v$ is excessive.

\begin{lem} \label{lem:exc}
Consider the optimal stopping problem formulated in
Section~\ref{sec:pr-form} and suppose that its value function
is real-valued.
The value function $v$ is $r(\cdot)$-excessive, i.e.,
\ben
\e_x \left[ e^{-\Lambda _\tau} v(X_\tau) {\bf 1} _{\{ \tau < \infty
\}} \right] \leq v(x) , \label{exc}
\een
for all initial conditions $x \in \CI$ and every stopping strategy
$(\str_x, \tau) \in \T_x$.
Also,
\ben
v(x) = \sup _{(\str_x, \tau) \in \T_x} \e_x \bigl[ e^{-\Lambda
_{\tau \wedge T_\alpha \wedge T_\beta}} \overline{f}
(X_{\tau \wedge T_\alpha \wedge T_\beta}) {\bf 1}
_{\{ \tau < \infty \}} \bigr] \quad \text{for all } x \in \CI ,
\label{v-over}
\een
where $\overline{f}$ is given by (\ref{f-over}).
\end{lem}
{\bf Proof.}
To prove the $r(\cdot)$-excessivity of $v$, we first show that
$v$ is continuous in $\intI$ and lower semicontinuous in $\CI$.
To this end, we consider any points $x, y \in \intI$.
Given the stopping strategy $(\str_x, T_y) \in \T_x$ and any
stopping strategy $(\str_y, \tau) \in \T_y$, we denote by
$(\hat{\str}_x, \hat{\tau})$ a stopping strategy that is as in
Corollary~\ref{cor:paste}, so that
\begin{align}
v(x) \geq J(\hat{\str}_x, \hat{\tau}) & = \e_x \left[
e^{-\Lambda _{T_\alpha \wedge T_\beta}} f(X_{T_\alpha \wedge
T_\beta}) {\bf 1} _{\{ T_\alpha \wedge T_\beta < T_y \}} \right]
+ \e_x \left[ e^{-\Lambda _{T_y}} {\bf 1} _{\{ T_y < T_\alpha
\wedge T_\beta \}} \right] J(\str_y, \tau) \nonumber \\
& \geq \e_x \left[ e^{-\Lambda _{T_y}} {\bf 1} _{\{ T_y < T_\alpha
\wedge T_\beta \}} \right] J(\str_y, \tau) . \nonumber
\end{align}
Since $(\str_y, \tau)$ is arbitrary, we can use the dominated
convergence theorem to see that this inequality implies that
\be
v(x) \geq \lim _{y \rightarrow x} \e_x \left[ e^{-\Lambda _{T_y}}
{\bf 1} _{\{ T_y < T_\alpha \wedge T_\beta \}}\right] \limsup
_{y \rightarrow x} v(y) = \limsup _{y \rightarrow x} v(y) ,
\ee
which proves that $v$ is upper semicontinuous in $\intI$.

Repeating the same arguments with the roles of the
points $x, y \in \intI$ reversed, we can see that
\be
\liminf _{y \rightarrow x} v(y) \geq \lim _{y \rightarrow x}
\e_y \left[ e^{-\Lambda _{T_x}} {\bf 1} _{\{ T_x < T_\alpha
\wedge T_\beta \}}\right] v(x) = v(x) .
\ee
If both $\alpha$ and $\beta$ are absorbing, then we can
use (\ref{phi})--(\ref{phi-psi-babs}) to calculate
\begin{align}
\liminf _{x \in \intI , \, x \downarrow \alpha} v(x) & \geq \liminf
_{x \in \intI , \, x \downarrow \alpha} \left( f(\alpha) \e_x \left[
e^{-\Lambda _{T_\alpha}} {\bf 1} _{\{ T_\alpha < \T_\beta \}}
\right] + f(\beta) \e_x \left[ e^{-\Lambda _{T_\beta}} {\bf 1}
_{\{ T_\beta < T_\alpha \}} \right] \right) \nonumber \\
& = \liminf _{x \downarrow \alpha} \left( \frac{f(\alpha) \varphi
(x)}{\varphi (\alpha)} + \frac{f(\beta) \psi (x)}{\psi (\beta)} \right)
= f(\alpha) = v(\alpha) , \nonumber
\end{align}
while, if $\alpha$ is absorbing and $\beta$ is inaccessible,
then
\be
\liminf _{x \in \intI , \, x \downarrow \alpha} v(x) \geq \liminf
_{x \in \intI , \, x \downarrow \alpha} f(\alpha) \e_x \left[
e^{-\Lambda _{T_\alpha}} \right] = \liminf _{x \downarrow \alpha}
\frac{f(\alpha) \varphi (x)}{\varphi (\alpha)}  = f(\alpha) = v(\alpha).
\ee
If $\beta$ is absorbing, then we can see that $\liminf _{x \in
\intI , \, x \uparrow \beta} v(x) \geq v(\beta)$ similarly.
It follows that $v$ is lower semicontinuous in $\CI$.

To show that $v$ satisfies (\ref{exc}), we consider any stopping
strategy $(\str_x, \tau) \in \T_x$.
We assume that $X_\tau {\bf 1} _{\{ \tau < T_\alpha \wedge T_\beta \}}$
takes values in a finite set $\{ a_1 , \ldots , a_n \} \subset \intI$.
For each $i = 1 , \ldots , n$, we consider an $\varepsilon$-optimal
strategy $(\str_{a_i}^\varepsilon, \tau_i^\varepsilon) \in
\T_{a_i}$.
If we denote by $(\str_x^\varepsilon, \tau^\varepsilon) \in \T_x$
a stopping strategy that  is as in Corollary~\ref{cor:paste}, then
\begin{align}
v(x) \geq J(\str_x^\varepsilon, \tau^\varepsilon) & = \e_x \left[
e^{-\Lambda _{T_\alpha \wedge T_\beta}} f(X_{T_\alpha \wedge
T_\beta}) {\bf 1} _{\{ T_\alpha \wedge T_\beta < \tau \}} \right]
+ \sum _{i=1}^n \e_x \left[ e^{-\Lambda _\tau} {\bf 1} _{\{ X_\tau = a_i \}}
\right] J(\str_{a_i}^\varepsilon, \tau_i^\varepsilon) \nonumber \\
& \geq \e_x \left[ e^{-\Lambda _{T_\alpha \wedge T_\beta}}
v(X_{T_\alpha \wedge T_\beta}) {\bf 1} _{\{ T_\alpha \wedge
T_\beta < \tau \}} \right] + \sum _{i=1}^n \e_x \left[ e^{-\Lambda
_\tau} {\bf 1} _{\{ X_\tau = a_i \}} \right] \left[ v(a_i) - \varepsilon
\right] , \nonumber
\end{align}
where the last inequality follows from the fact that $f(X_{T_\alpha
\wedge T_\beta}) = v(X_{T_\alpha \wedge T_\beta})$ and the
$\varepsilon$-optimality of the strategies $(\str_{a_i}^\varepsilon,
\tau_i^\varepsilon)$.
Since $\varepsilon > 0$ is arbitrary, it follows that
\be
v(x) \geq \e_x \left[ e^{-\Lambda _{\tau \wedge T_\alpha \wedge
T_\beta}} v(X_{\tau \wedge T_\alpha \wedge T_\beta}) {\bf 1}
_{\{ \tau < \infty \}} \right] ,
\ee
and (\ref{exc}) follows in this case.

Now, we consider any stopping strategy $(\str_x, \tau) \in \T_x$,
and we define
\be
\tau _n = \inf \left\{ t \geq \tau \mid \ X_t \in \{ a_1 ,
\ldots , a_n \} \right\} ,
\ee
where $(a_n)$ is any sequence that is dense in $\intI$.
Such a sequence of stopping times is such that
\be
\tau_n {\bf 1} _{\{ T_\alpha \wedge T_\beta \leq \tau \}} = \infty
{\bf 1} _{\{ T_\alpha \wedge T_\beta \leq \tau \}} \text{ for all }
n \geq 1 \quad \text{and} \quad
\lim _{n \rightarrow \infty} \tau_n {\bf 1} _{\{ \tau < T_\alpha
\wedge T_\beta \}} = \tau {\bf 1} _{\{ \tau < T_\alpha \wedge
T_\beta \}} .
\ee
Therefore, $\lim _{n \rightarrow \infty} \tau_n \wedge T_\alpha
\wedge T_\beta= \tau \wedge T_\alpha \wedge T_\beta$.
Our analysis above has established that (\ref{exc}) holds true
for each of the stopping strategies $(\str_x, \tau_n) \in \T_x$.
Combining this observation with Fatou's lemma and the fact
that $v$ is lower semicontinuous, we can see that 
\be
v(x) \geq \liminf _{n \rightarrow \infty} \e_x \left[ e^{-\Lambda
_{\tau_n \wedge T_\alpha \wedge T_\beta}} v(X_{\tau_n \wedge
T_\alpha \wedge T_\beta}) {\bf 1} _{\{ \tau _n < \infty \}} \right]
\geq \e_x \left[ e^{-\Lambda _{\tau \wedge T_\alpha \wedge
T_\beta}} v(X_{\tau \wedge T_\alpha \wedge T_\beta}) {\bf 1}
_{\{ \tau < \infty \}} \right] ,
\ee
which establishes (\ref{exc}).

Finally, we note that the continuity properties of $v$ and the
inequality $v \geq f$ imply that $v \geq \overline{f}$.
This observation and the $r(\cdot)$-excessivity of $v$ imply
that
\begin{align}
v(x) & =  \sup _{(\str_x, \tau) \in \T_x} \e_x
\left[ e^{-\Lambda _\tau} v(X_\tau) {\bf 1} _{\{ \tau
< \infty \}} \right] \nonumber \\
& \geq \sup _{(\str_x, \tau) \in \T_x} \e_x \left[
e^{-\Lambda _\tau} \overline{f} (X_\tau) {\bf 1} _{\{ \tau <
\infty \}} \right] \geq  \sup _{(\str_x, \tau) \in \T_x} \e_x
\left[ e^{-\Lambda _\tau} f(X_\tau) {\bf 1} _{\{ \tau
< \infty \}} \right] = v(x) , \nonumber
\end{align}
and (\ref{v-over}) follows.
\mbox{}\hfill$\Box$
\vspace{5mm}

Our main results in this section involve solutions to the variational
inequality
\ben
\max \left\{ \Lop v , \ \overline{f} - v \right\} = 0 \label{QVI}
\een
in the following sense.

\begin{defn} \label{sense} {\rm
A function $v: \CI \rightarrow \R_+$ satisfies the variational
inequality (\ref{QVI}) if its restriction in $\intI$ is the difference
of two convex functions,
\begin{gather}
- \Lop v \text{ is a positive measure on } \bigl( \intI,
{\mathcal B} (\intI) \bigr) , \label{QVI1} \\
f(x) \leq v(x) \text{ for all } x \in \intI , \label{QVI2} \\
\text{and the measure } \Lop v \text{ does not charge the
open set } \{ x \in \intI \mid \ v(x) > \overline{f}(x) \} ,
\label{QVI3}
\end{gather}
where $\Lop$ is defined by (\ref{lop}) and $\overline{f}$ is
defined by (\ref{f-over}).
} \mbox{}\hfill$\Box$ \end{defn}

We now prove that the value function $v$ satisfies the
variational inequality (\ref{QVI}) in the sense of this definition.
Also, we establish sufficient conditions for the existence of
$\varepsilon$-optimal as well as optimal stopping strategies.
It is worth noting that the requirements
(\ref{main-abseq3})--(\ref{main-abseq4}) are not really
needed: the only reason we have adopted them is to simplify
the exposition of the proof.

\begin{thm} \label{thm:main1}
Consider the optimal stopping problem formulated in
Section~\ref{sec:pr-form}.
The following statements are true.

\noindent
{\rm (I)}
If the problem data is such that
\be
\overline{f} (y) = \infty, \text{ for some } y \in \CI , \quad
\text{or} \quad \limsup _{y \downarrow \alpha} \frac{f(y)}
{\varphi (y)} = \infty \quad \text{or} \quad
\limsup _{y \uparrow \beta} \frac{f(y)}{\psi (y)} = \infty ,
\ee
then $v(x) = \infty$ for all $x \in \CI$, otherwise, $v(x) <
\infty$ for all $x \in \CI$.

\noindent
{\rm (II)}
If the problem data is such that
\ben
\overline{f} (y) < \infty \text{ for all } y \in \CI , \quad
\limsup _{y \downarrow \alpha} \frac{f(y)}{\varphi (y)} < \infty
\quad \text{and} \quad
\limsup _{y \uparrow \beta} \frac{f(y)}{\psi (y)} < \infty ,
\label{gr-cond}
\een
then the value function $v$ satisfies the variational inequality
(\ref{QVI}) in the sense of Definition~\ref{sense},
\begin{gather}
\lim _{y \in \intI , \, y \downarrow \alpha} \frac{v(y)}{\varphi (y)}
= \limsup _{y \downarrow \alpha} \frac{f(y)}{\varphi (y)} ,
\qquad \lim _{y \in \intI , \, y \uparrow \beta} \frac{v(y)}{\psi (y)}
= \limsup _{y \uparrow \beta} \frac{f(y)}{\psi (y)} \label{v-lims} \\
\intertext{and}
v(\alpha) = f (\alpha) \ \bigl( \text{resp., } v(\beta) = f (\beta)
\bigr) \text{ if } \alpha \ \bigl( \text{resp., } \beta \bigr)
\text{ is absorbing} . \label{v-absBC}
\end{gather}

\noindent
{\rm (III)}
Suppose that (\ref{gr-cond}) is true and that
$f = \overline{f}$.
Given an initial condition $x \in \intI$ consider any monotone
sequences $(\alpha_n)$, $(\beta_n)$ in $\CI$ such that
\begin{gather}
\alpha_1 < x < \beta_1 , \quad \lim _{n \rightarrow \infty} \alpha_n
= \alpha , \quad \lim _{n \rightarrow \infty} \beta_n = \beta
, \label{main-abseq1} \\
\lim _{n \rightarrow \infty} \frac{f(\alpha_n)} {\varphi (\alpha_n)}
= \limsup _{y \downarrow \alpha} \frac{f(y)}{\varphi (y)} , \quad
\lim _{n \rightarrow \infty} \frac{f(\beta_n)} {\psi (\beta_n)} =
\limsup _{y \uparrow \beta} \frac{f(y)}{\psi (y)}
, \label{main-abseq2} \\
\text{if } \alpha \text{ is absorbing and } f(\alpha) = \limsup
_{y \downarrow \alpha} f(y) , \text{ then } \alpha_n = \alpha
\text{ for all } n \geq 1 , \label{main-abseq3} \\
\text{and if } \beta \text{ is absorbing and } f(\beta) = \limsup
_{y \uparrow \beta} f(y) , \text{ then } \beta_n = \beta
\text{ for all } n \geq 1 . \label{main-abseq4}
\end{gather}
Also, let $\str_x$ be any weak solution to (\ref{SDE}), and define
the associated stopping times
\ben
\tau ^\star = \inf \left\{ t \geq 0 \mid \ v(X_t) = f(X_t)
\right\} \quad \text{and} \quad \tau _n^\star = \tau ^\star \wedge
T_{\alpha_n} \wedge T_{\beta_n} . \label{t*}
\een
Then
\ben
v(x) = \lim _{n \rightarrow \infty} \e_x \left[ e^{-\Lambda _{\tau
_n^\star}} f(X_{\tau _n^\star}) \right] . \label{eps-opt}
\een
Furthermore, the stopping strategy $(\str_x, \tau^\star) \in
\T_x$ is optimal if
\begin{gather}
\limsup _{y \downarrow \alpha} \frac{f(y)}{\varphi (y)} = 0
\text{ if } \alpha \text{ is inaccessible} , \quad
\limsup _{y \uparrow \beta} \frac{f(y)}{\psi (y)} = 0 \text{ if }
\beta \text{ is inaccessible} , \label{t*-cond121} \\
f(\alpha) = \limsup _{y \downarrow \alpha} f(y) \text{ if }
\alpha \text{ is absorbing} \quad \text{and} \quad
f(\beta) = \limsup _{y \uparrow \beta} f(y) \text{ if }
\beta \text{ is absorbing} . \label{t*-cond122}
\end{gather}
\end{thm}
{\bf Proof.}
We have established part~(I) of the theorem in Lemma~\ref{lem:v<oo},
so we assume that (\ref{gr-cond}) holds in what follows.
In view of (\ref{v-phi-psi-lims}) and the fact that $-\Lop v$ is a
positive measure on $\bigl( \intI , {\mathcal B} (\intI) \bigr)$
(see Theorem~\ref{prop:excessive}.(I) and
Lemma~\ref{lem:exc}), we can see that the restriction of $v$
in $\intI$ satisfies all of the assumptions of
Theorem~\ref{cor:F-repres}.
Therefore, the limits of $v/\varphi$ and $v/\psi$ in (\ref{v-lims})
exist,
\ben
v(x) = \lim _{y \in \intI , \, y \downarrow \alpha} \frac{v(y)}
{\varphi (y)} \varphi (x) + R_{-\Lop v} (x) + \lim _{y \in \intI ,
\, y \uparrow \beta} \frac{v(y)}{\psi (y)} \psi (x) \quad
\text{for all } x \in \intI , \label{v-repr}
\een
and, given any stopping strategy $(\str_x, \tau) \in \T_x$,
\ben
\e_x \left[ e^{-\Lambda _{\tau \wedge T_{\bar{\alpha}} \wedge
T_{\bar{\beta}}}} \tilde{v} (X_{\tau \wedge T_{\bar{\alpha}} \wedge
T_{\bar{\beta}}}) \right] = v(x) + \e_x \left[ \int _0^{\tau \wedge
T_{\bar{\alpha}} \wedge T_{\bar{\beta}}} e^{-\Lambda_u} \,
dA_u^{\Lop v} \right] \label{Dynkin-prop-tilde}
\een
for all $\bar{\alpha} < x < \bar{\beta}$ in $\CI$, where
\ben
\tilde{v} (x) = \begin{cases} v(x) , & \text{if } x \in \intI , \\
\lim _{y \in \intI , \, y \downarrow \alpha} v(y) , & \text{if }
\alpha \text{ is absorbing and } x=\alpha , \\ \lim _{y \in \intI
, \, y \downarrow \alpha} v(y) , & \text{if } \beta \text{ is
absorbing and } x=\beta . \end{cases} \label{v-tilde}
\een
If $\alpha$ (resp., $\beta$) is absorbing, then (\ref{v-absBC})
plainly holds true and
\be
f(\alpha) = v(\alpha) \leq \liminf _{y \in \intI , \, y \downarrow \alpha}
v(y) \quad \left( \text{resp., } f(\beta) = v(\beta) \leq \liminf _{y \in
\intI , \, y \uparrow \beta} v(y) \right) ,
\ee
thanks to the $r(\cdot)$-excessivity of $v$ (see
Theorem~\ref{prop:excessive}.(II) and Lemma~\ref{lem:exc}).
Combining this observation with the fact that the limit of
$v/\varphi$ in (\ref{v-lims}) exists, we can see that
\be
\lim _{y \in \intI , \, y \downarrow \alpha} \frac{v(y)}{\varphi (y)}
= \limsup _{y \downarrow \alpha} \frac{v(y)}{\varphi (y)}
\stackrel{(\ref{v-phi-psi-lims})}{=} \limsup _{y \downarrow
\alpha} \frac{f(y)}{\varphi (y)} .
\ee
We can establish the second identity in (\ref{v-lims}) similarly.

With each initial condition $x \in \intI$, we associate any monotone
sequences $(\alpha_n)$, $(\beta_n)$ in $\CI$ such that
(\ref{main-abseq1})--(\ref{main-abseq4}) hold true.
If $\alpha$ (resp., $\beta$) is absorbing and $\alpha_n = \alpha$
(resp., $\beta_n = \beta$), then (\ref{v-lims})--(\ref{v-absBC})
and (\ref{main-abseq3})--(\ref{main-abseq4}) imply that
\be
v(\alpha) = \lim _{y \in \intI , \, y \downarrow \alpha} v(y)
\quad \left( \text{resp., } v(\beta) = \lim _{y \in \intI , \, y
\uparrow \beta} v(y) \right) .
\ee
This observation, the definition of $\tilde{v}$ in (\ref{v-tilde})
and (\ref{Dynkin-prop-tilde}) imply that
\ben
\e_x \left[ e^{-\Lambda _{\tau \wedge T_{\alpha_n} \wedge
T_{\beta_n}}} v (X_{\tau \wedge T_{\alpha_n} \wedge
T_{\beta_n}}) \right] = v(x) + \e_x \left[ \int _0^{\tau \wedge
T_{\alpha_n} \wedge T_{\beta_n}} e^{-\Lambda_u} \, dA_u^{\Lop
v} \right] \label{Dynkin-prop}
\een
for every stopping strategy $(\str_x, \tau) \in \T_x$.
Furthermore, (\ref{v-lims}) and (\ref{main-abseq2}) imply that
\begin{gather}
\lim _{n \rightarrow \infty} \frac{v(\alpha_n)} {\varphi (\alpha_n)}
= \lim _{n \rightarrow \infty} \frac{\overline{f}(\alpha_n)}
{\varphi (\alpha_n)} = \lim _{n \rightarrow \infty}
\frac{f(\alpha_n)} {\varphi (\alpha_n)} \label{vf-lims1} \\
\intertext{and}
\lim _{n \rightarrow \infty} \frac{v(\beta_n)} {\psi (\beta_n)}
= \lim _{n \rightarrow \infty} \frac{f(\beta_n)} {\psi (\beta_n)}
= \lim _{n \rightarrow \infty} \frac{\overline{f}(\beta_n)}
{\psi (\beta_n)} . \label{vf-lims2}
\end{gather}

Given a stopping strategy $(\str_x, \tau) \in \T_x$ such that
$\tau = \tau \wedge T_\alpha \wedge T_\beta$, we can
use (\ref{phi})--(\ref{psi}) to calculate
\begin{align}
\e_x \biggl[ & e^{-\Lambda _{\tau \wedge T_{\alpha_n} \wedge
T_{\beta_n}}} \left[ v (X_{\tau \wedge T_{\alpha_n} \wedge
T_{\beta_n}}) - \overline{f}(X_{\tau \wedge T_{\alpha_n} \wedge
T_{\beta_n}}) \right] \biggr] \nonumber \\
= \mbox{} & \e_x \left[ e^{-\Lambda _\tau} \left[ v (X_\tau) -
\overline{f}(X_\tau) \right] {\bf 1} _{\{ \tau \leq T_{\alpha_n}
\wedge T_{\beta_n} \}} \right] + \left[ v(\alpha_n) - \overline{f}
(\alpha_n) \right] \e_x \left[ e^{-\Lambda _{T_{\alpha_n}}} {\bf 1}
_{\{ T_{\alpha_n} < \tau \wedge T_{\beta_n} \}} \right] \nonumber \\
& + \left[ v(\beta_n) - \overline{f}(\beta_n) \right] \e_x \left[
e^{-\Lambda _{T_{\beta_n}}} {\bf 1} _{\{ T_{\beta_n} < \tau
\wedge T_{\alpha_n} \}} \right] \nonumber \\
= \mbox{} & \e_x \left[ e^{-\Lambda _\tau} \left[ v(X_\tau) -
\overline{f}(X_\tau) \right] {\bf 1} _{\{ \tau \leq T_{\alpha_n}
\wedge T_{\beta_n} \}} \right] + \varphi (x) \frac{v(\alpha_n) -
\overline{f} (\alpha_n)} {\varphi (\alpha_n)} \frac{\e_x \left[
e^{-\Lambda _{T_{\alpha_n}}} {\bf 1} _{\{ T_{\alpha_n} < \tau
\wedge T_{\beta_n} \}} \right]} {\e_x \left[ e^{-\Lambda
_{T_{\alpha_n}}} \right]} \nonumber \\
& + \psi (x) \frac{v(\beta_n) - \overline{f} (\beta_n)} {\psi
(\beta_n)} \frac{\e_x \left[ e^{-\Lambda _{T_{\beta_n}}} {\bf 1}
_{\{ T_{\beta_n} < \tau \wedge T_{\alpha_n} \}} \right]} {\e_x
\left[ e^{-\Lambda _{T_{\beta_n}}} \right]} . \nonumber
\end{align}
Combining this calculation with (\ref{vf-lims1})--(\ref{vf-lims2})
and the monotone convergence theorem, we can see that
\begin{align}
\lim _{n \rightarrow \infty} & \e_x \left[ e^{-\Lambda
_{\tau \wedge T_{\alpha_n} \wedge T_{\beta_n}}} \bigl[
v( X_{\tau \wedge T_{\alpha_n} \wedge T_{\beta_n}}) -
\overline{f} (X_{\tau \wedge T_{\alpha_n} \wedge T_{\beta_n}})
\bigr] \right] \nonumber \\
& \hspace{57mm} = \e_x \left[ e^{-\Lambda _\tau} \bigl[
v(X_\tau) - \overline{f} (X_\tau) \bigr] {\bf 1} _{\Gamma
(\tau)} \right] , \label{VF-lims}
\end{align}
where
\be
\Gamma (\tau) = \begin{cases}
\{ \tau < T_\alpha \wedge T_\beta \} , & \text{if } \alpha
< \alpha_n < \beta_n < \beta , \\ \{ \tau < T_\beta \} , &
\text{if } \alpha_n = \alpha \text{ and } \beta_n < \beta
, \\ \{ \tau < T_\alpha \} , & \text{if } \alpha < \alpha_n
\text{ and } \beta_n = \beta , \\ \Omega , & \text{if }
\alpha_n = \alpha \text{ and } \beta_n = \beta , \end{cases}
\ee
(see also (\ref{main-abseq3})--(\ref{main-abseq4})).

With each initial condition $x \in \intI$, we associate any
sequence of stopping strategies $(\str_x^\ell, \tau _\ell)
\in \T_x$ such that $\tau_\ell = \tau_\ell \wedge T_\alpha
\wedge T_\beta$ and
\be
v(x) -\frac{1}{2\ell} \leq \e_x^\ell \left[ e^{-\Lambda _{\tau_\ell}}
\overline{f} (X_{\tau_\ell}) {\bf 1} _{\{ \tau _\ell < \infty \}}
\right] \quad \text{for all } \ell \geq 1
\ee
(see (\ref{v-over}) in Lemma~\ref{lem:exc}).
If $\alpha$ is absorbing and $\alpha < \alpha_n$ (see
(\ref{main-abseq3})), then we may assume without
loss of generality that $\tau_\ell < T_\alpha$, $\p_x^\ell$-a.s..
To see this claim, suppose that $\alpha$ is absorbing and
$\alpha < \alpha_n$, which is the case when
$f(\alpha) < \limsup _{y \downarrow \alpha} f(y)$.
Since $\tau_\ell = \tau_\ell \wedge T_\alpha \wedge T_\beta$,
\begin{align}
\bigcap _{n=1}^\infty \{ T_{\alpha_n} < \tau_\ell \} & = \bigcap
_{n=1}^\infty \{ T_{\alpha_n} < \tau_\ell \wedge T_\alpha \}
\cap \{ T_{\alpha_n} < T_\beta \} \nonumber \\
& = \{ T_\alpha \leq \tau_\ell \wedge T_\alpha \} \cap \{
T_\alpha < T_\beta \} \nonumber \\
& = \{ \tau_\ell = T_\alpha \} \cap \{ T_\alpha < T_\beta \} .
\nonumber
\end{align}
In view of  this observation and the dominated convergence
theorem, we can see that
\begin{align}
\lim _{n \rightarrow \infty} \e \bigl[ e^{- \Lambda_{\tau_\ell}}
\bigl[ f(\alpha_n) - f(X_{\tau_\ell}) \bigr] & {\bf 1} _{\{
T_{\alpha_n} < \tau_\ell \}} {\bf 1} _{\{ \tau _\ell < \infty \}} \bigr]
\nonumber \\
& = \left[ \limsup _{y \downarrow \alpha} f(y) - f(\alpha) \right]
\e \left[ e^{- \Lambda_{\tau_\ell}} {\bf 1} _{\{ \tau_\ell = T_\alpha
\} \cap \{ T_\alpha < T_\beta \}} \right] . \nonumber
\end{align}
If $\p_x^\ell ( \tau_\ell = T_\alpha ) > 0$, then the right-hand
side of this identity is strictly positive, and there exists $k \geq 1$
such that
\be
f(\alpha_k) \e \bigl[ e^{- \Lambda_{\tau_\ell}} {\bf 1}
_{\{ T_{\alpha_k} < \tau_\ell \}} \bigr] \geq \e \bigl[
e^{- \Lambda _{\tau_\ell}} f(X_{\tau_\ell}) {\bf 1}
_{\{ T_{\alpha_k} < \tau_\ell \}} {\bf 1} _{\{ \tau _\ell < \infty \}}
\bigr] .
\ee
Given such a $k$, we can see that
\begin{align}
\e \bigl[ e^{- \Lambda _{\tau_\ell \wedge T_{\alpha_k}}}
& f(X_{\tau_\ell \wedge T_{\alpha_k}}) {\bf 1} _{\{ \tau _\ell
\wedge T_{\alpha_k} < \infty \}} \bigr] \nonumber \\
& = \e \left[ e^{- \Lambda _{\tau_\ell}} f(X_{\tau_\ell}) {\bf 1}
_{\{ \tau_\ell \leq T_{\alpha_k} \} \cap \{ \tau _\ell < \infty \}}
\right] + f(\alpha_k) \e \left[ e^{- \Lambda _{T_{\alpha_k}}}
{\bf 1} _{\{ T_{\alpha_k} < \tau _\ell \}} \right] \nonumber \\
& \geq \e \left[ e^{- \Lambda _{\tau_\ell}} f(X_{\tau_\ell}) {\bf 1}
_{\{ \tau_\ell \leq T_{\alpha_k} \} \cap \{ \tau _\ell < \infty \}}
\right] + \e \bigl[ e^{- \Lambda_{\tau_\ell}} f(X_{\tau_\ell})
{\bf 1} _{\{ T_{\alpha_k} < \tau_\ell \}} {\bf 1} _{\{ \tau _\ell
< \infty \}} \bigr] \nonumber \\
& = \e \left[ e^{- \Lambda _{\tau_\ell}} f(X_{\tau_\ell}) {\bf 1}
_{\{ \tau _\ell < \infty \}} \right] , \nonumber
\end{align}
and the claim follows.
Similarly, we may assume that $\tau_\ell < T_\beta$, $\p_x^\ell$-a.s.,
if $\beta$ is absorbing and $\beta_n < \beta$.

In light of the above observations and
(\ref{main-abseq3})--(\ref{main-abseq4}), we can use the
monotone convergence theorem to calculate
\begin{align}
\liminf _{n \rightarrow \infty} \e_x^\ell \left[ e^{-\Lambda
_{\tau_\ell \wedge T_{\alpha_n} \wedge T_{\beta_n}}} \overline{f}
(X_{\tau_\ell \wedge T_{\alpha_n} \wedge T_{\beta_n}}) \right]
& \geq \lim _{n \rightarrow \infty} \e_x^\ell \left[ e^{-\Lambda
_{\tau_\ell}} \overline{f} (X_{\tau_\ell}) {\bf 1} _{\{ \tau _\ell
\leq T_{\alpha_n} \wedge T_{\beta_n} \}} \right] \nonumber \\
& = \e_x^\ell \left[ e^{-\Lambda _{\tau_\ell}} \overline{f}
(X_{\tau_\ell}) {\bf 1} _{\{ \tau _\ell < \infty \}} \right] ,
\nonumber
\end{align}
which implies that, for all $\ell \geq 1$, there exists $n_\ell$ such
that
\be
\e_x^\ell \left[ e^{-\Lambda _{\tau_\ell}} \overline{f} (X_{\tau_\ell})
{\bf 1} _{\{ \tau _\ell < \infty \}} \right] \leq \e_x^\ell \left[
e^{-\Lambda _{\tau_\ell \wedge T_{\alpha_{n_\ell}} \wedge
T_{\beta_{n_\ell}}}} \overline{f} (X_{\tau_\ell \wedge
T_{\alpha_{n_\ell}} \wedge T_{\beta_{n_\ell}}})\right] +
\frac{1}{2\ell} .
\ee
It follows that, if we define
\ben
\tau_\ell^\circ = \tau_\ell \wedge T_{\alpha_{n_\ell}} \wedge
T_{\beta_{n_\ell}} , \label{toll}
\een
then the stopping strategy $(\str_x^\ell, \tau_\ell^\circ)
\in \T_x$ satisfies
\begin{align}
v(x) - \e_x^\ell \left[ e^{-\Lambda _{\tau_\ell^\circ}} \overline{f}
(X_{\tau_\ell^\circ}) \right] \leq \frac{1}{\ell} . \label{t*1}
\end{align}
In view of (\ref{Dynkin-prop}) and (\ref{toll}), we can see that
\ben
v(x) - \e_x^\ell \left[ e^{-\Lambda _{\tau_\ell^\circ}} \overline{f}
(X_{\tau_\ell^\circ}) \right] = \e_x^\ell \left[ e^{-\Lambda
_{\tau_\ell^\circ}} \bigl[ v (X_{\tau_\ell^\circ}) - \overline{f}
(X_{\tau_\ell^\circ}) \bigr] \right] + \e_x^\ell \left[ - \int
_0^{\tau_\ell^\circ} e^{-\Lambda _u} \, dA_u^{\Lop v} \right] .
\label{t*2}
\een
The first term on the right-hand side of this identity is
clearly positive, while, the second one is positive because
$-\Lop v$ is a positive measure and $-A^{\Lop v}$ is an
increasing process (see also (\ref{A-lin}) in
Lemma~\ref{lem:Amu}).
This observation and (\ref{t*1})--(\ref{t*2})  imply that
\ben
\lim _{\ell \rightarrow \infty} \e_x^\ell \left[ e^{-\Lambda
_{\tau_\ell^\circ}} \bigl[ v (X_{\tau_\ell^\circ}) - \overline{f}
(X_{\tau_\ell^\circ}) \bigr] \right] = \lim _{\ell \rightarrow \infty}
\e_x^\ell \left[ - \int _0^{\tau_\ell^\circ} e^{-\Lambda _u} \,
dA_u^{\Lop v} \right] = 0 . \label{t*3}
\een

{\em Proof of} (II).
To prove that $v$ satisfies the variational inequality
(\ref{QVI}) in the sense of Definition~\ref{sense}, and thus
complete the proof of part~(II) of the theorem, we
have to show that (\ref{QVI3}) holds true because
$v \geq \overline{f}$ and $-\Lop v$ is a positive measure. 
To this end, we consider any interval $[\tilde{\alpha},
\tilde{\beta}] \subseteq \{ x \in \intI \mid \ v(x) > \overline{f}
(x) \}$ and we note that there exists $\xi>0$ such that
\be
\xi \leq \min _{x \in [\tilde{\alpha}, \tilde{\beta}]} \bigl[ v(x)
- \overline{f} (x) \bigr] \leq \max _{x \in [\tilde{\alpha},
\tilde{\beta}]} v(x) \leq \xi^{-1}
\ee
because the restrictions of $v-\overline{f}$ and $v$ in
$\intI$ are lower semicontinuous and continuous, respectively.
In view of this observation, we can see that
\begin{gather}
e^{-\Lambda _{\tau_\ell^\circ}} \bigl[ v (X_{\tau_\ell^\circ})
- \overline{f} (X_{\tau_\ell^\circ}) \bigr] \geq \xi
e^{-\Lambda _{\tau_\ell^\circ}} {\bf 1} _{\{ \tau_\ell^\circ
< T_{\tilde{\alpha}} \wedge T_{\tilde{\beta}} \}}
\geq \xi e^{-\Lambda _{T_{\tilde{\alpha}}}} {\bf 1} _{\{
\tau_\ell^\circ < T_{\tilde{\alpha}} < T_{\tilde{\beta}} \}}
+ \xi e^{-\Lambda _{T_{\tilde{\beta}}}} {\bf 1} _{\{
\tau_\ell^\circ < T_{\tilde{\beta}} < T_{\tilde{\alpha}} \}}
\nonumber \\
\intertext{and}
e^{-\Lambda _{\tau_\ell^\circ}} \bigl[ v (X_{\tau_\ell^\circ})
- \overline{f} (X_{\tau_\ell^\circ}) \bigr] \geq \xi
e^{-\Lambda _{\tau_\ell^\circ}} {\bf 1} _{\{ \tau_\ell^\circ
< T_{\tilde{\alpha}} \wedge T_{\tilde{\beta}} \}} \geq
\xi^2 e^{-\Lambda _{\tau_\ell^\circ}} v(X_{\tau_\ell^\circ})
{\bf 1} _{\{ \tau_\ell^\circ < T_{\tilde{\alpha}} \wedge
T_{\tilde{\beta}} \}} . \nonumber
\end{gather}
These inequalities and (\ref{t*3}) imply that
\begin{gather}
\lim _{\ell \rightarrow \infty} \e_x^\ell \left[ e^{-\Lambda
_{T_{\tilde{\alpha}}}} {\bf 1} _{\{ \tau_\ell^\circ < T_{\tilde{\alpha}}
< T_{\tilde{\beta}} \}} \right] = 0 , \qquad
\lim _{\ell \rightarrow \infty} \e_x^\ell \left[ e^{-\Lambda
_{T_{\tilde{\beta}}}} {\bf 1} _{\{ \tau_\ell^\circ < T_{\tilde{\beta}}
< T_{\tilde{\alpha}} \}} \right] = 0 , \nonumber \\
\lim _{\ell \rightarrow \infty} \e_x^\ell \left[ e^{-\Lambda
_{\tau_\ell^\circ}} v(X_{\tau_\ell^\circ}) {\bf 1} _{\{ \tau_\ell^\circ
< T_{\tilde{\alpha}} \wedge T_{\tilde{\beta}} \}} \right] = 0
\quad \text{and} \quad \lim _{\ell \rightarrow \infty} \e_x^\ell \left[ - \int
_0^{\tau_\ell^\circ \wedge T_{\tilde{\alpha}} \wedge
T_{\tilde{\beta}}} e^{-\Lambda _u} \, dA_u^{\Lop v} \right] = 0 .
\label{t*4}
\end{gather}
The first of these limits implies that
\ben
\lim _{\ell \rightarrow \infty} \e_x^\ell \left[ e^{-\Lambda
_{T_{\tilde{\alpha}}}} {\bf 1} _{\{ T_{\tilde{\alpha}} \leq
\tau_\ell^\circ \wedge T_{\tilde{\beta}} \}} \right] =
\lim _{\ell \rightarrow \infty} \e_x^\ell \left[ e^{-\Lambda
_{T_{\tilde{\alpha}}}} {\bf 1} _{\{ T_{\tilde{\alpha}} <
T_{\tilde{\beta}} \}} \right] \label{t*5}
\een
because $\{ \tau_\ell^\circ < T_{\tilde{\alpha}} <
T_{\tilde{\beta}} \} = \{ T_{\tilde{\alpha}} < T_{\tilde{\beta}}
\} \setminus \{ T_{\tilde{\alpha}} \leq \tau_\ell^\circ \wedge
T_{\tilde{\beta}} \}$.
Similarly, the second limit implies that
\ben
\lim _{\ell \rightarrow \infty} \e_x^\ell \left[ e^{-\Lambda
_{T_{\tilde{\beta}}}} {\bf 1} _{\{ T_{\tilde{\beta}} \leq
\tau_\ell^\circ \wedge T_{\tilde{\alpha}} \}} \right] =
\lim _{\ell \rightarrow \infty} \e_x^\ell \left[ e^{-\Lambda
_{T_{\tilde{\beta}}}} {\bf 1} _{\{ T_{\tilde{\beta}} <
T_{\tilde{\alpha}} \}} \right] . \label{t*6}
\een
Now, (\ref{Dynkin-prop}) and (\ref{toll}) imply that
\begin{align}
v(x) = \mbox{} & \e_x^\ell \left[ e^{-\Lambda _{\tau_\ell^\circ
\wedge T_{\tilde{\alpha}} \wedge T_{\tilde{\beta}}}}
v (X_{\tau_\ell^\circ \wedge T_{\tilde{\alpha}} \wedge
T_{\tilde{\beta}}}) \right] + \e_x^\ell \left[ - \int
_0^{\tau_\ell^\circ \wedge T_{\tilde{\alpha}} \wedge
T_{\tilde{\beta}}} e^{-\Lambda _u} \, dA_u^{\Lop v}
\right] \nonumber \\
= \mbox{} & \e_x^\ell \left[ e^{-\Lambda _{\tau_\ell^\circ}}
v (X_{\tau_\ell^\circ}) {\bf 1} _{\{ \tau_\ell^\circ <
T_{\tilde{\alpha}} \wedge T_{\tilde{\beta}} \}} \right] +
v(\tilde{\alpha}) \e_x^\ell \left[ e^{-\Lambda _{T_{\tilde{\alpha}}}}
{\bf 1} _{\{ T_{\tilde{\alpha}} \leq \tau_\ell^\circ \wedge
T_{\tilde{\beta}} \}} \right] \nonumber \\
& + v(\tilde{\beta}) \e_x^\ell \left[ e^{-\Lambda
_{T_{\tilde{\beta}}}} {\bf 1} _{\{ T_{\tilde{\beta}} \leq
\tau_\ell^\circ \wedge T_{\tilde{\alpha}} \}} \right]
+ \e_x^\ell \left[ - \int _0^{\tau_\ell^\circ \wedge
T_{\tilde{\alpha}} \wedge T_{\tilde{\beta}}} e^{-\Lambda _u}
\, dA_u^{\Lop v} \right] . \nonumber
\end{align}
In view of (\ref{t*4})--(\ref{t*6}), we can pass to the limit
as $\ell \rightarrow \infty$ to obtain
\begin{align}
v(x) & = \lim _{\ell \rightarrow \infty} \left\{ v(\tilde{\alpha})
\e_x^\ell \left[ e^{-\Lambda _{T_{\tilde{\alpha}}}} {\bf 1}
_{\{ T_{\tilde{\alpha}} < T_{\tilde{\beta}} \}} \right] + v(\tilde{\beta})
\e_x^\ell \left[ e^{-\Lambda _{T_{\tilde{\beta}}}} {\bf 1}
_{\{ T_{\tilde{\beta}} < T_{\tilde{\alpha}} \}} \right] \right\}
\nonumber \\
& = v(\tilde{\alpha}) \frac{\varphi (\tilde{\beta}) \psi (x) -
\varphi (x) \psi (\tilde{\beta})} {\varphi (\tilde{\beta}) \psi
(\tilde{\alpha}) - \varphi (\tilde{\alpha}) \psi (\tilde{\beta})}
+ v(\tilde{\beta}) \frac{\varphi (x) \psi (\tilde{\alpha}) - \varphi
(\tilde{\alpha}) \psi (x)} {\varphi (\tilde{\beta}) \psi (\tilde{\alpha})
- \varphi (\tilde{\alpha}) \psi (\tilde{\beta})} , \nonumber
\end{align}
the second identity following from
(\ref{ELam_alpha})--(\ref{ELam_beta}) in
Lemma~\ref{lem:phi-psi-expect}.
Since this identity is true for all $x \in \mbox{}
]\tilde{\alpha}, \tilde{\beta}[$ and $\Lop \varphi = \Lop \psi
= 0$, it follows that the restriction of the measure $\Lop v$
in $x \in \mbox{} ]\tilde{\alpha}, \tilde{\beta}[$ vanishes,
which establishes (\ref{QVI3}).

{\em Proof of} (III).
We now assume that $f = \overline{f}$ and we consider the
stopping times $\tau^\star$ and $\tau_n^\star$ that are defined
by (\ref{t*}) on any given weak solution $\str_x$ to (\ref{SDE}).
In view of (\ref{Dynkin-prop}) and the fact that $v$ satisfies
(\ref{QVI3}), we can see that
\be
v(x) - \e_x \left[ e^{-\Lambda _{\tau_n^\star}} f(X_{\tau_n^\star})
\right] = \e_x \left[ e^{-\Lambda _{\tau_n^\star}} \left[
v(X_{\tau_n^\star}) - f(X_{\tau_n^\star}) \right] \right] .
\ee
Combining this result with the identities
\be
\lim _{n \rightarrow \infty} \e_x \bigl[ e^{-\Lambda
_{\tau_n^\star}} \left[ v(X_{\tau_n^\star}) - f(X_{\tau_n^\star})
\right] \bigr] \stackrel{(\ref{VF-lims})}{=} \e_x \left[ e^{-\Lambda
_{\tau^\star}} \left[ v(X_{\tau^\star}) - f(X_{\tau^\star})
\right] {\bf 1} _{\Gamma (\tau^\star)} \right] = 0 ,
\ee
we obtain (\ref{eps-opt}).

To establish the optimality of $(\str_x, \tau^\star)$ if $f =
\overline{f}$ and (\ref{t*-cond121})--(\ref{t*-cond122}) are
satisfied, we first note that if $\alpha$ is inaccessible, then
\be
0 \leq \lim _{n \rightarrow \infty} f(\alpha_n) \e_x \left[
e^{-\Lambda _{T_{\alpha_n}}} {\bf 1} _{\{ T_{\alpha_n} < \tau^\star
\wedge T_{\beta_n} \}} \right] \leq \lim _{n \rightarrow \infty}
f(\alpha_n) \e_x \left[ e^{-\Lambda _{T_{\alpha_n}}} \right]
\stackrel{(\ref{phi})}{=} \lim _{n \rightarrow \infty}
\frac{f(\alpha_n) \varphi (x)}{\varphi (\alpha_n)} = 0 .
\ee
Similarly, if $\beta$ is inaccessible, then
\be
\lim _{n \rightarrow \infty} f(\beta_n) \e_x \left[
e^{-\Lambda _{T_{\beta_n}}} {\bf 1} _{\{ T_{\beta_n} < \tau^\star
\wedge T_{\alpha_n} \}} \right] = 0 .
\ee
In view of (\ref{main-abseq3})--(\ref{main-abseq4}) and
(\ref{t*-cond122}), we can see that, if $\alpha$ (resp., $\beta$)
is absorbing, then $\alpha_n = \alpha$ (resp., $\beta_n =
\beta$) and
\be
\{ T_{\alpha_n} < \tau^* \wedge T_{\beta_n} \} =
\{ T_\alpha < \tau^* \wedge T_{\beta_n} \} = \emptyset
\quad \bigl( \text{resp., } \{ T_{\beta_n} < \tau^\star \wedge
T_{\alpha_n} \} = \emptyset \bigr) .
\ee
In light of these observations and the monotone convergence theorem,
we can see that
\begin{align}
\lim _{n \rightarrow \infty} \e_x \left[ e^{-\Lambda _{\tau_n^\star}}
f(X_{\tau_n^\star}) \right] = \mbox{} & \lim _{n
\rightarrow \infty} \biggl( \e_x \left[ e^{-\Lambda _{\tau^\star}}
f(X_{\tau^\star}) {\bf 1} _{\{ \tau^\star \leq T_{\alpha_n}
\wedge T_{\beta_n} \}} \right] \nonumber \\
& + f(\alpha_n) \e_x \left[ e^{-\Lambda _{T_{\alpha_n}}} {\bf 1}
_{\{ T_{\alpha_n} < \tau^\star \wedge T_{\beta_n} \}} \right]
+ f(\beta_n) \e_x \left[ e^{-\Lambda _{T_{\beta_n}}} {\bf 1} _{\{
T_{\beta_n} < \tau^\star \wedge T_{\alpha_n} \}} \right] \biggr)
\nonumber \\
= \mbox{} & \e_x \bigl[ e^{-\Lambda _{\tau^\star}} f(X_{\tau^\star})
{\bf 1} _{\{ \tau^\star < \infty \}} \bigr] , \nonumber
\end{align}
and the optimality of $(\str_x, \tau^\star)$ follows thanks to
(\ref{eps-opt}).
\mbox{}\hfill$\Box$
\vspace{5mm}

It is straightforward to see that the variational inequality
(\ref{QVI}) does not have a unique solution.
In the previous result, we proved that the value function $v$
satisfies (\ref{QVI}) as well as the boundary / growth conditions
(\ref{v-lims}).
We now establish a converse result, namely a verification
theorem, which shows that $v$ is the minimal solution to
(\ref{QVI}).

\begin{thm} \label{thm:main2}
Consider the optimal stopping problem formulated in
Section~\ref{sec:pr-form} and suppose that (\ref{gr-cond})
holds true.
The following statements are true.

\noindent
{\rm (I)}
If a function $w: \intI \rightarrow \R_+$ is the difference of two
convex functions such that $-\Lop w$ is a positive measure,
$w(x) \geq f(x)$ for all $x \in \intI$,
\be
\limsup _{y \downarrow \alpha} \frac{w(y)}{\varphi (y)} < \infty
\quad \text{and} \quad \limsup _{y \uparrow \beta} \frac{w(y)}
{\psi (y)} < \infty ,
\ee
then $v(x) \leq w(x)$ for all
$x \in \intI$.

\noindent
{\rm (II)}
If a function $w: \intI \rightarrow \R_+$ is a solution to
the variational inequality (\ref{QVI}) in the sense of
Definition~\ref{sense} that satisfies
\ben
\limsup _{y \in \intI, \, y \downarrow \alpha} \frac{w(y)}{\varphi (y)}
= \limsup _{y \in \CI, \, y \downarrow \alpha} \frac{f(y)}{\varphi (y)}
\quad \text{and} \quad
\limsup _{y \in \intI, \, y \uparrow \beta} \frac{w(y)}{\psi (y)} =
\limsup _{y \in \CI, \, y \uparrow \beta} \frac{f(y)}{\psi (y)} ,
\label{w-lims}
\een
then $v(x) = w(x)$ for all $x \in \intI$.

\noindent
{\rm (III)}
The value function $v$ admits the characterisation
\ben
v(x) = \inf \bigl\{ A \varphi (x) + B \psi (x) \mid \ A, B \geq 0
\text{ and } A \varphi (y) + B \psi (y) \geq f(y) \text{ for all }
y \in \intI \bigr\} \label{vLP}
\een
for all $x \in \intI$.
Furthermore, if $c<d$ are any points in $\CI$ such that
$v(x) > \overline{f} (x)$ for all $x \in \mbox{} ]c,d[$, then
there exist constants $\tilde{A}$, $\tilde{B}$ such that
\ben
v(x) = \tilde{A} \varphi (x) + \tilde{B} \psi (x) \text{ and }
\tilde{A} \varphi (y) + \tilde{B} \psi (y) \geq f(y)
\text{ for all } x \in \mbox{} ]c,d[ \text{ and } y \in \intI .
\label{vLP+}
\een
\end{thm}
{\bf Proof.}
A function $w: \intI \rightarrow \R_+$ that is as in the
statement of part~(I) of the theorem satisfies all of the
requirements of Theorem~\ref{cor:F-repres}.
Therefore, if $\CI$ is not open and we identify $w$ with its
extension on $\CI$ that is given by $w(\alpha) = \lim
_{y \downarrow \alpha} w(y)$ (resp., $w(\beta) = \lim
_{y \uparrow \beta} w(y)$) if $\alpha$ (resp., $\beta$) is
absorbing, then
\ben
\e_x \left[ e^{-\Lambda _{\tau \wedge T_{\alpha_n} \wedge
T_{\beta_n}}} w (X_{\tau \wedge T_{\alpha_n} \wedge
T_{\beta_n}}) \right] = w(x) + \e_x \left[ \int _0^{\tau \wedge
T_{\alpha_n} \wedge T_{\beta_n}} e^{-\Lambda_u} \, dA_u^{\Lop
w} \right] \label{Dynkin-prop-w}
\een
for every stopping strategy $(\str_x, \tau) \in \T_x$, where
$(\alpha_n)$, $(\beta_n)$ are any monotone sequences
in $\CI$ satisfying (\ref{main-abseq1}).
Combining this identity with the fact that $- \Lop w$ is a positive
measure, which implies that $-A^{\Lop w}$ is an increasing process,
we can see that
\ben
\e_x \left[ e^{-\Lambda _{\tau \wedge T_{\alpha_n} \wedge
T_{\beta_n}}} w (X_{\tau \wedge T_{\alpha_n} \wedge
T_{\beta_n}}) \right] \leq w(x).
\een
This inequality and Fatou's lemma imply that
\be
\e_x \left[ e^{-\Lambda _{\tau \wedge T_\alpha \wedge T_\beta}}
w (X_{\tau \wedge T_\alpha \wedge T_\beta}) \right] \leq \liminf
_{n \rightarrow \infty} \e_x \left[ e^{-\Lambda _{\tau \wedge
T_{\alpha_n} \wedge T_{\beta_n}}} w(X_{\tau \wedge T_{\alpha_n}
\wedge T_{\beta_n}}) \right] \leq w(x) ,
\ee
which, combined with the inequality $w \geq f$, proves that
$v(x) \leq w(x)$.

If the function $w$ satisfies (\ref{w-lims}) as well, then we
choose any monotone sequences $(\alpha_n)$, $(\beta_n)$
as in (\ref{main-abseq1})--(\ref{main-abseq4}) and we note
that (\ref{vf-lims1})--(\ref{vf-lims2}) hold true with the extension
of $w$ on $\CI$ considered at the beginning of the proof in
place of $v$.
If we consider the stopping strategies $(\str_x, \tau_n^\star)
\in \T_x$, where
\be
\tau_n^\star = \bigl( \inf \left\{ t \geq 0 \mid \ w(X_t) =
\overline{f} (X_t) \right\} \bigr) \wedge T_{\alpha_n} \wedge
T_{\beta_n} ,
\ee
then we can see that (\ref{VF-lims}) with $w$ in place of $v$
and (\ref{Dynkin-prop-w}) imply that
\begin{align}
\lim _{n \rightarrow \infty} \e_x \left[ e^{-\Lambda _{\tau_n^\star}}
\overline{f} (X_{\tau_n^\star}) \right] & = w(x) + \lim _{n \rightarrow
\infty} \e_x \left[ e^{-\Lambda _{\tau_n^\star}} \bigl[ \overline{f}
(X_{\tau_n^\star}) - w (X_{\tau_n^\star}) \bigr] \right]
\nonumber \\
& = w(x) + \e_x \left[ e^{-\Lambda _{\tau^\star}} \left[ \overline{f}
(X_{\tau^\star}) - w (X_{\tau^\star}) \right] {\bf 1} _{\Gamma
(\tau^\star)} \right] \nonumber \\
& = w(x) . \nonumber
\end{align}
It follows that $v(x) \geq w(x)$ thanks to (\ref{v-over}) in
Lemma~\ref{lem:exc}, which, combined with the inequality
$v(x) \leq w(x)$ that we have established above, implies that
$v(x) = w(x)$.

To show part~(III) of the theorem, we first note that,
given any constants $A, B \in \R$, the function $A\varphi
+ B\psi$ satisfies the variational inequality (\ref{QVI}) if
and only if $A\varphi + B\psi \geq f$.
Combining this observation with part~(I) of the theorem,
we can see that $v(x)$ is less than or equal to the
right-hand side of (\ref{vLP}).
To establish the reverse inequality, we first use
(\ref{phi-psi-expect}) in Lemma~\ref{lem:phi-psi-expect}
and (\ref{Dynkin-prop}) with $\tau \equiv \infty$ to obtain
\begin{align}
&\e_{\bar{x}} \left[ e^{-\Lambda _{T_{\bar{\alpha}} \wedge
T_{\bar{\beta}}}} \bigl[ v(X_{T_{\bar{\alpha}} \wedge
T_{\bar{\beta}}}) - A \varphi (X_{T_{\bar{\alpha}} \wedge
T_{\bar{\beta}}}) - B \psi (X_{T_{\bar{\alpha}} \wedge
T_{\bar{\beta}}}) \bigr] \right] \nonumber \\
& \hspace{50mm} = v(\bar{x}) - A \varphi (\bar{x}) - B \psi
(\bar{x}) + \e_{\bar{x}} \left[ \int _0^{T_{\bar{\alpha}} \wedge
T_{\bar{\beta}}} e^{-\Lambda_u} \, dA_u^{\Lop v} \right]
\label{vLP0}
\end{align}
for all points $\bar{\alpha} < \bar{x} < \bar{\beta}$ in
$\intI$ and all constants $A,B \in \R$.
Also, we fix any point $x \in \intI$ and we consider any
monotone sequences $(\alpha_n)$, $(\beta_n)$ in $\intI$
such that
\ben
\alpha_n < x < \beta_n \, \text{ for all } n \geq 1 \quad
\text{and} \quad \lim _{n \rightarrow \infty} \alpha_n =
\lim _{n \rightarrow \infty} \beta_n =x . \label{vLP1}
\een
If we define
\be
A_n = \frac{v(\beta_n) \psi (\alpha_n) - v(\alpha_n) \psi
(\beta_n)} {\varphi (\beta_n) \psi (\alpha_n) - \varphi
(\alpha_n) \psi (\beta_n)} \quad \text{and} \quad
B_n = \frac{\varphi (\beta_n) v(\alpha_n) - \varphi
(\alpha_n) v(\beta_n)} {\varphi (\beta_n) \psi (\alpha_n)
- \varphi (\alpha_n) \psi (\beta_n)} ,
\ee
then we can check that
\be
A_n \varphi (\alpha_n) + B_n \psi (\alpha_n) = v(\alpha_n)
\quad \text{and} \quad
A_n \varphi (\beta_n) + B_n \psi (\beta_n) = v(\beta_n) ,
\ee
and observe that the identity
\be
0 = v(x) - A_n \varphi (x) - B_n \psi (x) + \e_x \left[
\int _0^{T_{\alpha_n} \wedge T_{\beta_n}} e^{-\Lambda_u}
\, dA_u^{\Lop v} \right]
\ee
follows immediately from (\ref{vLP0}) for $\bar{\alpha}
= \alpha_n$, $\bar{x} = x$ and $\bar{\beta} = \beta_n$.
Since $- A^{\Lop v} = A^{-\Lop v}$ is a continuous increasing
process, this identity, (\ref{vLP1}) and the dominated
convergence theorem imply that
\ben
v(x) \geq A_n \varphi (x) + B_n \psi (x) \quad \text{and} \quad
v(x) = \lim _{n \rightarrow \infty} \bigl[ A_n \varphi (x) +
B_n \psi (x) \bigr] . \label{vLP4}
\een
Also, given any $y \in \mbox{} ]\beta_n, \beta[$, we can
see that (\ref{vLP0}) with $\bar{\alpha} = \alpha_n$,
$\bar{x} = \beta_n$ and $\bar{\beta} = y$ yields
\be
\bigl[ v(y) - A_n \varphi (y) - B_n \psi (y) \bigr] \e_{\beta_n}
\left[ e^{-\Lambda _{T_y}} {\bf 1} _{\{ T_y < T_{\alpha_n} \}}
\right] = \e_{\beta_n} \left[ \int _0^{T_{\alpha_n} \wedge T_y}
e^{-\Lambda_u} \, dA_u^{\Lop v} \right] ,
\ee
which implies that
\ben
A_n \varphi (y) + B_n \psi (y) \geq v(y) \quad \text{for all }
y \in \mbox{} ]\beta_n, \beta[ . \label{vLP5}
\een
Similarly,
\ben
A_n \varphi (y) + B_n \psi (y) \geq v(y) \quad \text{for all }
y \in \mbox{} ]\alpha, \alpha_n[ . \label{vLP6}
\een
Combining these results with (\ref{phi-psi->0}), 
we can see that
\be
A_n \geq \lim _{y \in \intI , \, y \downarrow \alpha}
\frac{v(y)}{\varphi (y)} \geq 0 \quad \text{and}
\quad B_n \geq \lim _{y \in \intI , \, y \uparrow \beta}
\frac{v(y)}{\psi (y)} \geq 0 \quad \text{for all } n
\geq 1 .
\ee
If we consider any sequence $(n_\ell)$ such that
$\lim _{\ell \rightarrow \infty} A_{n_\ell}$ exists,
then the positivity of the constants $A_n$, $B_n$ and
(\ref{vLP4}) imply that $\lim _{\ell \rightarrow \infty}
B_{n_\ell}$ also exists and that both limits are
positive and finite.
In particular, (\ref{vLP4}) and (\ref{vLP5})--(\ref{vLP6})
imply that
\begin{align}
v(x) & = \lim _{\ell \rightarrow \infty} A_{n_\ell} \varphi
(x) + \lim _{\ell \rightarrow \infty} B_{n_\ell} \psi (x)
\nonumber \\
\intertext{and}
v(y) & \leq \lim _{\ell \rightarrow \infty} A_{n_\ell}
\varphi (y) + \lim _{\ell \rightarrow \infty} B_{n_\ell}
\psi (y) \quad \text{for all } y \in \intI \setminus \{ x \}
. \nonumber
\end{align}
It follows that $v(x)$ is greater than or equal to the
right-hand side of (\ref{vLP}).

The existence of constants $\tilde{A}$, $\tilde{B}$
such that the identity in (\ref{vLP+}) holds true follows
from the fact that the measure $\Lop v$ does not charge
the interval $]c,d[$.
If $[d,\beta[$ is not empty, then, given any $\bar{\alpha}
< \bar{x}$ in $]c,d[$ and $y \in [d, \beta[$, we can see
that (\ref{vLP0}) with $\bar{\beta} = y$ yields
\be
\bigl[ v(y) - \tilde{A} \varphi (y) - \tilde{B} \psi (y) \bigr]
\e_{\bar{x}} \left[ e^{-\Lambda _{T_y}} {\bf 1}
_{\{ T_y < T_{\bar{\alpha}} \}} \right] = \e_{\bar{x}} \left[
\int _0^{T_{\bar{\alpha}} \wedge T_y} e^{-\Lambda_u}
\, dA_u^{\Lop v} \right] .
\ee
It follows that $\tilde{A} \varphi (y) + \tilde{B} \psi (y)
\geq v(y) \geq f(y)$ because $- A^{\Lop v} = A^{-\Lop v}$
is a continuous increasing process.
We can show that $\tilde{A} \varphi (y) + \tilde{B} \psi (y)
\geq f(y)$ for all $y \in \mbox{} ]\alpha, c]$, if $]\alpha, c]
\neq \emptyset$, similarly, and the inequality in
(\ref{vLP+}) has been established.
\mbox{}\hfill$\Box$

\section{Ramifications including a generalisation of the
``principle of smooth fit''}
\label{sec:sf}

Throughout the section, we assume that (\ref{gr-cond}) is true, so
that the value function is real-valued, and that $f = \overline{f}$.
We can express the so-called waiting region $\mathcal W$ as
a countable union of pairwise disjoint open intervals because
it is an open subset of $\intI$.
In particular, we write
\ben
{\mathcal W} = \left\{ x \in \CI \mid \ v(x) > f(x) \right\}
= \bigcup _{\ell =1}^\infty {\mathcal W}_\ell , \label{W-ints}
\een
where
\be
{\mathcal W}_\ell = \mbox{} ]c_\ell , d_\ell[ , \quad \text{for some }
c_\ell, d_\ell \in \CI \cup \{ \alpha, \beta \} \text{ such that }
c_\ell \leq d_\ell ,
\ee
and we adopt the usual convention that $]c,c[ \mbox{} =
\emptyset$ for $c \in \CI \cup \{ \alpha, \beta \}$.
Since the measure $\Lop v$ does not charge the waiting region
$\mathcal W$,
\ben
v(x) = A_\ell \varphi (x) + B_\ell \psi (x) \quad \text{for
all } x \in {\mathcal W} _\ell , \label{v-rams}
\een
for some constants $A_\ell$ and $B_\ell$.

Our first result in this section is concerned with a characterisation
of the value function if the problem data is such that
${\mathcal W} = \intI$.
Example~\ref{ex1} in Section~\ref{sec:ex} provides an illustration
of this case.

\begin{cor} \label{no-stopping}
Consider the optimal stopping problem formulated in
Section~\ref{sec:pr-form}, and suppose that (\ref{gr-cond}) is
true and $f = \overline{f}$.
If ${\mathcal W}_1 = \mbox{} ]\alpha, \beta[$ and ${\mathcal W}
_\ell = \emptyset$ for $\ell > 1$, then
\ben
A_1 = \limsup _{y \downarrow \alpha} \frac{f(y)}{\varphi (y)}
\quad \text{and} \quad
B_1 = \limsup _{y \uparrow \beta} \frac{f(y)}{\psi (y)} .
\label{A1B1}
\een
\end{cor}
{\bf Proof.}
The result follows immediately from the fact that $v(x) = A_1
\varphi (x) + B_1 \psi (x)$ for all $x \in \intI$, (\ref{phi-psi->0})
and (\ref{v-lims}).
\mbox{}\hfill$\Box$
\vspace{5mm}

We next study the special case that arises when a portion of
the general problem's value function has the features of the
value function of a perpetual American call option, which has
been extensively studied in the literature.

\begin{cor} \label{call}
Consider the optimal stopping problem formulated in
Section~\ref{sec:pr-form}, and suppose that (\ref{gr-cond}) is
true and $f = \overline{f}$.
If ${\mathcal W}_\ell = \mbox{} ]\alpha, d_\ell[$, for some
$\ell \geq 1$ and $d_\ell \in \intI$, then
\begin{gather}
A_\ell = \limsup _{y \downarrow \alpha} \frac{f(y)}{\varphi (y)}
, \qquad B_\ell = \frac{1}{\psi (d_\ell)} \bigl[ f (d_\ell) -
A_\ell \varphi (d_\ell) \bigr] \label{ABcall} \\
\intertext{and}
\frac{f(x)}{A_\ell \varphi (x) + B_\ell \psi (x)} \begin{cases} < 1 &
\text{for all } x \in \mbox{} ]\alpha, d_\ell[ , \\ = 1 & \text{for } x =
d_\ell , \\ \leq 1 & \text{for all } x > d_\ell . \end{cases}
\label{BL1}
\end{gather}
\end{cor}
{\bf Proof.}
The identities in (\ref{ABcall}) follow immediately from the
fact that $v(x)$ is given by (\ref{v-rams}) for all $x \in
{\mathcal W}_\ell = \mbox{} ]\alpha, d_\ell[$, (\ref{phi-psi->0})
and (\ref{v-lims}).
The first two inequalities in (\ref{BL1}) are trivial.
Given any $x \in \mbox{} ]\alpha, d_\ell[$, the fact that $v(x)$
is given by (\ref{v-rams}) and part~(III) of Theorem~\ref{thm:main2}
imply that
\be
A_\ell \varphi (x) + B_\ell \psi (x) \geq f(y) \quad \text{for all }
y \in \intI ,
\ee
and the last inequality in (\ref{BL1}) follows.
\mbox{}\hfill$\Box$
\vspace{5mm}

Using similar symmetric arguments, we can establish the
following result that arises in the context of a perpetual
American put option.

\begin{cor} \label{put}
Consider the optimal stopping problem formulated in
Section~\ref{sec:pr-form}, and suppose that (\ref{gr-cond})
is true and $f = \overline{f}$.
If ${\mathcal W}_\ell = \mbox{} ]c_\ell, \beta[$, for some
$\ell \geq 1$ and $c_\ell \in \intI$, then
\begin{gather}
A_\ell =  \frac{1}{\varphi (c_\ell)} \bigl[ f (c_\ell) - B_\ell \psi
(c_\ell) \bigr] , \qquad
B_\ell = \limsup _{y \uparrow \beta} \frac{f(y)}{\psi (y)} \\
\intertext{and}
\frac{f(x)}{A_\ell \varphi (x) + B_\ell \psi (x)} \begin{cases} \leq
1 & \text{for all } x < c_\ell , \\ = 1 & \text{for } x = c_\ell , \\ < 1
& \text{for all } x \in \mbox{} ]c_\ell, \beta[ . \end{cases}
\end{gather}
\end{cor}

The final result in this section focuses on a special case in
which a component of the waiting region $\mathcal W$
has compact closure in $\intI$, which is a case that can arise
in the context of the valuation of a perpetual American straddle
option.

\begin{cor} \label{straddle}
Consider the optimal stopping problem formulated in
Section~\ref{sec:pr-form}, and suppose that (\ref{gr-cond}) is
true and $f = \overline{f}$.
If ${\mathcal W}_\ell = \mbox{} ]c_\ell, d_\ell[$, for some
$\ell \geq 1$ and $c_\ell, d_\ell \in \intI$, then
\begin{gather}
A_\ell =  \frac{f(d_\ell) \psi (c_\ell) - f(c_\ell) \psi(d_\ell)} {\varphi
(d_\ell) \psi (c_\ell) - \varphi (c_\ell) \psi (d_\ell)} , \qquad
B_\ell = \frac{\varphi (d_\ell) f(c_\ell) - \varphi (c_\ell) f(d_\ell)}
{\varphi (d_\ell) \psi (c_\ell) - \varphi (c_\ell) \psi (d_\ell)}
\label{AB-craddle} \\
\intertext{and}
\frac{f(x)}{A_\ell \varphi (x) + B_\ell \psi (x)} \begin{cases} <
1 & \text{for all } x \in \mbox{} ]c_\ell, d_\ell[ , \\ = 1 & \text{for }
x = c_\ell \text{ and } x = d_\ell , \\ \leq 1 & \text{for all } x \leq
c_\ell \text{ and } x \geq d_\ell . \end{cases} \label{BL3}
\end{gather}
\end{cor}
{\bf Proof.}
The expressions in (\ref{AB-craddle}) follow immediately
from the continuity of the value function.
The first two inequalities in (\ref{BL3}) are a consequence
of the definition of the waiting region $\mathcal W$, while
the last one is an immediate consequence of part~(II) of
Theorem~\ref{thm:main2}.
\mbox{}\hfill$\Box$
\vspace{5mm}

Our final result is concerned with a generalisation of the
``principle of smooth fit''.

\begin{cor} \label{smooth-fit}
Consider the optimal stopping problem formulated in
Section~\ref{sec:pr-form}, and suppose that (\ref{gr-cond}) is
true and $f = \overline{f}$.
Also, consider any point $y \in \intI$ such that $y \notin
{\mathcal W}$.
If $f$ admits right and left-hand derivatives at $y$, then
\ben
f_+' (y) \leq v_+' (y) \leq v_-' (y) \leq f_-' (y) .
\een
\end{cor}
{\bf Proof.}
The inequality $v_+' (y) \leq v_-' (y)$ is an immediate
consequence of the fact that $\Lop v\leq 0$.
The inequalities $f_+' (y) \leq v_+' (y)$ and $v_-' (y) \leq f_-' (y)$
follow from the fact that $v-f$ has a local minimum at $y$.
\mbox{}\hfill$\Box$

\section{Examples}
\label{sec:ex}

We assume that an appropriate weak solution $\str_x$ to
(\ref{SDE}) has been associated with each initial condition
$x \in \intI$ in all of the examples that we discuss in this
section.
The following example shows that an optimal stopping time
may not exist if (\ref{t*-cond121}) is not satisfied.
In this example, the stopping region $\CI \setminus
{\mathcal W}$ is empty.

\begin{exa} \label{ex1} {\rm
Suppose that $\CI = \mbox{} ]0,\infty[$ and $X$ is a geometric
Brownian motion, so that
\be
dX_t = bX_t \, dt + \sigma X_t \, dW_t ,
\ee
for some constants $b$ and $\sigma$.
Also, suppose that $r$ is a constant.
In this case, it is well-known that
\be
\varphi (x) = x^m \quad \text{and} \quad \psi (x) = x^n ,
\ee
where $m<0<n$ are the solutions to the quadratic equation
\ben
\half \sigma ^2 k^2 + \left( b - \half \sigma ^2 \right) k  - r = 0
. \label{mn}
\een
In this context, if the reward function $f$ is given by
\be
f(x) = \begin{cases} \kappa (x^m -x) , & \text{if } x \in \mbox{}
]0,1] , \\ \lambda (x^n - x^{-1}) , & \text{if } x > 1 , \end{cases}
\ee
for some constants $\kappa, \lambda > 0$, then
\begin{align}
v(x) & = \kappa x^m + \lambda x^n \nonumber \\
& = \lim _{j \rightarrow \infty} \e_x \left[ e^{-r (T_{\alpha_j}
\wedge T_{\beta_j})} f(X_{T_{\alpha_j} \wedge T_{\beta_j}})
\right] \quad \text{for all } x>0 ,
\end{align}
where $(\alpha_j)$ and $(\beta_j)$ are any sequences
in $]0,\infty[$ such that
\ben
\alpha_j < x < \beta_j \text{ for all } j , \quad \lim
_{j \rightarrow \infty} \alpha_j = 0 \quad \text{and} \quad
\lim _{j \rightarrow \infty} \beta_j = \infty . \label{abj-seq}
\een
In particular, there exists no optimal stopping time.
} \mbox{}\hfill$\Box$ \end{exa}

The next example shows that an optimal stopping time
may not exist if (\ref{t*-cond121}) is not satisfied, while the
stopping region $\CI \setminus {\mathcal W}$ is not empty.

\begin{exa} \label{ex2} {\rm
In the context of the previous example, suppose that
the reward function $f$ is given by
\be
f(x) = \begin{cases} 0 , & \text{if } x \in \mbox{} ]0,1[ , \\
1 , & \text{if } x = 1, \\ x^n - x^{-1} , & \text{if } x > 1 .
\end{cases}
\ee
In view of straightforward considerations, we can see that
\be
v(x) = x^n \quad \text{for all } x > 0 .
\ee
In this case,
\be
\tau^\star \equiv \inf \left\{ t \geq 0 \mid \ v(X_t) = f(X_t)
\right\} = T_1 ,
\ee
i.e., $\tau^\star$ is the first hitting time of $\{ 1 \}$, and
\be
v(x) = \lim _{j \rightarrow \infty} \e_x \left[ e^{- r \left(T_1
\wedge T_{\beta_j} \right)} f(X_{T_1 \wedge T_{\beta_j}})
\right] > x^m = \e_x \left[ e^{- r \tau_\star} f(X_{\tau_\star})
\right] \quad \text{for all } x>1 ,
\ee
where $(\beta_j)$ is any sequence in $]x,\infty[$ such
that $\lim _{j \rightarrow \infty} \beta_j = \infty$.
} \mbox{}\hfill$\Box$ \end{exa}

The following example shows that an optimal stopping time
may not exist if (\ref{t*-cond122}) is not satisfied.
In this example, the stopping region $\CI \setminus
{\mathcal W}$ is empty.

\begin{exa} \label{ex2+} {\rm
Suppose that $\CI = \R_+$, $X$ is a standard one-dimensional
Brownian motion starting from $x>0$ and absorbed at 0 and
$r$ is a constant.
In this case, we can see that
\be
\varphi (x) = e^{-\sqrt{2r} x} \quad \text{and} \quad \psi (x) =
e^{\sqrt{2r} x} - e^{-\sqrt{2r} x} .
\ee
If the reward function $f$ is given by
\be
f(x) = \begin{cases} 0, & \text{if } x=0 , \\ e^{-2 \sqrt{2r} x} ,
& \text{if } x>0 , \end{cases}
\ee
then we can see that
\be
v(x) = \begin{cases} 0 , & \text{if } x=0 , \\ e^{-\sqrt{2r} x} ,
& \text{if } x>0 . \end{cases}
\ee
In particular,
\be
v(x) = \lim _{j \rightarrow \infty} \e_x \left[ e^{-r (T_{\alpha_j}
\wedge T_{\beta_j})} f(X_{T_{\alpha_j} \wedge T_{\beta_j}})
\right] \quad \text{for all } x>0 ,
\ee
where $(\alpha_j)$, $(\beta_j)$ are any sequences in
$]0,\infty[$ satisfying (\ref{abj-seq}), and there exists no optimal
stopping time.
} \mbox{}\hfill$\Box$ \end{exa}

The following example shows that an optimal stopping time
may not exist if $f \neq \overline{f}$.
In particular, the first hitting time $\tau^\star$ of the stopping
region $\CI \setminus {\mathcal W}$ may not be optimal.

\begin{exa} \label{ex3} {\rm
Suppose that $X$ is a standard Brownian motion, namely,
$\CI = \R$ and $dX_t = dW_t$, and that $r = \frac{1}{2}$.
In this context, it is straightforward to verify that
\be
\varphi (x) = e^{-x} \quad \text{and} \quad \psi (x) = e^x .
\ee
Also, consider the reward function
\be
f(x) = \begin{cases} 0 , & \text{if } x \leq 0 , \\ 1 , & \text{if }
x \in \mbox{} ]0,1] , \\ 2 , & \text{if } x > 1 , \end{cases}
\ee
which is not upper semicontinuous.
In this case, we can see that
\be
v(x) = \begin{cases} e^x , & \text{if } x \leq 0 , \\
\frac{e-2}{e-e^{-1}} e^{-x} + \frac{2-e^{-1}}{e-e^{-1}}
e^x , & \text{if } x \in \mbox{} ]0,1] , \\ 2 , & \text{if } x > 1 .
\end{cases}
\ee
Given an initial condition $x<1$ and an associated solution
$\str_x$ to (\ref{SDE}), we note that
\be
\tau^\star \equiv \inf \bigl\{ t \geq 0 \mid \ v(X_t) = f(X_t) \bigr\}
= \inf \bigl\{ t \geq 0 \mid \ X_t > 1 \bigr\}
\ee
defines an $(\F_t)$-stopping time because we have
assumed that the filtration in $\str_x$ satisfies the
usual conditions.
However, $X_{\tau^\star} = 1$, $\p_x$-a.s., and
\be
\e_x \left[ e^{-r \tau^\star} f(X_{\tau^\star}) \right] =
e^{x-1} < v(x) \quad \text{for all } x< 1 .
\ee
In view of these considerations, we can see that there is
no optimal stopping time for initial conditions $x<1$.
} \mbox{}\hfill$\Box$ \end{exa}

The final example that we consider illustrates that a
characterisation such as the one provided by (\ref{BL1})
in Corollary~\ref{call} has a local rather than global
character.

\begin{exa} \label{ex5} {\rm
In the context of the previous example, we consider
the reward function
\be
f(x) = \begin{cases} e^{2x} , & \text{if } x<0 , \\ 1 , & \text{if }
x \in [0,1] , \\ 1+ (x-1)^2 , & \text{if } x>1 , \end{cases}
\ee
and we note that the calculation
\be
\frac{d}{dx} \frac{f(x)}{\psi (x)} = \begin{cases} e^{x} , &
\text{if } x<0 , \\ - e^{-x} , & \text{if } x \in [0,1] , \\ - (x-2)^2
e^{-x} , & \text{if } x>1 \end{cases}
\ee
implies that the function $f/\psi$ is strictly increasing in
$]-\infty, 0[$ and strictly decreasing in $]0,\infty[$.
A first consideration of the associated optimal stopping
problem suggests that the value function $v$ could
identify with the function $u$ given by
\be
u(x) = \begin{cases} e^x , & \text{if } x<0 , \\ 1 , & \text{if }
x \in [0,1] , \\ 1+ (x-1)^2 , & \text{if } x>1 . \end{cases}
\ee
In particular, we can check that
\be
\frac{u(x)}{u(y)} \geq \min \left\{ \frac{\varphi (x)}{\varphi (y)}
, \frac{\psi (x)}{\psi (y)} \right\} \quad \text{for all } x,y \in \R .
\ee
However, the function $u$ is not excessive because
\be
\Lop u(dx) \equiv \half u''(dx) - \half u(x) \, dx = - \delta _0
(dx) - \half \left( {\bf 1} _{[0,1]} (x) + x (x-2) {\bf 1} _{[1,\infty[}
(x) \right) dx , 
\ee
where $\delta_0$ is the Dirac probability measure that assigns
mass 1 on the set $\{ 0 \}$, which implies that
\be
\Lop u \bigl( [c,d] \bigr) > 0 \quad \text{for all } 1 \leq c < d
\leq 2 ,
\ee
and suggests that $[1,2]$ should be a subset of the
waiting region $\mathcal W$.
In this example, the value function $v$ is given by
\ben
v(x) = \begin{cases} e^x , & \text{if } x<0 , \\ 1 , & \text{if }
x \in [0,a_l] , \\ \half e^{a_l -  x} + \half e^{-a_l + x} , &
\text{if } x \in \mbox{} ]a_l, a_r[ , \\ 1+ (x-1)^2 , & \text{if }
x > a_r , \end{cases} \label{v-exa}
\een
where
\be
a_l = 1 + \sqrt{2} + 2 \ln \left( \sqrt{2} - 1 \right) \simeq
0.651 \quad \text{and} \quad a_r = 1 + \sqrt{2} \simeq
2.414 .
\ee
These values for the boundary points $a_l$, $a_r$
arise by the requirements that $a_l \in \mbox{} ]0,1[$,
$a_r > 2$ and $v$ should be $C^1$ along $a_l$, $a_r$
(see Corollary~\ref{smooth-fit}), which are associated with
the system of equations
\begin{gather}
a_l = a_r + 2 \ln \left( a_r - 2 \right) , \nonumber \\
a_r^4 - 4a_r^3 + 4a_r^2 -1 \equiv (a_r - 1)^2 \left(a_r -1
- \sqrt{2} \right) \left( a_r -1 + \sqrt{2} \right) = 0 . \nonumber
\end{gather}
In particular, we can check that the function given by the
right-hand side of (\ref{v-exa}) satisfies all of the requirements
of the verification Theorem~\ref{thm:main2}.(II) and therefore
identifies with the value function $v$.
} \mbox{}\hfill$\Box$ \end{exa}

\section*{Appendix: pasting weak solutions of SDEs}

The next result is concerned with aggregating two filtrations, one of
which ``switches on'' at a stopping time of the other one.

\begin{thm} \label{prop:augm-filtr}
Consider a measurable space $(\Om, \F)$ and two filtrations $(\Hc_t)$,
$(\G_t)$ such that $\Hc_\infty \cup \G_\infty \subseteq \F$.
Also, suppose that $(\G_t)$ is right-continuous and let $\tau$ be an
$(\Hc_t)$-stopping time.
If we define
\begin{align}
\F_t = \bigl\{ A \in \Hc_\infty \vee \G_\infty \mid \ & A \cap \{
t < \tau \} \in \Hc_t \vee \G_0 \text{ and } \nonumber \\
& A \cap \{ s \leq \tau \} \in \Hc_t \vee \G_{t-s} \text { for all }
s \in [0,t] \bigr\} , \label{hatFdef}
\end{align}
then $(\F_t)$ is a filtration such that
\begin{gather}
\F_{\tau + t} = \Hc_{\tau + t} \vee \G_t \quad \text{for all } t
\geq 0 \label{hatFt+tau} \\
\intertext{and}
\F_{t \wedge \tau} = \Hc_{t \wedge \tau} \vee \G_0 \quad
\text{for all } t \geq 0 . \label{hatFttau}
\end{gather}
\end{thm}
{\bf Proof.}
First, it is straightforward to verify that
\ben
\Hc_t \subseteq \F_t \quad \text{and} \quad \G_0 \subseteq
\F_t \quad \text{for all } t \geq 0 . \label{PP0}
\een
To prove that $(\F_t)$ is indeed a filtration, we consider any times
$u < t$ and any event $A \in \F_u$.
Using the definition of $\F_u$, we can see that
\begin{gather}
A \cap \{ t < \tau \} = A \cap \{ u < \tau \} \cap \{ t < \tau \} \in
\Hc_u \vee \G_0 \vee \Hc_t \subseteq \Hc_t \vee \G_0 ,
\nonumber \\
A \cap \{ s \leq \tau \} \in \Hc_u \vee \G_{u-s} \subseteq \Hc_t
\vee \G_{t-s} \quad \text{for all } s \in [0,u] , \nonumber
\intertext{and}
A \cap \{ s \leq \tau \} = A \cap \{ u \leq \tau \} \cap \{ s \leq \tau \}
\in \Hc_u \vee \G_0 \vee \Hc_s \subseteq \Hc_t \vee \G_{t-s}
\quad \text{for all } s \in [u,t] . \nonumber
\end{gather}
It follows that $A \in \F_t$.

To establish (\ref{hatFt+tau}), we first show that $\G_t \subseteq
\F_{\tau + t}$, which amounts to proving that, given any
$t \geq 0$ and $A \in \G_t$,
\ben
A \cap \{ \tau + t \leq u \} = A \cap \{ \tau \leq u-t \} \in \F_u
\quad \text{for all } u \geq 0 . \label{GcF}
\een
To this end, we note that
\be
A \cap \{ \tau \leq u-t \} \cap \{ u < \tau \} = \emptyset
\in \Hc_u  \vee \G_0 \quad \text{for all } u \geq 0 .
\ee
Also, given any $s, u \geq 0$ such that $s \in \mbox{} ]u-t, u]$,
\be
A \cap \{ \tau \leq u-t \} \cap \{ s \leq \tau \} = \emptyset
\in \Hc_u  \vee \G_{u-s} ,
\ee
while, given any $s, u \geq 0$ such that $s \in [0, u-t]$,
\be
A \cap \{ \tau \leq u-t \} \cap \{ s \leq \tau \} = A \cap \{ s \leq \tau
\leq u-t \} \in \Hc_{u-t}  \vee \G_t
\subseteq \Hc_u  \vee \G_{u-s} .
\ee
These observations and the definition (\ref{hatFdef}) of $(\F_t)$
imply that (\ref{GcF}) holds true and $\G_t \subseteq
\F_{\tau + t}$.
Combining this result with the fact that $\Hc_{\tau + t}
\subseteq \F_{\tau + t}$, which follows from (\ref{PP0}), we can
see that $\Hc_{\tau + t} \vee \G_t \subseteq \F_{\tau + t}$.

To prove that $\F_{\tau + t} \subseteq \Hc_{\tau + t} \vee \G_t$,
we consider any $A \in \F_{\tau + t}$ and we show that
\ben
A \cap \{ \tau +t \leq u \} \in \Hc_u \vee \G_t \quad \text{for
all } u \geq 0 . \label{PP-30}
\een
Since $A \cap \{ \tau + t \leq \bar{u} \} \in \F_{\bar{u}}$ for
all $\bar{u} \geq 0$, $A \cap \{ \tau \leq \bar{u} \} \in
\F_{\bar{u}+t}$ for all $\bar{u} \geq 0$.
Combining this observation with the definition (\ref{hatFdef}) of
$(\F_t)$ we can see that
\ben
A \cap \{ \tau \leq \bar{u} \} \cap \{ s \leq \tau \} \in
\Hc_{\bar{u}+t} \vee \G_{\bar{u}+t-s} \quad \text{for all } \bar{u}
\geq 0 \text{ and } s \leq \bar{u}+t . \label{PP-3}
\een
In particular,
\be
A \cap \{ \bar{u} - \varepsilon \leq \tau \leq \bar{u} \} \in
\Hc_{\bar{u}+t} \vee \G_{t+\varepsilon} \quad \text{for all } \bar{u}
> 0 \text{ and } \varepsilon \in [0, \bar{u}] .
\ee
In view of this result, we can see that, given any $u>t$,
\be
A \cap \{ \tau + t \leq u \} = \bigcup _{j=0}^{n-1} A \cap \left\{
\frac{j(u-t)}{n} \leq \tau \leq \frac{(j+1)(u-t)}{n} \right\} \in
\Hc_u \vee \G_{t+\frac{(u-t)}{n}} .
\ee
It follows that
\be
A \cap \{ \tau + t \leq u \} \in \bigcap _{n=1}^\infty \Hc_u
\vee \G_{t+\frac{(u-t)}{n}} = \Hc_u \vee \G_t \quad
\text{for all } u>t ,
\ee
the equality being true thanks to the right continuity of $(\G_t)$.
Combining this result with the fact that
\be
A \cap \{ \tau + t \leq t \} \in \Hc_t \vee \G_t ,
\ee
which follows from (\ref{PP-3}) for $\bar{u} = s = 0$, we obtain
(\ref{PP-30}).

To prove (\ref{hatFttau}), we first note that (\ref{PP0}) implies that
$\Hc_{t \wedge \tau} \vee \G_0 \subseteq \F_{t \wedge \tau}$.
To establish the reverse inclusion, we consider any $A \in \F
_{t \wedge \tau}$ and we show that
\be
A \cap \{ t \wedge \tau \leq u \} \in \Hc_u \vee \G_0 \quad \text{for
all } u \geq 0 .
\ee
Since $A \cap \{ t \wedge \tau \leq \bar{u} \} \in \F_{\bar{u}}$
for all $\bar{u} \geq 0$, the definition (\ref{hatFdef}) of
$(\F_t)$ implies that
\ben
A \cap \{ t \wedge \tau \leq \bar{u} \} \cap \{ s \leq \tau \} \in
\Hc_{\bar{u}} \vee \G_{\bar{u}-s} \quad \text{for all } \bar{u} \geq 0
\text{ and } s \in [0,\bar{u}] . \label{PP2}
\een
For $\bar{u} = s = 0$, this implies immediately that
\ben
A \cap \{ t \wedge \tau \leq 0 \} \in \Hc_0 \vee \G_0 .
\label{PP3}
\een
Also, it implies that
\be
A \cap \{ \bar{u} - \varepsilon \leq \tau \leq \bar{u} \} \in
\Hc_{\bar{u}} \vee \G_\varepsilon \quad \text{for all } \bar{u}
\in [0, t[ \mbox{} \text{ and } \varepsilon \in [0,\bar{u}] .
\ee
In view of this observation, we can see that
\be
A \cap \{ t \wedge \tau \leq u \} = \bigcup _{j=0}^{n-1} A \cap
\left\{ \frac{ju}{n} \leq \tau \leq \frac{(j+1)u}{n} \right\} \in
\Hc_u \vee \G_{\frac{u}{n}} \quad \text{for all } u \in \mbox{}
]0,t[ .
\ee
It follows that
\be
A \cap \{ t \wedge \tau \leq u \} \equiv A \cap \{ \tau \leq u \}
\in \bigcap _{n=1}^\infty \Hc_u \vee \G_{\frac{u}{n}} = \Hc_u
\vee \G_0 \quad \text{for all } u \in \mbox{} ]0,t[
\ee
because $(\G_t)$ is right-continuous.
In particular, this implies that
\be
A \cap \{ \tau < t \} = \bigcup _{n=1}^\infty A \cap \left\{ \tau \leq
\frac{n}{n+1} t \right\} \in \Hc_t \vee \G_0 .
\ee
Combining this observation with the fact that
\be
A \cap \{ t \leq \tau \} \in \Hc_t \vee \G_0 ,
\ee
which follows from (\ref{PP2}) for $\bar{u} = s = t$, we can see that
\ben
A \cap \{ t \wedge \tau \leq u \} = A = A \cap \{ \tau < t \} \cup A
\cap \{ t \leq \tau \} \in \Hc_t \vee \G_0 \subseteq \Hc_u \vee \G_0
\quad \text{for all } u \geq t . \label{PP4}
\een
From (\ref{PP3})--(\ref{PP4}), it follows that $A \in \Hc_{t \wedge
\tau} \vee \G_0$.
\mbox{}\hfill$\Box$
\vspace{5mm}

The following result is concerned with the pasting of two
stopping strategies, in particular, two weak solutions to (\ref{SDE}),
at an appropriate stopping time.

\begin{thm} \label{prop:pasting}
Consider initial conditions $x_0, x_1 \in \intI$ and stopping strategies
\be
(\str_{x_0}^0, \tau^0) = \bigl( \bigl( \Om^0, \F^0, \F_t^0,
\p_{x_0}^0, W^0, X^0 \bigr) , \tau^0 \bigr) \text{ and }
(\str_{x_1}^1, \tau^1) = \bigl( \bigl( \Om^1, \F^1, \F_t^1,
\p_{x_1}^1, W^1, X^1 \bigr) , \tau^1 \bigr) .
\ee
Given any event $A \in \F_{\tau^0}^0$ such that $A \subseteq
\{ X_{\tau^0}^0 {\bf 1} _\{ \tau^0 < \infty \} = x_1 \}$, there exists a
stopping strategy
$(\str_{x_0}, \tau^{0,1}) = \bigl( \bigl( \Om, \F, \F_t, \p_{x_0},
W, X \bigr) , \tau^{0,1} \bigr) \in \T_{x_0}$ such that
\ben
J(\str_{x_0}, \tau^{0,1}) = J(\str_{x_0}^0, \tau_{A^c}^0) +
\e_{x_0}^0 \left[ e^{-\Lambda _{\tau^0} (X^0)} {\bf 1} _A \right]
J(\str_{x_1}, \tau^1) , \label{J-paste}
\een
where $\tau_{A^c}^0$ is the $(\F_t^0)$-stopping time defined by
$\tau_{A^c}^0 = \tau^0 {\bf 1} _{A^c} + \infty {\bf 1} _A$.
\end{thm}
{\bf Proof.}
Without loss of generality, we assume that $\{ X_{\tau^0}^0 =
x_1 \} \neq \emptyset$.
For $j=0,1$, we define on the product space $\bigl( \Om, \F,
\p_{x_0} \bigr) = \bigl( \Om^0 \times \Om^1 , \F^0 \otimes \F^1
, \p_{x_0}^0 \otimes \p_{x_1}^1 \bigr)$ the independent filtrations
$(\tilde{\F}_t^j)$ given by $\tilde{\F}_t^0 = \F_t^0 \otimes \left\{
\Om^1, \emptyset \right\}$ and $\tilde{\F}_t^1 = \left\{ \Om^0,
\emptyset \right\} \otimes \F_t^1$, the $(\tilde{\F}_t^j)$-stopping
times $\tilde{\tau}^j$ given by $\tilde{\tau}^j (\om^0, \om^1) =
\tau^j (\om^j)$, the $(\tilde{\F}_t^j)$-Brownian motions
$\tilde{W}^j$ given by $\tilde{W}_t^j (\om^0, \om^1) = W_t^j
(\om^j)$, and the continuous $(\tilde{\F}_t^j)$-adapted processes
$\tilde{X}^j$ given by $\tilde{X}_t^j (\om^0, \om^1) = X_t^j
(\om^j)$.
Also, we denote by $\tilde{T}_y^j$ the first hitting time of
$\{ y \}$ by $\tilde{X}^j$, for $y \in \CI$ and $j=0,1$.
In particular, we note that each of the collections $\bigl( \Om,
\F, \tilde{\F}^j, \p_{x_0}, \tilde{W}^j, \tilde{X}^j \bigr)$ is a weak
solution to the SDE (\ref{SDE}) with initial condition $x_j$.

We next consider the filtration $(\F_t)$ that is defined by
(\ref{hatFdef}) in Proposition~\ref{prop:augm-filtr} above with
$(\Hc_t) = (\tilde{\F}_t^0)$, $(\G_t) = (\tilde{\F}_t^1)$
and $\tau = \tilde{\tau}^0$, so that
\ben
\F_{\tilde{\tau}^0 + t} = \tilde{\F}_{\tilde{\tau}^0 + t}^0 \vee
\tilde{\F}_t^1 \quad \text{and} \quad \F_{t \wedge \tilde{\tau}^0}
= \tilde{\F}_{t \wedge \tilde{\tau}^0}^0 \vee \tilde{\F}_0^1
, \label{paste12-filtrations}
\een
and we define
\ben
\tilde{A} = A \times \Om^1 \quad \text{and} \quad \tilde{A}^c
= A^c \times \Om^1 . \label{A-tilde}
\een

The independence of $(\tilde{\F}_t^0)$, $(\tilde{\F}_t^1)$
and (\ref{hatFttau}) imply that the processes 
$\bigl(\tilde{W}_{t \wedge \tilde{\tau}^0}^0, \ t \geq 0 \bigr)$ and
$\bigl((\tilde{W}_{t \wedge \tilde{\tau}^0}^0)^2 - t \wedge
\tilde{\tau}^0 , \ t \geq 0 \bigr)$ are $(\F_t)$-martingales.
On the other hand, (\ref{paste12-filtrations}) and the fact that
$(\tilde{\F}_t^0)$, $(\tilde{\F}_t^1)$ are independent imply
that $\tilde{W}^1$ is an $(\F_{\tilde{\tau}^0+t})$-Brownian motion.
Since $(t-\tilde{\tau}^0)^+$ is an $(\F_{\tilde{\tau}^0+t})$-stopping
time for all $t \geq 0$ and $\tilde{\tau}^0 + (t-\tilde{\tau}^0)^+ =
\tilde{\tau}^0 \vee t$, the time-changed processes $\bigl(
\tilde{W}_{(t-\tilde{\tau}^0)^+}^1, \ t \geq 0 \bigr)$ and $\bigl(
(\tilde{W}_{(t-\tilde{\tau}^0)^+}^1)^2 - (t-\tilde{\tau}^0)^+ , \ t
\geq 0 \bigr)$ are $(\F_{\tilde{\tau}^0 \vee t})$-martingales, while
the $(\F_{\tilde{\tau}^0 \vee t})$-adapted process $\bigl(
\tilde{X}_{(t-\tilde{\tau}^0)^+}^1, \ t \geq 0 \bigr)$ satisfies
\begin{align}
\tilde{X} _{(t-\tilde{\tau}^0)^+}^1 & = x_1 + \int
_0^{(t-\tilde{\tau}^0)^+} b \bigl( \tilde{X} _s^1 \bigr) \, ds + \int
_0^{(t-\tilde{\tau}^0)^+} \sigma \bigl( \tilde{X} _s^1 \bigr) \,
d\tilde{W}_s^1 \nonumber \\
& = x_1 + \int _0^t b \bigl( \tilde{X} _{(s-\tilde{\tau}^0)^+}^1
\bigr) \, d{(s-\tilde{\tau}^0)^+} + \int _0^t \sigma \bigl( \tilde{X}
_{(s-\tilde{\tau}^0)^+}^1 \bigr) \,
d\tilde{W}_{(s-\tilde{\tau}^0)^+}^1 \nonumber \\
& = x_1 + \int _0^t {\bf 1} _{\{ \tilde{\tau}^0 \leq s \}} b \bigl(
\tilde{X} _{(s-\tilde{\tau}^0)^+}^1 \bigr) \, ds + \int _0^t {\bf 1}
_{\{ \tilde{\tau}^0 \leq s \}} \sigma \bigl( \tilde{X}
_{(s-\tilde{\tau}^0)^+}^1 \bigr) \, d\tilde{W}_s^1 . \label{pasteSDE1}
\end{align}
In fact, all of these processes are $(\F_t)$-adapted.
To see this, we consider, e.g., the process $\bigl( \tilde{X}
_{(t-\tilde{\tau}^0)^+}^1 , \ t \geq 0)$, and we note that $\bigl\{
\tilde{X} _{(t-\tilde{\tau}^0)^+}^1 \in C \} \in \F
_{\tilde{\tau}^0 \vee t}$ implies that
\be
\bigl\{ \tilde{X} _{(t-\tilde{\tau}^0)^+}^1 \in C \} \cap \{
\tilde{\tau}^0 \vee t \leq u \} \in \F_u \quad \text{for all } u \geq
0 .
\ee
Therefore,
\be
\bigl\{ \tilde{X} _{(t-\tilde{\tau}^0)^+}^1 \in C \} \cap \{
\tilde{\tau}^0 < t \} = \bigcup _{n=1}^\infty \bigl\{ \tilde{X}
_{(t-\tilde{\tau}^0)^+}^1 \in C \} \cap \left\{ \tilde{\tau}^0 \vee t
\leq \frac{nt}{n+1} \right\} \in \F_t .
\ee
It follows that
\be
\bigl\{ \tilde{X} _{(t-\tilde{\tau}^0)^+}^1 \in C \} = \bigl\{
\tilde{X} _0^1 \in C \} \cap \{ t \leq \tilde{\tau}^0 \} \cup \bigl\{
\tilde{X} _{(t-\tilde{\tau}^0)^+}^1 \in C \} \cap \{ \tilde{\tau}^0 <
t \} \in \F_t ,
\ee
because $\tilde{X} _0^1 = x_1$ is a constant, which establishes the
claim.
Furthermore, $\bigl(
\tilde{W}_{(t-\tilde{\tau}^0)^+}^1, \ t \geq 0 \bigr)$ and $\bigl(
(\tilde{W}_{(t-\tilde{\tau}^0)^+}^1)^2 - (t-\tilde{\tau}^0)^+ , \ t
\geq 0 \bigr)$ are in fact $(\F_t)$-martingales.
Indeed, given $s<t$, we can check, e.g., that
\begin{align}
\e_{x_0} & \left[ \tilde{W}_{(t-\tilde{\tau}^0)^+}^1 \mid \F_s \right]
\nonumber \\
& = \e_{x_0} \left[ \e_{x_0} \left[ \tilde{W}_{(t-\tilde{\tau}^0)^+}^1
\mid \F_{\tilde{\tau}^0 \vee s} \right] \mid \F_s \right]
= \e_{x_0} \left[ \tilde{W}_{(s-\tilde{\tau}^0)^+}^1 \mid \F_s
\right] = \tilde{W}_{(s-\tilde{\tau}^0)^+}^1 , \label{paste-mart1}
\end{align}
the last equality following because $\bigl( \tilde{W}
_{(t-\tilde{\tau}^0)^+}^1, \ t \geq 0 \bigr)$ is $(\F_t)$-adapted.
For future reference, we also note that
\begin{align}
\e_{x_0} & \left[ \tilde{W}_{t \wedge \tilde{\tau}^0}^0 \tilde{W}
_{(t-\tilde{\tau}^0)^+}^1 \mid \F_s \right] \nonumber \\
= \mbox{} & \e_{x_0} \left[ \e_{x_0} \left[ \tilde{W}_{\tilde{\tau}^0}^0
\tilde{W} _{(t-\tilde{\tau}^0)^+}^1 \mid \F_{\tilde{\tau}^0 \vee s}
\right] \mid \F_s \right] = \e_{x_0} \left[ \tilde{W}_{\tilde{\tau}^0}^0
\tilde{W} _{(s-\tilde{\tau}^0)^+}^1 \mid \F_s \right] = \tilde{W}
_{s \wedge \tilde{\tau}^0}^0 \tilde{W} _{(s-\tilde{\tau}^0)^+}^1
. \label{paste-mart2}
\end{align}

In view of (\ref{paste12-filtrations}) and (\ref{A-tilde}),
the process $Y$ defined by $Y_t = {\bf 1} _{\tilde{A}} {\bf 1}
_{\{ \tilde{\tau}^1 < t \}}$ is $(\F_{\tilde{\tau}^0+t})$-adapted.
Using arguments similar to the ones we developed above, we can
see that the time-changed process given by
\be
Y_{(t-\tilde{\tau}^0)^+} = {\bf 1} _{\tilde{A}} {\bf 1}
_{\{ \tilde{\tau}^0 + \tilde{\tau}^1 < t \}} , \quad t \geq 0 ,
\ee
is $(\F_t)$-adapted, which proves that the random variable
$(\tilde{\tau}^0 + \tilde{\tau}^1) {\bf 1} _{\tilde{A}} + \infty
{\bf 1} _{\tilde{A}^c}$ is an $(\F_t)$-stopping time.
It follows that the random variable
\ben
\tau^{0,1} = \min \left\{ \tilde{\tau}^0 {\bf 1} _{\tilde{A}^c} +
\infty {\bf 1} _{\tilde{A}} , \ (\tilde{\tau}^0 + \tilde{\tau}^1) {\bf 1}
_{\tilde{A}} + \infty {\bf 1} _{\tilde{A}^c} \right\} = \tilde{\tau}^0
{\bf 1} _{\tilde{A}^c} + (\tilde{\tau}^0 + \tilde{\tau}^1) {\bf 1}
_{\tilde{A}} \label{pasteSDEtau}
\een
is an $(\F_t)$-stopping time.

To proceed further, we define
\be
W_t = \tilde{W}_{t \wedge \tilde{\tau}^0}^0 + \tilde{W}
_{(t-\tilde{\tau}^0)^+}^1
\ee
and
\begin{align}
X_t & = \tilde{X} _{t \wedge \tilde{\tau}^0}^0 + \left( \tilde{X} _t^0
- \tilde{X} _{\tilde{\tau}^0}^0 \right) {\bf 1} _{\tilde{A}^c} {\bf 1}
_{\{ \tilde{\tau}^0 \leq t \}} + \left( \tilde{X}
_{(t-\tilde{\tau}^0)^+}^1 - x_1 \right) {\bf 1} _{\tilde{A}} {\bf 1}
_{\{ \tilde{\tau}^0 \leq t \}} \nonumber \\
& \equiv \tilde{X} _t^0 {\bf 1} _{\{ t < \tilde{\tau}^0 \}} +
\tilde{X} _t^0 {\bf 1} _{\tilde{A}^c} {\bf 1} _{\{ \tilde{\tau}^0
\leq t \}} + \tilde{X} _{(t-\tilde{\tau}^0)^+}^1 {\bf 1} _{\tilde{A}}
{\bf 1} _{\{ \tilde{\tau}^0 \leq t \}} , \label{pasteSDE2}
\end{align}
and we note that
\ben
X_{\tau^{0,1}} = \tilde{X} _{\tilde{\tau}^0}^0 {\bf 1} _{\tilde{A}^c}
+ \tilde{X} _{\tilde{\tau}^1}^1 {\bf 1} _{\tilde{A}} . \label{pasteSDE3}
\een
Given $y \in \CI$, if we denote by $T_y$ the first hitting time of
$\{ y \}$ by $X$, then
\begin{gather}
T_\alpha {\bf 1} _{\tilde{A}^c} = \tilde{T}_\alpha^0 {\bf 1} _{\tilde{A}^c}
, \quad T_\beta {\bf 1} _{\tilde{A}^c} = \tilde{T}_\beta^0 {\bf 1}
_{\tilde{A}^c} ,  \label{pasteTab1} \\
T_\alpha {\bf 1} _{\tilde{A}} = \left( \tilde{\tau}^0 + \tilde{T}_\alpha^1
\right) {\bf 1} _{\tilde{A}} \quad \text{and} \quad T_\beta {\bf 1}
_{\tilde{A}} = \left( \tilde{\tau}^0 + \tilde{T} _\beta^1 \right) {\bf 1}
_{\tilde{A}} \label{pasteTab2}
\end{gather}
because $\tilde{\tau}^0 {\bf 1} _{\tilde{A}} < \bigl( \tilde{T}_\alpha^0
+ \tilde{T} _\beta^1 \bigr) {\bf 1} _{\tilde{A}}$.
Since $W$ is the sum of two $(\F_t)$-martingales, it is an
$(\F_t)$-martingale.
Furthermore, since
\be
W_t^2 - t = \left\{  \left( \tilde{W}_{t\wedge \tilde{\tau}^0}^0 \right)
^2 - t \wedge \tilde{\tau}^0 \right\} + \left\{ \left( \tilde{W}
_{(t-\tilde{\tau}^0)^+}^1 \right) ^2 - (t-\tilde{\tau}^0)^+ \right\}
+ 2 \tilde{W}_{t\wedge \tilde{\tau}^0}^0 \tilde{W}
_{(t-\tilde{\tau}^0)^+}^1 ,
\ee
and the three processes identified on the right-hand side of this
identity are $(\F_t)$-martingales (see
(\ref{paste-mart1})--(\ref{paste-mart2})), the process $(W_t^2 - t)$
is an $(\F_t)$-martingale.
From L\'{e}vy's characterisation theorem, it follows that $W$ is an
$(\F_t)$-Brownian motion.
Also, combining (\ref{pasteSDE2}) with (\ref{pasteSDE1}) and
the fact that $\tilde{X}^0$ satisfies the SDE (\ref{SDE}), we can
see that
\begin{align}
X_t = \mbox{} & x_0 + \int _0^t {\bf 1} _{\{ s < \tilde{\tau}^0 \}} b
\bigl( \tilde{X} _s^0 \bigr) \, ds + \int _0^t {\bf 1} _{\{ s <
\tilde{\tau}^0 \}} \sigma \bigl( \tilde{X} _s^0 \bigr) \,
d\tilde{W}_s^0 \nonumber \\
& + {\bf 1} _{\tilde{A}^c} \int _0^t {\bf 1} _{\{ \tilde{\tau}^0 \leq
s \}} b \bigl( \tilde{X} _s^0 \bigr) \, ds + {\bf 1} _{\tilde{A}^c}
\int _0^t {\bf 1} _{\{ \tilde{\tau}^0 \leq s \}} \sigma \bigl(
\tilde{X} _s^0 \bigr) \, d\tilde{W}_s^0 \nonumber \\
& + {\bf 1} _{\tilde{A}} \int _0^t {\bf 1} _{\{ \tilde{\tau}^0 \leq
s \}} b \bigl( \tilde{X} _{(s-\tilde{\tau}^0)^+}^1 \bigr) \, ds +
{\bf 1} _{\tilde{A}} \int _0^t {\bf 1} _{\{ \tilde{\tau}^0 \leq s \}}
\sigma \bigl( \tilde{X} _{(s-\tilde{\tau}^0)^+}^1 \bigr) \,
d\tilde{W}_s^1 \nonumber \\
= \mbox{} & x_0 + \int _0^t b(X_s) \, ds + \int _0^t \sigma (X_s) \,
dW_s . \nonumber
\end{align}
This calculation and the preceding considerations show that
\be
(\str_{x_0}, \tau^{0,1}) = \bigl( \bigl( \Om, \F, \F_t, \p_{x_0},
W, X \bigr), \tau^{0,1} \bigr) \in \T_{x_0} .
\ee

To complete the proof, we use the definition (\ref{Lambda})
of $\Lambda$, (\ref{pasteSDEtau})--(\ref{pasteSDE2}) and
(\ref{pasteTab1})--(\ref{pasteTab2}) to calculate
\be
{\bf 1} _{\tilde{A}} \Lambda _{\tau^{0,1} \wedge T_\alpha \wedge
T_\beta} (X) = {\bf 1} _{\tilde{A}} \int _0^{\left( \tilde{\tau}^0
+ \tilde{\tau}^1 \right) \wedge \left( \tilde{\tau}^0 + \tilde{T}
_\alpha^1 \right) \wedge \left( \tilde{\tau}^0 + \tilde{T} _\beta^1
\right)} r(X_s) \, ds = {\bf 1} _{\tilde{A}} \left[ \Lambda
_{\tilde{\tau}^0} (\tilde{X}^0) + \Lambda _{\tilde{\tau}^1 \wedge
\tilde{T} _\alpha^1 \wedge \tilde{T} _\beta^1} (\tilde{X}^1) \right]
\ee
In view of this observation, (\ref{paste12-filtrations})--(\ref{A-tilde}),
(\ref{pasteSDEtau})--(\ref{pasteTab2}) and the independence
of $(\tilde{\F}_t^0)$, $(\tilde{\F}_t^1)$, we can see that
\begin{align}
J(\str_{x_0}, \tau^{0,1}) = \mbox{} & \e_{x_0} \left[ e^{-\Lambda
_{\tau^{0,1} \wedge T_\alpha \wedge T_\beta} (X)}
f(X_{\tau^{0,1} \wedge T_\alpha \wedge T_\beta}) {\bf 1}
_{\{ \tau^{0,1} < \infty \}} \right] \nonumber \\
= \mbox{} & \e_{x_0} \left[ e^{-\Lambda _{\tilde{\tau}^0 \wedge
\tilde{T} _\alpha^0 \wedge \tilde{T} _\beta^0} (\tilde{X}^0)}
f(\tilde{X}_{\tilde{\tau}^0 \wedge \tilde{T} _\alpha^0 \wedge
\tilde{T} _\beta^0}^0) {\bf 1} _{\tilde{A}^c} {\bf 1} _{\{
\tilde{\tau}^0 < \infty \}} \right] \nonumber \\
& + \e_{x_0} \left[ e^{-\Lambda _{\tilde{\tau}^0} (\tilde{X}^0)}
\e_{x_0} \left[ e^{-\Lambda _{\tilde{\tau}^1 \wedge \tilde{T}
_\alpha^1 \wedge \tilde{T} _\beta^1} (\tilde{X}^1)}
f(\tilde{X}_{\tilde{\tau}^1 \wedge \tilde{T} _\alpha^1 \wedge
\tilde{T} _\beta^1}^1) {\bf 1} _{\{ \tilde{\tau}^1 < \infty
\}} \mid \tilde{\F}_{\tilde{\tau}^0}^0 \vee \tilde{\F}_0^1 \right]
{\bf 1} _{\tilde{A}} \right] \nonumber \\
= \mbox{} & \e_{x_0}^0 \left[ e^{-\Lambda _{\tau^0 \wedge
T_\alpha^0 \wedge T_\beta^0} (X^0)} f(X_{\tau^0 \wedge
T_\alpha^0 \wedge T_\beta^0}^0) {\bf 1} _{A^c} {\bf 1}
_{\{ \tau^0 < \infty \}} \right] \nonumber \\
& + \e_{x_0} \left[ e^{-\Lambda _{\tilde{\tau}^0} (\tilde{X}^0)}
\e_{x_0} \left[ e^{-\Lambda _{\tilde{\tau}^1 \wedge \tilde{T}
_\alpha^1 \wedge \tilde{T} _\beta^1} (\tilde{X}^1)}
f(\tilde{X}_{\tilde{\tau}^1 \wedge
\tilde{T} _\alpha^1 \wedge \tilde{T} _\beta^1}^1) {\bf 1}
_{\{ \tilde{\tau}^1 < \infty \}} \mid \tilde{\F}_0^1 \right] {\bf 1}
_{\tilde{A}} \right] \nonumber \\
= \mbox{} & J(\str _{x_0}^0, \tau_{A^c}^0) + \e_{x_0} \left[
e^{-\Lambda _{\tilde{\tau}^0} (\tilde{X}^0)} {\bf 1} _{\tilde{A}}
\right] \e_{x_0} \left[ e^{-\Lambda _{\tilde{\tau}^1 \wedge
\tilde{T} _\alpha^1 \wedge \tilde{T} _\beta^1} (\tilde{X}^1)}
f(\tilde{X}_{\tilde{\tau}^1 \wedge \tilde{T} _\alpha^1 \wedge
\tilde{T} _\beta^1}^1) {\bf 1} _{\{ \tilde{\tau}^1 < \infty \}} \right]
\nonumber \\
= \mbox{} & J(\str _{x_0}^0, \tau_{A^c}^0) + \e_{x_0}^0 \left[
e^{-\Lambda _{\tau^0} (X^0)} {\bf 1} _A \right]
\e_{x_1}^1 \left[ e^{-\Lambda _{\tau^1 \wedge T_\alpha^1
\wedge T_\beta^1} (X^1)} f(X_{\tau^1 \wedge T_\alpha^1
\wedge T_\beta^1}^1) {\bf 1} _{\{ \tau^1 < \infty \}} \right] ,
\nonumber
\end{align}
and (\ref{J-paste}) follows.
\mbox{}\hfill$\Box$
\vspace{5mm}

Iterating the construction above, we obtain the following result.

\begin{cor} \label{cor:paste}
Fix an initial condition $x \in \intI$ and any distinct points
$a_1, \ldots , a_n \in \intI$.
Given stopping strategies
\be
(\str_x^0, \tau^0) = \bigl( \bigl( \Om^0, \F^0, \F_t^0, \p_x^0, W^0,
X^0 \bigr),\tau^0 \bigr) \quad \text{and} \quad
(\str_{a_i}^i, \tau^i) = \bigl( \bigl( \Om^i, \F^i, \F_t^i,
\p_{a_i}^i, W^i, X^i \bigr), \tau^i \bigr) ,
\ee
$i = 1, \ldots, n$, define $A = \left\{ X_{\tau^0}^0 {\bf 1} _{\{ \tau^0 <
\infty \}} \in \{ a_1, \ldots , a_n \} \right\} \in \F_{\tau^0}^0$.
Then, there exists a stopping strategy $(\str_x, \tau) \in \T_x$
such that
\be
J(\str_x, \tau) = J(\str_x^0, \tau_{A^c}^0) + \sum _{i=1}^n
\e_x^0 \left[ e^{-\Lambda_{\tau^0} (X^0)} {\bf 1} _{\{ X_{\tau^0}^0
= a_i \}} \right] J(\str_{a_i}^i, \tau^i) ,
\ee
where $\tau_{A^c}^0$ is the $(\F_t^0)$-stopping time defined by
$\tau_{A^c}^0 = \tau^0 {\bf 1} _{A^c} + \infty {\bf 1} _A$.
\end{cor}


\end{document}